\documentclass[11pt, leqno]{article}

\usepackage{amssymb}

\newtheorem{thm}{Theorem}[section]
\newtheorem{prop}[thm]{Proposition}
\newtheorem{lemma}[thm]{Lemma}

\newtheorem{Definition}[thm]{Definition}

\newtheorem{Remark}[thm]{Remark}

\newcounter{ex}[section]

\newcommand{\G}{{\bf G}}

\newcommand{\no}{\noindent}

\newcommand{\ov}{\overline}
\newcommand{\GL}{{\rm GL}}
\newcommand{\SL}{{\rm SL}}
\newcommand{\Sp}{{\rm Sp}}

\newcommand{\C}{{\bf C}}
\newcommand{\R}{{\bf R}}
\newcommand{\Q}{{\bf Q}}

\newcommand{\Ff}{{\bf F}}
\newcommand{\Gg}{{\cal G}}

\newcommand{\Gm}{{\bf G}_m}
\newcommand{\Z}{{\bf Z}}
\newcommand{\F}{{\cal F}}
\newcommand{\ti}{\tilde}
\newcommand{\Spec}{{\rm Spec }\, }
 \renewcommand{\O}{{\cal O}}

\newcommand{\M}{{\cal M}}
\renewcommand{\L}{{\cal L}}

\newcommand{\Res}{{\rm Res}}

\newcommand{\Tr}{{\rm Tr}}

\newcommand{\N}{{\cal N}}

\newcommand{\grass}{{\bf {Gr}}}

\newcommand{\flag}{{\bf {Fl}}}
\newcommand{\sflag}{{\bf {SFl}}}

\def\thfill{\null\nobreak\hfill}

\def\endproof{\thfill\vbox{\hrule
  \hbox{\vrule\hbox to 5pt{\vbox to 5pt{\vfil}\hfil}\vrule}\hrule}}

\newenvironment{eq}{\addtocounter{thm}{0}\begin{equation} }{\end{equation}}

\begin{document}
\title{Local models in the ramified case\\
II. Splitting models}
\author{G. Pappas\footnote{Supported in part by NSF Grant \# DMS02-01140 and by a Sloan Research Fellowship.}
\  and \ M. Rapoport}
\date{\ }
\maketitle
\begin{abstract}We study the reduction of certain PEL Shimura varieties
with parahoric level structure at primes $p$ at which the group that defines
the Shimura variety ramifies. We describe ``good" $p$-adic integral models
of these Shimura varieties and study their \'etale local structure.
In particular, we exhibit a stratification of their (singular) special fibers
and give a partial calculation of the sheaf of nearby cycles.
\medskip

\noindent {\sl 2000 Mathematics Subject Classification:} Primary  14G35, 11G18; Secondary  14M15. 
\end{abstract}

\medskip
\bigskip
\bigskip

\section{Introduction} \label{intro}

\ \ \ The problem of constructing ``good'' models of Shimura varieties
over the ring of integers of the reflex field, or over its
completion at some prime ideal, has a long history. For Shimura
varieties of PEL type, which can be described as moduli spaces of
abelian varieties, it is desirable to define such a model by a
suitable extension of the moduli problem. The ultimate goal is to
calculate the Hasse-Weil zeta function.
For the local factor of the zeta function at a prime ideal
$\mathfrak p$, this comes down to a counting problem.
Namely, in the case of good reduction, one counts the number of points of the reduction modulo $\mathfrak p$
of the model over finite extensions of the residue field. 
In the case of bad reduction, one has to weight those points by the trace of Frobenius 
on the sheaf of nearby cycles.  

In the case of good reduction, good models of PEL Shimura varieties
were constructed by Kottwitz [K], when the group $G$ defining the Shimura variety has as
simple factors only groups of type $A$ or $C$. Furthermore, Kottwitz 
only considers hyperspecial level structure at $\mathfrak p$. In the more
general case of a parahoric level structure, integral models were
proposed in [RZ]. In two papers, G\"ortz ([G1], [G2]) proved that
these models are flat and have reasonable singularities, provided
the group $G$ only involves factors of type $A$ or $C$ and splits
over an unramified extension of $\Q_p$. This follows
work of Deligne-Pappas [DP] and de Jong [dJ] (we also mention related work by Chai-Norman [CN],
Faltings [F2] and Genestier [Ge]). On the negative side, it
was shown in [P] that when this last condition fails, the models
of [RZ] are not flat in general. The aim of this series of papers
is to find the correct modification of the proposed models of
[RZ] in the ramified case, to show that these modified models are
flat and have reasonable singularities, and to calculate the
weighting factor in the counting problem mentioned above, i.e.\
the trace of the Frobenius on the sheaf of nearby cycles. By a
procedure that is by now well-known, these questions reduce to
problems on the local models of Shimura varieties. The advantage
of this approach is that we are then dealing with varieties which
can be defined in terms of linear algebra. This is the point of
view taken in this paper. In the last section we indicate the
implications of our results for the original problem of
constructing suitable models of Shimura varieties.

Assume that the group $G$ defining the
Shimura variety splits over a ramified extension of $\Q_p$. As
was mentioned above, in this case the ``naive'' local models of
[RZ] are not flat in general. One wants to define closed
subschemes of the naive local models which are flat and to
understand the structure of their special fibers. Two typical
cases in which ramification occurs for a PEL Shimura variety are
the following:
\smallskip

(i) $G_{\Q_p}$ is of the form $G_{\Q_p}={\rm Res}_{F/\Q_p}G'$, where $G'$
is a quasi-split group over $F$  which splits over an unramified
extension of $F$ and $F/\Q_p$ is a ramified extension.

\smallskip

(ii) $G_{\Q_p}$ is the group of unitary similitudes corresponding
to a ramified quadratic extension of $\Q_p$.

\smallskip

The case (ii), first addressed in [P], presents  challenges
that are of different nature from those in case (i) (but compare Remark \ref{unitary}
below).
We intend to take up this case in subsequent work. Here we will be concerned
with the case (i). Although we will only consider the cases where
$G'={\rm GL}_d$ or $G'={\rm GSp}_{2g}$, our method applies more
generally, comp.\ section \ref{general}. Loosely speaking, the
method developed here allows us to deal with ramification caused
by restriction of scalars. In this introduction we will
concentrate on the case $G'={\rm GL}_d$ which brings out better
the outlines of our approach.  

\smallskip

Let $F_0$ be a complete discretely valued field with ring of
integers ${\cal O}_{F_0}$ and perfect residue field. Let $F$ be a
totally ramified extension of degree $e$ contained in a fixed
separable closure $F_0^{\rm sep}$ of $F_0$. Let ${\cal O}_F$ be
the ring of integers of $F$ and $\pi$ a uniformizer which is the
root of an Eisenstein polynomial $Q(T)\in {\cal O}_{F_0}[T]$. Let
$K$ be the Galois hull of $F$ in $F_0^{\rm sep}$, with ring of
integers ${\cal O}_K$ and residue field $k'$.

Let $V$ be an $F$-vector space of dimension $d$. Fix an $F$-basis
$e_1,\ldots, e_d$ of $V$ and let $\Lambda_i$, $0\leq i\leq d-1$
be the ${\cal O}_F$-lattice in $V$ spanned by
$\pi^{-1}e_1,\ldots, \pi^{-1}e_i, e_{i+1},\ldots, e_d$. For a
subset $I=\{ i_0<\cdots <i_{m-1}\} \subset \{ 0,\ldots, d-1\}$ we
obtain a periodic lattice chain $\Lambda_I$ in $V$ which is given
by all multiples of the lattices $\Lambda_i$ with $i\in I$.

Choose for each embedding $\varphi: F\to F_0^{\rm sep}$ an integer
$r_{\varphi}$ with $0\leq r_{\varphi}\leq d$. Set $r=\Sigma
r_{\varphi}$. Then the naive local model $M_I^{\rm naive}=M^{\rm
naive}({\cal O}_F, \Lambda_I, {\bf r})$ associated to the lattice
chain $\Lambda_I$ and to ${\bf r}=(r_{\varphi})$ parametrizes
points ${\cal F}_j$ in the Grassmannian of subspaces of rank $r$
of $\Lambda_{i_j, S}=\Lambda_{i_j}\otimes_{{\cal O}_{F_0}}{\cal
O}_S$ which are ${\cal O}_F$-stable and on which the
representation of ${\cal O}_F$ is prescribed in terms of ${\bf
r}$, and which are compatible with varying $j=0,\ldots, m-1$. It
is a projective scheme defined over $\Spec {\cal O}_E$, where
$E=E(V, {\bf r})$ is the reflex field. Let $k$ be the residue
field of ${\cal O}_E$. A fact which is important for
our analysis is that the special fiber of $M^{\rm naive}_I$ can
be considered as a closed subscheme of the affine partial flag
variety $\flag_I= {\rm GL}_d(k[[\Pi]])/P_I$, 
$$
i:M_I^{\rm
naive}\otimes_{{\cal O}_E}k\hookrightarrow \flag_I\ \ .
$$ 
We
choose an ordering of the embeddings $\varphi$. The basic new
ingredient of the present paper is the splitting model ${\cal
M}_I= {\cal M}({\cal O}_F, \Lambda_I, {\bf r})$. This is a
projective scheme over $\Spec {\cal O}_K$ representing the functor 
which to  an ${\cal
O}_K$-scheme $S$ associates the set of commutative diagrams of
${\cal O}_F\otimes_{{\cal O}_{F_0}}{\cal O}_S$-morphisms resp.\
-inclusions,

$$ \matrix{\Lambda_{i_0, S}&\rightarrow&\Lambda_{i_1,S}&\rightarrow&\cdots
&\rightarrow
&\Lambda_{i_{m-1},S}&\buildrel\pi\over\rightarrow&\Lambda_{i_0,S}\cr
\cup&&\cup&&&&\cup&&\cup\cr
\F^{e}_0&\rightarrow&\F^{e}_1&\rightarrow&\cdots
&\rightarrow&\F^{e}_{m-1}&\rightarrow&\F^{e}_0\cr
\cup&&\cup&&&&\cup&&\cup\cr
\F^{e-1}_0&\rightarrow&\F^{e-1}_1&\rightarrow&\cdots
&\rightarrow&\F^{e-1}_{m-1}&\rightarrow&\F^{e-1}_0\cr
\cup&&\cup&&&&\cup&&\cup\cr \vdots&&\vdots&&&&\vdots&&\vdots\cr
\cup&&\cup&&&&\cup&&\cup\cr
\F^{1}_0&\rightarrow&\F^{1}_1&\rightarrow&\cdots
&\rightarrow&\F^{1}_{m-1}&\rightarrow&\F^{1}_0\cr} $$

\medskip\noindent
which satisfy the following conditions:
\smallskip
\begin{itemize}
\item[a)] ${\cal F}_t^j$ is locally on $S$ an ${\cal O}_S$-direct summand of $\Lambda_{i_t,S}\ \ \hbox{of rank}\ \sum_{l=1}^j r_{l}\ \ .$

\smallskip

\item[b)] For each $a\in {\cal O}_F$ and $j=1,\ldots, e$ $$(a\otimes 1-1\otimes
\varphi_j(a)) ({\cal F}_t^j)\subset {\cal F}_t^{j-1}\ \ .$$ Here
we have set ${\cal F}^0_t=(0)$.
\end{itemize}
We obtain an ${\cal O}_K$-morphism $$\pi_I:{\cal
M}_I\longrightarrow M_I^{\rm naive}\otimes_{{\cal O}_E}{\cal
O}_K\ \ ,$$ given by $\{ {\cal F}_t^j\}_{t,j}\longmapsto \{ {\cal
F}_t^{e}\}_t$.

Our crucial observation is that ${\cal M}_I$ can be identified
with a twisted product of unramified local models for ${\rm
GL}_d$ over ${\cal O}_K$, $M^{l}_I= M({\cal O}_K,
\Lambda_I\otimes_{{\cal O}_F, \varphi_{l}}{\cal O}_K, r_{l})$,
for $l=1,\ldots, e$: $${\cal M}_I=
M^1_I\tilde\times\ldots\tilde\times M^e_I\ \ .$$ Let us define the
canonical local model $M_I^{\rm can}$ as the scheme-theoretic
image of the composed morphism $${\cal
M}_I\buildrel\pi_I\over\longrightarrow M_I^{\rm
naive}\otimes_{{\cal O}_E}{\cal O}_K\longrightarrow M_I^{\rm
naive}\ \ .$$ We may then state the following result:

\medskip\noindent
{\bf Theorem A:} {\it a) $M_I^{\rm can}$ is the flat closure of
the generic fiber $M_I^{\rm naive}\otimes_{{\cal O}_E}E$ in
$M_I^{\rm naive}$, and it coincides with the local model
$M_I^{\rm loc}$ defined in [PR], \S 8. Its special fiber is
reduced, and all its irreducible components are normal and with
rational singularities.

\smallskip

b) The special fiber $M^{\rm can}_I\otimes_{{\cal O}_E}k$ is the
union of the Schubert strata in $\flag_I$ for all $w\in \tilde
W_I\setminus \tilde W/\tilde W_I$ in the $\mu$-admissible set
([KR]), for $\mu=\omega_{r_1}+\ldots+\omega_{r_e}$, $$M^{\rm
can}_I\otimes_{{\cal O}_E}k= \bigcup\limits_{w\in {\rm
Adm}_I(\mu)}{\cal O}_w\ \ .$$ }

\smallskip\noindent
Here $\tilde W$ denotes the extended affine Weyl group for ${\rm
GL}_d$ and $\tilde W_I$ the parabolic subgroup corresponding to
$I$.

\medskip
\noindent The basic ingredients of the proof of Theorem A are the
presentation of the splitting model as a twisted direct product of
unramified local models and the results of G\"ortz [G1] on these
unramified local models. We also need results of Haines-Ng\^o [HN2]
and G\"ortz [G3] on affine Weyl groups. When the integers $r_{\varphi}$ 
differ by at most one
amongst each other, we conjecture that $M_I^{\rm can}=M_I^{\rm
naive}$, comp.\ [PR]. Similarly, it seems reasonable to expect in
the symplectic case that the canonical local model coincides with
the naive local model, i.e., that the naive local model is flat
in this case, comp.\ [G3]. Theorem A seems to indicate that
$M^{\rm can}_I$ is the ``correct'' way to extend its generic fiber
into an integral model. Even though $M^{\rm can}_I$ does not
represent a good moduli problem, the geometric points of $M^{\rm
can}_I\otimes_{{\cal O}_E}k$ can be described as a subset of
$M^{\rm naive}_I\otimes_{{\cal O}_E}k$ and $M^{\rm can}_I$
satisfies a maximal property with respect to the morphism from
${\cal M}_I$ to $M^{\rm naive}_I$. The situation is therefore
quite similar to the solution of an orbit problem by its coarse
moduli space.

\smallskip

Our second use of the presentation of ${\cal M}_I$ as a twisted
direct product of unramified local models concerns the
calculation of the complex of nearby cycles of $M_I^{\rm can}$.
Let us assume that the residue field $k$ of ${\cal O}_E$ is
finite. Let us denote by $$R\Psi_K^{M_I^{\rm can}} =
R\Psi(M_I^{\rm can}\otimes_{{\cal O}_E}{\cal O}_K/ {\cal O}_K)
\Q_{\ell}[d]\left( \hbox{${d\over 2}$}\right)$$ the adjusted
complex of nearby cycles, where $d$ denotes the relative
dimension of $M_I^{\rm can}$. This is a  perverse
$\overline{\Q}_{\ell}$-sheaf on $M_I^{\rm
can}\otimes_{{\cal O}_E}\overline k$, which we may regard as a
$P_I$-equivariant perverse $\overline{\Q}_{\ell}$-sheaf on
$\flag_I\otimes_k\overline k$, equipped with an action by ${\rm
Gal}(F_0^{\rm sep}/K)$.

\medskip\noindent
{\bf Theorem B:} {\it There is an isomorphism of perverse
$\overline{\Q}_{\ell}$-sheaves with ${\rm Gal}(F_0^{\rm
sep}/K)$-action
$$R\Psi_K^{M_I^{\rm can}} = R\Psi_K^{M^1_I}*\cdots* R\Psi_K^{M^e_I}\ \ .$$
}

Here on the right hand side there appears the convolution in the
sense of Lusztig of the adjusted complexes of nearby cycles of
the unramified local models $M_I^j$, $j=1,\ldots, e$. The latter
perverse sheaves on $\flag_I\otimes_k\overline k$ are known due
to the solution of the Kottwitz conjecture by Haines and Ng\^o
[HN1].
\smallskip

In section \ref{general} we extend our construction of the
splitting model to the general ramified PEL case. The splitting
model comes equipped with a morphism to the naive PEL local model
of [RZ]. The scheme-theoretic image of this morphism is a closed
subscheme of the naive local model; as our examples indicate, it
is reasonable to expect that in many cases (see \S\ref{general}
for details) this is the ``canonical" flat local model. However,
as A.~Genestier pointed out to us, this is not true for the
 even orthogonal group. Also, this is not true in general in
the case of a unitary group corresponding to a ramified quadratic
extension. However, we believe that even in these cases the
methods of the present paper will turn out to be useful. In the
last section we briefly indicate how to construct integral models
of the relevant moduli spaces of abelian varieties with
additional structure, defined by the splitting local models and
canonical local models.

\medskip

We thank V.~Drinfeld, A.~Genestier, U.~G\"ortz and T.~Haines for
interesting discussions and the referees for their remarks and
suggestions.

\setcounter{equation}{0}

\section{General notations} \label{notations}

Most of the time, we will follow the notations and assumptions of [R-P],
\S 2. In particular, $F_0$ is a complete discretely valued field with ring
of integers $\O_{F_0}$, uniformizer $\pi_0$ and perfect residue field. We
fix a separable closure $F_0^{\rm sep}$ of $F_0$. Let $F$ be a totally
ramified separable extension of degree $e$ of $F_0$ with ring of integers
$\O_F$. Let $\pi$ be a uniformizer of ${\cal O}_F$ which is a root of the
Eisenstein polynomial
\begin{eq}\label{2.1}
Q(T)=T^e+\sum_{k=0}^{e-1}b_kT^k,\ \ b_0\in
\pi_0\cdot \O^{\times}_{F_0},\ b_k\in
(\pi_0).\end{eq}
Let us denote by $K$ the Galois hull of $F$ in $F_0^{\rm sep}$
and let $\O_K$ be the ring of integers of $K$; denote by $k'$
the residue field of $\O_K$. Let us choose an
ordering of the embeddings $\phi: F\to F^{\rm sep}_0$ and for
$i\in\{1,\cdots, e\}$ let us set $a_i=\phi_i(\pi)$.
We have
\begin{eq}
 \O_{F_0}[T]/(Q(T))\simeq \O_F\
\end{eq}
given by $T\mapsto\pi $.
For $i=1,\ldots, e$, we set
\begin{eq}
\ \ \ \ \ \ Q^i(T)=\prod^e_{j=i}(T-a_j),\quad
Q_i(T)=\prod^{i-1}_{j=1}(T-a_j)\in \O_K[T], \quad
\O_K^{(i)}=\O_K[T]/(Q^i(T))\ ,
\end{eq}
so that $Q^1(T)=Q(T)$, and $Q_i(T)Q^i(T)=Q(T)$. There are natural
surjective $\O_K$-algebra homomorphisms $$ \phi^i:
\O_F\otimes_{\O_{F_0}}\O_K\simeq \O_K[T]/(Q(T))\to
\O_K[T]/(Q^i(T))=\O^{(i)}_K $$ obtained by sending $\pi\otimes 1$ to $T$.

There are exact sequences $$ \O_K[T]/(Q(T))\buildrel Q_i(T)\over\to
\O_K[T]/(Q(T))\buildrel Q^i(T)\over\to \O_K[T]/(Q(T))\ , $$ $$
\O_K[T]/(Q(T))\buildrel Q^i(T)\over\to \O_K[T]/(Q(T))\buildrel
Q_i(T)\over\to \O_K[T]/(Q(T))\ , $$ with the image and the kernel of each
morphism $\O_K$-free. We conclude that if $S$ is an $\O_K$-scheme, there
are functorial isomorphisms
\begin{eq}\label{kerim1}
\ \ \ \ \ \O^{(i)}_K\otimes_{\O_K}\O_S\simeq {\rm Im}(Q_i(T)\ |\ \O_S[T]/(Q(T)))
=
\ker( Q^i(T) |\ \O_S[T]/(Q(T)))\ ,
\end{eq}
the first one obtained by multiplying with $Q_i(T)$.

\part{}

\section{The ``naive" local models for $G=\Res_{F/F_0}\GL_d$} \label{standGl}
\setcounter{equation}{0}

Now let $V$ be an $F$-vector space of dimension $d$.
Fix an $F$-basis $e_1, \ldots e_d$ of $V$ and let
$\Lambda_i$, $0\leq i\leq d-1$, be the free $\O_F$--module
of rank $d$ with basis $e^i_1:=\pi^{-1}e_1,\ldots, e^i_i:=\pi^{-1}e_i$,
$e^i_{i+1}:=e_{i+1},\ldots ,e^i_d:=e_d$. Let us choose a subset
$I=\{i_0<\cdots <i_{m-1}\}\subset
\{0,\ldots, d-1\}$ and consider the $\O_F$-lattice chain $\Lambda_I$
in $V$ which is given by all multiples of the lattices $\Lambda_i$ with
$i\in I$.
\medskip

Let us choose for each  embedding $\varphi : F\to F_0^{\rm sep}$
an integer $r_{\varphi}$ with
$0\leq r_{\varphi}\leq d$. Set $r=\sum_{\varphi}r_\phi$.
Associated to these data we have the {\it reflex field} $E$, a
finite extension of $F_0$ contained in $F_0^{\rm sep}$ with

\begin{eq}\label{2.3}
{\rm Gal}(F_0^{\rm sep}/E) =\{ \sigma\in {\rm Gal}(F_0^{\rm
sep}/F_0);\ r_{\sigma\varphi}=r_{\varphi},\ \forall \varphi\}\ \ .
\end{eq}

We also have a cocharacter $\mu: {\Gm}{/F_0^{\rm sep}}\rightarrow
(\Res_{F/F_0}\GL_d){/F_0^{\rm sep}}$ given by $(1^{r_\phi},
0^{d-r_\phi})_\phi$. The conjugacy class of $\mu$ is defined over the
reflex field $E$. Let $\O_E$ be the ring of integers in $E$ and $k$ its
residue field.
\smallskip

The ``naive" local model $M^{\rm naive}_I=M(\O_F, \Lambda_I, {\bf r})$ of
[RZ], Definition 3.27 for $G=\Res_{F/F_0}\GL_d$, the cocharacter $\mu$ and the lattice
chain $\Lambda_I$, is the $\O_E$--scheme which represents the following
functor: To each $\O_E$--scheme $S$, we associate the set $M^{\rm
naive}_I(S)$ of collections $\{\F_t\}_t$ of
$\O_F\otimes_{\O_{F_0}}\O_S$-submodules of $\Lambda_{i_t,
S}:=\Lambda_{i_t}\otimes_{\O_{F_0}}\O_S$ which fit into a commutative
diagram 
$$ 
\matrix{\Lambda_{i_0,
S}&\rightarrow&\Lambda_{i_1,S}&\rightarrow&\cdots &\rightarrow
&\Lambda_{i_{m-1},S}&\buildrel\pi\over\rightarrow&\Lambda_{i_0,S}\cr
\cup&&\cup&&&&\cup&&\cup\cr \F_0&\rightarrow&\F_1&\rightarrow&\cdots
&\rightarrow&\F_{m-1}&\rightarrow&\F_{0,S}\ ,\cr} 
$$ 
(with the morphisms
$\Lambda_{i_t, S}\to \Lambda_{i_{t+1}, S}$ of the first row induced by the
lattice inclusions $\Lambda_{i_t}\to \Lambda_{i_{t+1}}$). We require the
following conditions:

i) $\F_t$ is Zariski locally on $S$ a $\O_S$-direct summand
of $\Lambda_{i_t,S}$ of rank $r$,

ii) for $a\in \O_F$, we have $$ \det(a\ |\
\F_t)=\prod_\phi\phi(a)^{r_\phi}\ ,
$$
 where this last identity is meant as
an identity of polynomial functions on $\O_F$ (comp. [K], or [RZ], 3.23 (a)).

It is clear that
this functor is represented by a projective scheme
over $\Spec\O_E$.
\medskip

Consider the group scheme $\Gg_I$ over $\Spec\O_{F_0}$
\begin{eq}
\Gg_I:=\underline{\rm Aut}_{\O_F}(\Lambda_I)
\end{eq}
with $S$-valued points the $\O_F\otimes_{\O_{F_0}}\O_S$-automorphisms of
the lattice chain $\Lambda_I\otimes_{\O_{F_0}}\O_S$. A simple extension of
the arguments of [RZ] Appendix (see loc. cit.,  Proposition A.4, also [P] Theorem 2.2) shows that $\Gg_I$ is
smooth over $\Spec\O_{F_0}$, comp.\ Remark \ref{carry} below. Often we
will use the base change of $\Gg_I$ to $\Spec\O_E$, which we will denote
by the same symbol.

\begin{Remark} \label{carry}
{\rm The arguments of the proof of [RZ] Proposition A.4 carry over with essentially no changes to the
following situation. Let $(F_0, {\cal O}_{F_0}, \pi_0)$ be as in section
\ref{notations}. Let $\cal O$ be an ${\cal O}_{F_0}$-order in a
semi-simple $F_0$-algebra ${\cal O}\otimes_{{\cal O}_{F_0}}F_0$. For the
purposes of the present paper we may assume that ${\cal O}$ is
commutative, i.e.\ ${\cal O}\otimes_{{\cal O}_{F_0}}F_0$ is a product of
field extensions of $F_0$. 
Let $\Pi\in{\cal O}$ be an element with
$\pi_0\in (\Pi)$. ( Let us remark that in loc. cit. 
it is assumed in addition that $\cal O$ is a maximal order and $\Pi$ gives
a uniformizer in each component of ${\cal O}\otimes_{{\cal O}_{F_0}}F_0$.) 
Let $V$ be a finite-dimensional $F_0$-vector space which
is an ${\cal O}\otimes_{{\cal O}_{F_0}}F_0$-module. An ${\cal O}$-lattice
in $V$ is a ${\cal O}_{F_0}$-lattice in $V$ stable under ${\cal O}$. A
$({\cal O},\Pi)${\it -periodic lattice chain} is a chain of inclusions of
${\cal O}$-lattices in $V$,
$$\subset\Lambda_{i-1}\subset\Lambda_i\subset\ldots\quad ,\ i\in {\bf Z}\
\ ,$$ such that
\smallskip

(i) $\exists r:\  \Lambda_{i-r}=\Pi\Lambda_i\ \ ,\ \ \forall i\in {\bf
Z}$.
\smallskip

(ii) $\Lambda_i/\Lambda_{i-1}$ is a free ${\cal O}/\Pi {\cal O}$-module
$\forall i\in {\bf Z}$.

\smallskip

Let us fix a $({\cal O},\Pi)$-periodic lattice chain ${\cal L}$. Let $S$
be a ${\cal O}_{F_0}$-scheme such that $\pi_0$ is locally nilpotent on
$S$. A {\it chain of ${\cal O}\otimes_{{\cal O}_{F_0}} {\cal O}_S$-modules
of type $({\cal L})$ on $S$} is given by a chain of ${\cal
O}\otimes_{{\cal O}_{F_0}}{\cal O}_S$-module homomorphisms,
$$\ldots\buildrel\varrho\over\longrightarrow
M_{i-1}\buildrel\varrho\over\longrightarrow
M_i\buildrel\varrho\over\longrightarrow\ldots$$ such that the following
conditions are satisfied.

\smallskip

(i) $\varrho^r=\Pi$.

\smallskip

(ii) Locally on $S$ there exist isomorphisms of ${\cal O}\otimes_{{\cal
O}_{F_0}}{\cal O}_S$-modules, $$M_i\simeq \Lambda_i\otimes_{{\cal
O}_{F_0}} {\cal O}_S\ \ ,\ \ M_i/\varrho(M_{i-1})\simeq
\Lambda_i/\Lambda_{i-1}\otimes_{{\cal O}_{F_0}}{\cal O}_S\ \ .$$ The proof
of Prop.\ A.4 of loc.\ cit.\ shows then that any chain $\{ M_i\}$ of
${\cal O}\otimes_{{\cal O}_{F_0}}{\cal O}_S$-modules of type $({\cal L})$
on $S$ is locally on $S$ isomorphic to ${\cal L}\otimes_{{\cal
O}_{F_0}}{\cal O}_S$, and that the functor on $(Sch/S)$, $$S'\longmapsto
{\rm Aut}(\{ M_i\otimes_{{\cal O}_S}{\cal O}_{S'}\} )$$ is representable
by a smooth group scheme over $S$.

}
\end{Remark}

\section{Affine flag varieties for $\GL_d$} \label{aFlag}
\setcounter{equation}{0}

If $R$ is a $k$-algebra, a lattice in $R((\Pi))^d$ is by definition a
sub-$R[[\Pi]]$-module $\L$ of $R((\Pi))^d$ which is locally on $\Spec R$
free of rank $d$ and such that $\L\otimes_{R[[\Pi]]}R((\Pi))=R((\Pi))^d$.
(Here $R[[\Pi]]$, resp. $R((\Pi))$ denotes the power series ring, resp.
Laurent power series ring in the indeterminate $\Pi$ over $R$).
Equivalently, a lattice is a sub-$R[[\Pi]]$-module $\L$ of $R((\Pi))^d$
such that $\Pi^NR[[\Pi]]^d\subset \L\subset \Pi^{-N}R[[\Pi]]^d$ for some
$N$ and such that $\Pi^{-N}R[[\pi]]^d/\L$ is a locally free $R$-module.

Recall ( [BL] ) that the affine Grassmannian $\grass$ over $k$ associated to
$\GL_d$ is the Ind-scheme over $\Spec k$ which represents the functor on
$k$-algebras which to a $k$-algebra $R$ associates the set of lattices
$\L$ in $R((\Pi))^d$. The affine Grassmannian can be identified with the
fpqc quotient $\GL_d(k((\Pi)))/\GL_d(k[[\Pi]])$ where $\GL_d(k((\Pi)))$,
resp. $\GL_d(k[[\Pi]])$ is the Ind-group scheme, resp. group scheme over
$\Spec k$ whose $R$-rational points is $\GL_d(R((\Pi)))$, resp.
$\GL_d(R[[\Pi]])$.

For each $i\in \{0,\ldots, d-1\}$, we will denote by $\ti\Lambda_i$ the
$k[[\Pi]]$-lattice in $$ k((\Pi))^d=k((\Pi))\ti e_1\oplus\cdots \oplus
k((\Pi))\ti e_d $$ which is generated by $\Pi^{-1}\ti e_1,\ldots,
\Pi^{-1}\ti e_i,\ti e_{i+1},\ldots, \ti e_d$.

Denote by $P_I$, resp. $P'_I$, the parahoric
subgroup scheme of $\GL_d(k((\Pi)))$, resp. $\SL_d(k((\Pi)))$,
whose $k$-valued points stabilize the lattice chain
\begin{eq}
\ti \Lambda_{i_0} \subset \ti \Lambda_{i_1} \subset \cdots \subset \ti \Lambda_{i_{m-1}} \subset
\Pi^{-1}\ti \Lambda_{i_0} \ .
\end{eq}
If $I=\{0,\ldots, d-1\}$, then
$P_I$, resp. $P'_I$, is an Iwahori subgroup scheme of $\GL_d(k((\Pi)))$,
resp. $\SL_d(k((\Pi)))$.

\medskip

For every nonempty subset $I=\{i_0<\cdots <i_{m-1}\}\subset \{0,\ldots, d-1\}$, we have the partial
affine flag variety $\flag_I$ whose $R$-rational points
parametrize lattice chains in $R((\Pi))^d$
\begin{eq}
\L_{0}\subset \L_1\subset \cdots \subset \L_{m-1}\subset \Pi^{-1}\L_0
\end{eq}
with $\L_{t+1}/\L_{t}$, resp.\  $\Pi^{-1}\L_0/\L_{m-1}$ locally free
$R$-modules of rank $i_{t+1}-i_t$ for $t=0,\ldots , m-2$, resp.\
$(d+i_0)-i_{m-1}$.
\smallskip

The  affine Grassmannian variety corresponds to the choice $I=\{0\}$,
while the full affine flag variety corresponds to $I=\{0,\ldots, d-1\}$.
The Ind-group scheme $\GL_d(k((\Pi)))$ acts on the partial affine flag
variety $\flag_I$ and we can identify $\flag_I$
($\GL_d(k((\Pi)))$-equivariantly)
 with the fpqc quotient

\begin{eq}
\flag_I=\GL_d((k((\Pi)))/P_I\ .
\end{eq}
\medskip

Given $r\in\Z$, we may also consider the special partial
affine flag variety $\flag^r_I$ whose $R$-rational points
parametrize lattice chains in $R((\Pi))^d$
\begin{eq}
\L_{0}\subset \L_1\subset \cdots \subset \L_{m-1}\subset \Pi^{-1}\L_0
\end{eq}
such that:
\smallskip

i)  $\L_{t+1}/\L_{t}$, resp.\ $\Pi^{-1}\L_0/\L_{m-1}$ are locally free
$R$-modules of rank $i_{t+1}-i_t$ for $t=0\ldots m-2$, resp.\
$(d+i_0)-i_{m-1}$,

ii) $\wedge^d\L_0=\Pi^rR[[\Pi]]^d$ (as a submodule of $\wedge^dR((\Pi))^d=R((\Pi))$).
\medskip

The special affine flag varieties $\flag^r_I$ for various $r$ are all
isomorphic to the fpqc quotient 
$$
 \SL_d(k((\Pi)))/P'_I 
 $$ (as abstract
Ind-schemes but not $\SL_d(k((\Pi)))$-equivariantly, unless $r=-i_0$). For
$I=\{0\}$, we obtain the special affine Grassmannian $$ \grass^r\ \simeq\
\SL_d(k((\Pi)))/\SL_d(k[[\Pi]])\ . $$

\smallskip

Now fix an identification  $\O_F\otimes_{\O_{F_0}}k=k[[\Pi]]/(\Pi^e)$ and
$\O_F\otimes_{\O_{F_0}}k$-isomorphisms 
$$
\Lambda_{i_t}\otimes_{\O_{F_0}}k\simeq
\ti\Lambda_{i_t}\otimes_{k[[\Pi]]}k[[\Pi]]/(\Pi^e) 
$$ which induce a
$k[[\Pi]]/(\Pi^e)$-module chain isomorphism $$
\Lambda_{I}\otimes_{\O_{F_0}}k\simeq
\ti\Lambda_{I}\otimes_{k[[\Pi]]}k[[\Pi]]/(\Pi^e)\ . $$

Let $R$ be a $k$-algebra. For an $R$-valued point $\{\F_{t}\}_t$ of
$M^{\rm naive}_I$, we have
\begin{eq}
\F_{t}\subset \Lambda_{i_t}\otimes_{\O_{F_0}}R=\tilde\Lambda_{i_t}\otimes_{k[[\Pi]]}R[[\Pi]]/(\Pi^e)\ .
\end{eq}
Let $\L_t\subset \tilde \Lambda_{i_t}\otimes_{k[[\Pi]]}R[[\Pi]]$ be the
inverse image of $\F_t$ under the canonical projection 
$$
 \tilde
\Lambda_{i_t}\otimes_{k[[\Pi]]}{R[[\Pi]]}\to
\tilde\Lambda_{i_t}\otimes_{k[[\Pi]]}R[[\Pi]]/(\Pi^e) \ , 
$$
 so that we
have 
$$
 \Pi^e\ti \Lambda_{i_t}\otimes_{k[[\Pi]]}R[[\Pi]]\subset\L_t\subset
\tilde \Lambda_{i_t}\otimes_{k[[\Pi]]}R[[\Pi]]\ . 
$$
 Then $\{\L_t\}_t$
gives an $R$-valued point of $\flag_I$. In this way, we obtain a morphism
\begin{eq}
i: {M^{\rm naive}_I}\otimes_{\O_E} k\longrightarrow \flag_I
\end{eq}
which is a closed immersion (of Ind-schemes).

\section{The splitting model for $G=\Res_{F/F_0}\GL_d$} \label{splitGl}
\setcounter{equation}{0}

Fix $I=\{ i_0< i_1<\ldots < i_{m-1}\} \subset \{ 0,\ldots, d-1\}$.
Consider the functor $\M_I=\M(\O_F, \Lambda_I, {\bf r})$ on $({\rm
Schemes}/\Spec\O_K)$ which to a $\O_K$-scheme $S$ associates the set
$\M_I(S)$
 of collections $\{\F^j_t\}_{j,t}$ of $\O_F\otimes_{\O_{F_0}}\O_S$-submodules
of $\Lambda_{i_t, S}$ which fit into a commutative diagram 
$$
\matrix{\Lambda_{i_0, S}&\rightarrow&\Lambda_{i_1,S}&\rightarrow&\cdots
&\rightarrow
&\Lambda_{i_{m-1},S}&\buildrel\pi\over\rightarrow&\Lambda_{i_0,S}\cr
\cup&&\cup&&&&\cup&&\cup\cr
\F^{e}_0&\rightarrow&\F^{e}_1&\rightarrow&\cdots
&\rightarrow&\F^{e}_{m-1}&\rightarrow&\F^{e}_0\cr
\cup&&\cup&&&&\cup&&\cup\cr
\F^{e-1}_0&\rightarrow&\F^{e-1}_1&\rightarrow&\cdots
&\rightarrow&\F^{e-1}_{m-1}&\rightarrow&\F^{e-1}_0\cr
\cup&&\cup&&&&\cup&&\cup\cr \vdots&&\vdots&&&&\vdots&&\vdots\cr
\cup&&\cup&&&&\cup&&\cup\cr
\F^{1}_0&\rightarrow&\F^{1}_1&\rightarrow&\cdots
&\rightarrow&\F^{1}_{m-1}&\rightarrow&\F^{1}_0\cr} $$ and are such that:
\smallskip

a) $\F^j_t$ is Zariski locally on $S$ a $\O_S$-direct summand
of $\Lambda_{i_t, S}$ of rank $\sum_{l=1}^jr_l$.
\smallskip

b) For each $a\in \O_F$ and $j=1,\ldots, e$,
$$
(a\otimes 1-1\otimes \phi_j(a))(\F^{j}_{t})\subset
\F^{j-1}_{t}
$$
where the tensor products $a\otimes 1$, $1\otimes \phi_j(a)$ are in
$\O_F\otimes_{\O_{F_0}}\O_S$ and where, for each $t$, we set $\F^0_{t}=(0)$.
\smallskip

The functor $\M_I$ is obviously represented by a projective scheme over
$\Spec\O_K$. Note that there is an $\O_K$-morphism
\begin{eq}
\pi_I: \M_I\to M^{\rm naive}_I\otimes_{\O_E}\O_K
\end{eq}
given by $\{\F_{t}^j\}_{t, j}\mapsto \{\F^e_{t}\}_{t}$.
Indeed, if $\F_{t}=\F^e_{t}$ supports
a filtration $\{\F^j_{t}\}_j$
with the above properties, then the characteristic polynomial of the action of
$a\in \O_F$ on $\F_t$ is
\begin{eq}
\prod_{l=1}^e(T-\phi_l(a))^{r_l}
\end{eq}
and therefore $\F_{t}$ satisfies the condition ii) in the definition of
$M^{\rm naive}_I$.

\begin{prop} \label{generic}
The morphism $\pi_I$ induces an isomorphism
$$
\pi_I\otimes_{\O_K}K: \M_I\otimes_{\O_K}K\buildrel \sim\over \longrightarrow M^{\rm naive}_I\otimes_{\O_E}K
$$
on the generic fibers.
\end{prop}

\begin{Proof} To each $S$-valued point of $M^{\rm naive}_I$ with $S$ a $K$-scheme,
given by $\{\F_{t}\}_t$, we can associate an $S$-valued point of $\M_I$ by
considering, for each $k$, the filtration $\{\F^l_{t}\}_l$ associated  to
the grading on the $\O_F\otimes_{\O_{F_0}}K$-module $\F_{t}=\F^e_{t}$
given using the decomposition $$ \O_F\otimes_{\O_{F_0}}K\simeq
\oplus_{l=1}^eK, \quad a\otimes b\mapsto (b\phi_l(a))_{l=1,\ldots, e}\ .
$$ This gives a morphism inverse to $\pi_I\otimes_{\O_K}K$.\endproof
\end{Proof}
\bigskip

For each $l=1,\ldots ,e$, $t=0,\ldots, m-1$, set $\Xi^l_{i_t}=\Lambda_{i_t}\otimes_{\O_F,\phi_l}\O_K$
(an $\O_K$-lattice in $V\otimes_{F,\phi_l}K$). Denote by $\Xi^l_I$ the $\O_K$-lattice chain in
$V\otimes_{F, \phi_l}K$ given by the lattices $\{a_l^n\Xi^l_{i_t}\}_{t,n\in\Z}$. An ``essential"
part of the lattice chain $\Xi^l_I$ is
$$
\Xi^l_{i_0}\subset \Xi^l_{i_1}\subset \cdots\subset
\Xi^l_{i_{m-1}}\subset a_l^{-1}\Xi^l_{i_0}\ ,
$$
in the sense that each successive link $\Xi^l_i\subset \Xi^l_{i'}$
in the total lattice chain $\Xi^l_I$
is a multiple of one of the links in the part above.

Let $\Gg^l_I$ be the group scheme over $\Spec\O_K$ whose $S$-points are
the $\O_S$-automorphisms of the chain $\Xi^l_I\otimes_{\O_K}\O_S$ (once
again, a simple extension of [RZ] Prop. A.4 shows that this is a smooth
group scheme, comp.\ Remark \ref{carry}).

\medskip
Now if $S$ a scheme over $\Spec\O_K$, we obtain from $\Lambda_{I, S}$ a
$\O^{(l)}_K\otimes_{\O_{K}}\O_S$-lattice chain $\Lambda^l_{I, S}$ by
extending scalars via $$ \phi^{l}\otimes_{\O_K}\O_S:
\O_F\otimes_{\O_{F_0}}\O_S\simeq \O_S[T]/(Q(T))\to
\O_S[T]/(Q^{l}(T))=\O^{(l)}_K\otimes_{\O_K}\O_S\ . $$ An argument as in
the proof of (\ref{kerim1}), shows that we have functorial isomorphisms of
chains of $\O_S$-modules

\begin{eq}\label{caniso}
\Xi^{l}_{I, S}=\Lambda_I\otimes_{\O_F,\phi_l}\O_S
\simeq
\ker(\pi-a_{l}\ |\ \Lambda^{l}_I\otimes_{\O_K}\O_S)
\end{eq}
obtained by sending the element $\lambda\otimes 1$ of $\Lambda_I\otimes_{\O_F,\phi_l}\O_S
$
to the image of $Q^{l+1}(\pi)\cdot (\lambda\otimes 1)$
in $\Lambda^{l}_I\otimes_{\O_K}\O_S$.
\medskip

Denote by $\Gg^{(l)}_I$ the group scheme over $\Spec\O_K$ whose $S$-points
are the $\O^{(l)}_K\otimes_{\O_K}O_S$-automorphisms of the chain
${\bf \Lambda^l_{I, S}:=}\Lambda^l_I\otimes_{\O_K}\O_S$ (once again we can see that this is a
smooth group scheme, comp.\ Remark \ref{carry}). The isomorphism
(\ref{caniso}) induces a group scheme homomorphism
\begin{eq}\label{gpcaniso}
\Gg^{(l)}_I\to \Gg^{l}_I\ .
\end{eq}

\medskip

Now suppose that $\{\F^j_{t}\}_{j,t}$ is an $S$-valued point of $\M_I$.
For $l=1, \ldots ,e$, let us set
$$
\Psi^l_{i_t, S}=\ker(Q^{l}(\pi)\ |\ \Lambda_{i_t, S}/\F^{l-1}_t)\ ;
$$
this is an $\O^{(l)}_K\otimes_{\O_K}\O_S$-module. We also set
$$
\Upsilon^l_{i_t, S}:={\rm ker}(\pi-a_l\ |\ \Lambda_{i_t, S}/{\F^{l-1}_{t}})=
{\rm ker}(\pi-a_l\ |\ \Psi^l_{i_t, S})\ .
$$

We have $\O^{(l)}_K\otimes_{\O_K}\O_S$-module, resp. $\O_S$-module, homomorphisms
$$
\Psi^l_{i_t, S}\to \Psi^l_{i_{t+1}, S},\qquad \Psi^l_{i_{m-1}, S}
\buildrel T\over
\to \Psi^l_{i_0, S}
$$
resp.
$$
\Upsilon^l_{i_t, S}\to \Upsilon^l_{i_{t+1}, S},\qquad \Upsilon^l_{i_{m-1}, S}
\buildrel a_l\over
\to \Upsilon^l_{i_0, S}
$$
induced by the $\O_F\otimes_{\O_{F_0}}\O_S$-module homomorphisms
$$
\Lambda_{i_t, S}/\F^{l-1}_{t}\to \Lambda_{i_{t+1}, S}/\F^{l-1}_{{t+1}}, \qquad
\Lambda_{i_{m-1}, S}/\F^{l-1}_{m-1}\buildrel \pi\over
\to \Lambda_{i_{0}, S}/\F^{l-1}_{{0}}
$$
by taking the kernel of $Q^{l}(\pi)$, resp. of $\pi-\phi_l(\pi)=\pi-a_l$.

\begin{prop} \label{bundles}
a) The formation of $\Psi^l_{i_t, S}$, resp. of $\Upsilon^l_{i_t,S}$,
from $\{\F^j_{k}\}_{j,k}$ commutes with base change.

b) The $\O^{(l)}_K\otimes_{\O_K}\O_S$-module $\Psi^l_{i_t, S}$,
resp. the $\O_S$-module $\Upsilon^l_{i_t, S}$, is locally
on $S$ free of rank $d$.

c)  The chain $\Psi^l_{i_t, S}$ of $\O^{(l)}_K\otimes_{\O_K}\O_S$-modules
given by 
$$
 \cdots\to \Psi^l_{i_0, S}\to\cdots \to \Psi^l_{i_{m-1},
S}\buildrel \pi\over \to \Psi^l_{i_0, S}\to\cdots 
$$ 
is Zariski locally on
$S$ isomorphic to the chain of  $\O^{(l)}_K\otimes_{\O_K}\O_{S}$-modules
$\Lambda^l_I\otimes_{\O_K}\O_{S}$. Similarly, the chain $\Upsilon^l_{I,
S}$ of $\O_{S}$-modules on $S$ given by $$ \cdots\to \Upsilon^l_{i_0,
S}\to\cdots \to \Upsilon^l_{i_{m-1}, S}\buildrel a_l\over \to
\Upsilon^l_{i_0, S}\to\cdots $$ is Zariski locally on $S$ isomorphic to
the chain of $\O_{S}$-modules $\Xi^l_I\otimes_{\O_K}\O_{S}$.
\end{prop}

\begin{Proof}
The statements for the modules $\Upsilon^l_{i_t, S}$ follow from the
corresponding statements for the modules $\Psi^l_{i_t, S}$. Indeed, we can
see this fact using the functorial isomorphisms (\ref{caniso}) and the
fact that $$ \Upsilon^l_{i_t, S}= {\rm ker}(\pi-a_l\ |\ \Psi^l_{i_t, S})\
. $$ Write $Q^l(T)^{-1}(\F^{l-1}_{t})$ for the inverse image of
$\F^{l-1}_{t}\subset \Lambda_{i_t, S}$ under $\Lambda_{i_t,S}\to
\Lambda_{i_t, S}$ given by multiplication by $Q^l(T)$. Notice that since
$Q_l(T)(\F^{l-1}_t)=(0)$, by (\ref{kerim1}) we have $\F^{l-1}_t\subset
Q^l(T)(\Lambda_{i_t, S})$. Hence, there is an exact sequence 
$$ 0\to
\ker(Q^l(T)\ |\ \Lambda_{ i_t, S})\to Q^l(T)^{-1}(\F^{l-1}_{t})\to
\F^{l-1}_{t}\to 0\ . 
$$
 By (\ref{kerim1}), $\ker(Q^l(T)\ |\ \Lambda_{i_t,
S})\simeq \Lambda_{i_t, S}/Q_l(T)\Lambda_{i_t, S}$. Hence,
$Q^l(T)^{-1}(\F^{l-1}_{t})$ is a locally free $\O_S$-module of rank
$d(e-l+1)+\sum_{i=1}^{l-1}r_i$ whose formation commutes with base change.
The exact sequence 
$$ 
0\to\Lambda_{i_t,
S}/Q^l(T)^{-1}(\F^{l-1}_{t})\buildrel Q^l(T)\over\to \Lambda_{i_t,
S}/\F^{l-1}_{t}\to \Lambda_{i_k, S}/Q^l(T)\Lambda_{i_t, S}\to 0 
$$
 now
implies that $\Lambda_{i_t, S}/Q^l(T)^{-1}(\F^{l-1}_{t})$ is also
$\O_S$-locally free. Hence,  $Q^l(T)^{-1}(\F^{l-1}_{t})\subset
\Lambda_{i_t, S}$ is locally an $\O_S$-direct summand. Now 
$$ 
0\to
Q^l(T)^{-1}(\F^{l-1}_{t})/\F^{l-1}_t \to \Lambda_{i_t, S}/\F^{l-1}_{t}\to
\Lambda_{i_t, S}/Q^l(T)^{-1}(\F^{l-1}_{t})\to 0 
$$ 
implies that
$\Psi^l_{i_t,S}=Q^l(T)^{-1}(\F^{l-1}_{t})/\F^{l-1}_t$ is a
$\O^{(l)}_K\otimes_{\O_K}\O_S$-module which is locally free of rank
$d(e-l+1)$ as an $\O_S$-module and that its formation commutes with base
change in $S$. To show that $\Psi^l_{i_t,S}$ is locally on $S$ a free
$\O^{(l)}_K\otimes_{\O_K}\O_S$-module it is enough to show this for
$S=\Spec L$, $L$ a field. This is easy to see if $L$ is an extension of
$K$. If $L$ is an extension of $k'$, then
$\O^{(l)}_K\otimes_{\O_K}L=L[\pi]/(\pi^{e-l+1})$. In this case, there is a
$L[\pi]/(\pi^e)$-basis $f_1,\ldots, f_d$ of $\Lambda_{i_t}\otimes_{\O_K}L$
and $l-1\geq s_1\geq \cdots \geq s_d\geq 0$ such that
\begin{eq}
\F^{l-1}_k=L[\pi]/(\pi^e)\cdot \pi^{e-s_1}f_1\oplus
\cdots \oplus L[\pi]/(\pi^e)\cdot \pi^{e-s_d}f_d\ .
\end{eq}
Then
\begin{eq} \label{explicit}
\Psi^l_{i_t, S}={L[\pi]/(\pi^e)\cdot \pi^{l-1-s_1}f_1\oplus
\cdots \oplus L[\pi]/(\pi^e)\cdot \pi^{l-1-s_d}f_d
\over L[\pi]/(\pi^e)\cdot \pi^{e-s_1}f_1\oplus
\cdots \oplus L[\pi]/(\pi^e)\cdot \pi^{e-s_d}f_d}\ ,
\end{eq}
which is freely generated over $L[\pi]/(\pi^{e-l+1})$ by
the classes of $\pi^{l-1-s_1}f_1,\ldots ,\pi^{l-1-s_d}f_d$.

It remains to show (c) for $\Psi^l_{i_t, S}$, i.e that the chain
$\Psi^l_{I, S}$ is Zariski locally isomorphic to the chain
$\Lambda^l_I\otimes_{\O_K}\O_{S}$. Given (b), an extension of the
arguments in the proof of [RZ] Prop. A 4, p. 133 shows that it will be
enough to prove that the cokernels of $$ \Psi^l_{i_{t},S}\to
\Psi^l_{i_{t+1}, S}\ ,\quad \Psi^l_{i_{m-1},S}\buildrel \pi\over\to
\Psi^l_{i_{0}, S} $$ are Zariski locally on $S$ isomorphic to the
$\O^{(l)}_K\otimes_{\O_K}\O_S$-modules
$(\Lambda^l_{i_t}/\Lambda^l_{i_{t+1}})\otimes_{\O_K}\O_S$ and
$(\Lambda^l_{i_0}/T\Lambda^l_{i_{m-1}})\otimes_{\O_K}\O_S$ respectively,
comp.\ Remark \ref{carry}. In what follows, we will only deal with the
case of $\Psi^l_{i_{t},S}\to \Psi^l_{i_{t+1}, S}$, the case of
$\Psi^l_{i_{m-1},S}\buildrel \pi\over\to \Psi^l_{i_{0}, S}$ being similar.
Notice that $\Lambda^l_{i_t}/\Lambda^l_{i_{t+1}}$ is a module over
$\O^{(l)}_K/T\O^{(l)}_K=\O_K[T]/(T, Q^l(T))\simeq \O_K/(\varpi^{e-l+1})$
with $\varpi$ a uniformizer of $\O_K$. For simplicity of notation, set $$
R_l=\O_K/(\varpi^{e-l+1})\ . $$ The cokernel of $\Psi^l_{i_{t},S}\to
\Psi^l_{i_{t+1}, S}$ is a module over $R_1\otimes_{\O_K}\O_S$. Hence, it
is enough to assume that $S=\Spec R$ is an affine $R_1$-scheme and prove
the result in this case. In fact, since $\M_I$ is a Noetherian scheme, we
can also assume that $R$ is Noetherian. We can lift the $R_1$-chain
$\Lambda_{I, R_1}:=\Lambda_I\otimes_{\O_{F_0}}R_1$ to a chain of
$R_1[[T]]$-free modules $\ti\Lambda_I$ 
$$
 \ti\Lambda_{i_0}\subset
\ti\Lambda_{i_1}\subset \cdots\subset \ti\Lambda_{i_{m-1}}\subset
T^{-1}\ti\Lambda_{i_0} 
$$
 which are all $R_1[[T]]$-submodules of
$R_1((T))^d$ such that there is an isomorphism of
$R_1[T]/(Q(T))=\O_K\otimes_{\O_{F_0}}R_1$-chains $$
\ti\Lambda_{I}\otimes_{R_1[[T]]}R_1[[T]]/(Q(T))\simeq \Lambda_{I, R_1}\ .
$$ Now notice that since each $a_j$ is nilpotent in $R_1$, the elements
$T-a_j$ of $R_1[[T]]$ are invertible in $R_1((T))$ and hence the inverse
$Q^l(T)^{-1}$ makes sense in $R_1((T))$. The diagram corresponding to the
$R$-valued point $\{\F^j_{t}\}_{j,t}$ of $\M_I$ now provides us with a
diagram of $R[[T]]$-lattices in $R((T))^d$:

$$ \matrix{\ti\Lambda_{i_0, R}&\subset&\ti\Lambda_{i_1, R}&\subset &\cdots
&\subset& \ti\Lambda_{i_{m-1}, R}&\subset&T^{-1}\ti\Lambda_{i_0, R}\cr
\cup&&\cup&&&&\cup&&\cup\cr Q^l(T)^{-1}\L^{l-1}_0&\subset
&Q^l(T)^{-1}\L^{l-1}_1&\subset& \cdots&\subset
&Q^l(T)^{-1}\L^{l-1}_{m-1}&\subset &T^{-1}Q^l(T)^{-1}\L^{l-1}_0\cr
\cup&&\cup&&&&\cup&&\cup\cr \L^{l-1}_0&\subset &\L^{l-1}_1&\subset&
\cdots&\subset &\L^{l-1}_{m-1}&\subset &T^{-1}\L^{l-1}_0\cr
\cup&&\cup&&&&\cup&&\cup\cr Q(T)\ti\Lambda_{i_0,
R}&\subset&Q(T)\ti\Lambda_{i_1, R}&\subset &\cdots &\subset&
Q(T)\ti\Lambda_{i_{m-1}, R}&\subset&T^{-1}Q(T)\ti\Lambda_{i_0, R}\ ,\cr }
$$ 
where $\ti\Lambda_{i_t, R}=\ti\Lambda_{i_t}\otimes_{R_1[[T]]}R[[T]]$,
and $\L^{l-1}_t$, resp. $Q^l(T)^{-1} \L^{l-1}_t$, is the inverse image of
$\F^{l-1}_t$, resp. $Q^l(T)^{-1}(\F^{l-1}_t)$ under the surjection $$
\ti\Lambda_{i_t, R}\to \ti\Lambda_{i_t,
R}\otimes_{R[[T]]}R[[T]]/(Q(T))\simeq \Lambda_{i_t, R}\ . $$ Each quotient
created by the inclusion of any two modules in this diagram is a finitely
generated locally free $R$-module. In particular $$
Q^l(T)^{-1}\L^{l-1}_{t+1}/Q^l(T)^{-1}\L^{l-1}_t $$ is annihilated by $T$
and is $R$-locally free of rank equal to the rank of
$\ti\Lambda_{i_{t+1},S}/\ti\Lambda_{i_t, S}$. It now follows that the
cokernel of $$ \Psi^l_{i_{t},S}=Q^l(T)^{-1}\L^{l-1}_{t}/\L^{l-1}_t\to
Q^l(T)^{-1}\L^{l-1}_{t+1}/\L^{l-1}_{t+1}= \Psi^l_{i_{t+1},S} $$ is
isomorphic to $$
\left(Q^l(T)^{-1}\L^{l-1}_{t+1}/Q^l(T)^{-1}\L^{l-1}_t\right)\otimes_{R_1}R_l
$$ and therefore it is a locally free $R\otimes_{R_1}R_l$-module of the
expected rank. This concludes the proof. \endproof
\end{Proof}

\medskip
\bigskip

Now let $M^l_I:=M(\O_{K}, \Xi^l_I, r_l)$ be the (``unramified") local
model over $\Spec \O_K$ for $G={\GL_d}{/K}=\GL(V\otimes_{F,\phi_l}K)$,
$\mu$ given by $(1^{r_l}, 0^{d-r_l})$, and the lattice chain $\Xi^l_I$
([RZ]). By definition, $M^l_I=M(\O_{K}, \Xi^l_{I}, r_l)$ is the
projective scheme over $\O_{K}$ which classifies collections $\{\F_t\}_t$
of $\O_S$-submodules of $\Xi^l_{i_t, S}:=\Xi^l_{i_t}\otimes_{\O_K}\O_S$
which fit into a commutative diagram $$ \matrix{\Xi^l_{i_0,
S}&\rightarrow&\Xi^l_{i_1,S}&\rightarrow&\cdots &\rightarrow
&\Xi^l_{i_{m-1},S}&\buildrel a_l\over\rightarrow&\Xi^l_{i_0,S}\cr
\cup&&\cup&&&&\cup&&\cup\cr \F_0&\rightarrow&\F_1&\rightarrow&\cdots
&\rightarrow&\F_{{m-1}}&\rightarrow&\F_{0}\ \cr} $$ and are such that
$\F_t$ is Zariski locally on $S$ a $\O_S$-direct summand of
$\Xi^l_{i_t,S}$ of rank $r_l$.

\medskip

Let us denote by $\widetilde {\M}_I$ the scheme over $\Spec\O_K$ whose
$S$-points correspond to pairs $$ \widetilde{\M}_I(S):=(\{\F^l_{I}\}_{l}\
, \ \{\sigma^l_I\})_{l=2}^e\ , $$ where $\{\F^l_k\}_{l,k}$ is an
$S$-valued point of $\M_I$ and for $l=2,\ldots ,e$, $$ \sigma^l_I:
\Psi^l_{I, S}\buildrel\sim\over\to \Lambda^l_{I, S} $$ is an isomorphism
of chains of $\O^{(l)}_K\otimes_{\O_K}\O_S$-modules. The natural
projection morphism $$ q_I: \widetilde{\M}_I\to \M_I $$ is a torsor for
the smooth group scheme $\prod_{l=2}^e\Gg^{(l)}_I$ by the action
\begin{eq}\label{action1}
(g^l)_{l=2}^e\cdot  \left(\{\F^l_{t}\}_{l,t}\ , \ \{\sigma^l_I\}_{l=2}^e\right)=
\left(\{\F^l_{t}\}_{l,t}\ , \ \{g^l\cdot \sigma^l_I\}_{l=2}^e\right)\ .
\end{eq}
\medskip

Notice that an isomorphism $\sigma^l_I$ as above, in view of
(\ref{caniso}), induces an isomorphism of chains of $\O_S$-modules $$
\tau^l_I: \Upsilon^{l}_{I, S}\buildrel\sim\over\to \Xi^{l}_{I, S}\ ,
l=2,\ldots, e\ . $$ For $l=1$, $\Psi^1_{i_t,S}=\Lambda_{i_t, S}$ and
(\ref{caniso}) gives a canonical isomorphism $$ v_{I}:\Upsilon^1_{I,
S}\buildrel\sim\over\to \Xi^1_{I, S}\ . $$

Now if $\{\F^l_t\}_{l,t}$ is an $S$-valued point of
$\M_I$, then since $(\pi-a_l)\F^l_{t}\subset \F^{l-1}_t$
we can consider $\F^l_t/\F^{l-1}_t$ as an $\O_S$-submodule
of $\Upsilon^l_{i_t, S}=\ker(\pi-a_l\ |\ \Psi^l_{i_t, S})$.
Consequently, if $(\{\F^l_{t}\}_{l,t}\ , \{\ \sigma^l_I\}_{l=2}^e)$ is an $S$-valued point
of $\widetilde \M_I$, then we can consider
\begin{eq}
\tau^l_{i_t}(\F^l_{t}/ \F^{l-1}_t)\subset \Xi^{l}_{i_t, S}\ .
\end{eq}

For $l=2,\ldots, e$, the $\O_S$-modules $\tau^l_{i_t}(\F^l_t/ \F^{l-1}_t)$
are locally direct summands of $\Xi^{l}_{i_t, S}$ and they provide us with
an $S$-valued point of the ``unramified" local model $M^{l}_I$. For $l=1$,
the $\O_S$-modules $v_{i_t}(\F^1_t)$ are locally direct summands of
$\Xi^{1}_{i_t, S}$ and provide us with an $S$-valued point of the local
model $M^{1}_I$. We conclude that there is a morphism of $\O_K$-schemes 
$$
p_I: \widetilde \M_I\ \to\ \prod_{l=1}^eM^l_I 
$$ 
given by 
$$
(\{\F^l_{t}\}_{l,t}\ ,\ \{\sigma^l_{I}\}_{l=2}^e)\mapsto
\left(\{v_{t}(\F^1_t)\}_t, \{\sigma^2_t(\F^2_t/\F^{1}_t)\}_t, \ldots,
\{\sigma^{e}_{t}(\F^e_t/\F^{e-1}_t)\}_t \right)\ . 
$$

It is easy to see that the morphism $p_I$ is also a
$\prod_{l=2}^e\Gg^{(l)}_I$-torsor. Note that the corresponding
$\prod_{l=2}^e\Gg^{(l)}_I$-action on $\widetilde\M_I$ is different from
the action which produces the torsor $q_I: \widetilde\M_I\to \M_I$:

\begin{eq} \label{action2}
(g^l)_{l=2}^e\cdot  \left(\{\F^l_{t}\}_{l,t}\ , \
\{\sigma^l_I\}_{l=2}^e\right)= \left(\{ (\sigma^l_I)^{-1}\cdot g^l\cdot
\sigma^l_I( \F^l_{t})\}_{l,t}\ , \ \{g^l\cdot \sigma^l_I\}_{l=2}^e\right)\
.
\end{eq}

\medskip

In short, we have obtained a diagram of morphisms of schemes
over $\Spec\O_K$:

\begin{eq} \label{diagram}
\matrix{&&&\widetilde\M_I\ \ &&&&&&&\cr
&&&&&&&\cr
&&\ \ p_I\swarrow&&\searrow q_I\ \ &&&&\cr
&&&&&&&&&\cr
&&\prod_{l=1}^eM^l_I\ \ \ \ &&\ \ \ \ \M_I&\buildrel\pi_I\over\to&M^{\rm naive}_I\otimes_{\O_E}\O_K\cr}
\end{eq}

\no in which both of the slanted arrows are torsors for the smooth group
scheme $\prod_{l=2}^e\Gg^{(l)}_I$. This diagram allows us to think of the
splitting model as a twisted product of the ``unramified" local models
$M^l_I$. In the next sections, we will see that the special fiber of this
diagram coincides with a certain geometric convolution diagram ([Lu],
[HN1]). By the main result of [G1] the schemes $M^l_I$ are flat over
$\Spec\O_K$. The existence of such a diagram of torsors for a smooth group
scheme therefore implies:

\begin{thm}\label{flat}
The scheme $\M_I$ is flat over $\Spec\O_K$.\endproof
\end{thm}

\medskip

\section{Local models and affine flag varieties} \label{relaflag}
\setcounter{equation}{0}

We continue with the notation of the previous sections. Recall that there
is a closed immersion of Ind-schemes $$ i: {M^{\rm naive}_I}\otimes_{\O_E}
k\longrightarrow \flag_I $$ which is described in \S \ref{aFlag}. This
immersion is equivariant for the action of $P_I\subset \GL_d(k((\Pi)))$ in
the following sense. The special fiber $\bar\Gg_I:=\Gg_I\otimes_{\O_E}k$
of the group scheme $\Gg_I$ defined in \S \ref{standGl} acts on ${M^{\rm
naive}_I}\otimes_{\O_E} k$. The isomorphism
$\Lambda_I\otimes_{\O_{F_0}}k\simeq
\ti\Lambda_I\otimes_{k[[\Pi]]}k[[\Pi]]/(\Pi^e)$ allows us to identify
$\bar\Gg_I$ with the group scheme giving the
$k[[\Pi]]/(\Pi^e)$-automorphisms of the chain
$\ti\Lambda_I\otimes_{k[[\Pi]]}k[[\Pi]]/(\Pi^e)$.  The immersion $i$ is
$P_I$-equivariant in the sense that the action of $P_I$ on $\flag_I$
stabilizes the image of $i$, the action on this image factors through the
natural group scheme homomorphism $P_I\to \bar\Gg_I$ and $i$ is
$\bar\Gg_I$-equivariant. As a result, the image of $i$ is a (finite) union
of $P_I$-orbits in $\flag_I=\GL_d(k((\Pi)))/P_I$. In fact, if $R$ is a
$k$-algebra, the $R$-rational points of the image of $i$ correspond to the
lattice chains $$ \L_0\subset \L_1\subset \cdots\subset \L_{m-1}\subset
\Pi^{-1}\L_0 $$ which fit into a diagram 
$$ 
\matrix{\ti\Lambda_{i_0,
R}&\subset&\ti\Lambda_{i_1, R}&\subset &\cdots &\subset&
\ti\Lambda_{i_{m-1}, R}&\subset&\Pi^{-1}\ti\Lambda_{i_0, R}\cr
\cup&&\cup&&&&\cup&&\cup\cr \L_0&\subset &\L_1&\subset& \cdots&\subset
&\L_{m-1}&\subset &\Pi^{-1}\L_0\cr \cup&&\cup&&&&\cup&&\cup\cr
\Pi^e\ti\Lambda_{i_0, R}&\subset&\Pi^e\ti\Lambda_{i_1, R}&\subset &\cdots
&\subset& \Pi^e\ti\Lambda_{i_{m-1}, R}&\subset&\Pi^{e-1}\ti\Lambda_{i_0,
R}\cr } 
$$
 and are such that $\L_t/\Pi^e\ti\Lambda_{i_t, R}$, and
$\ti\Lambda_{i_t,R}/\L_t$ are $R$-locally free of rank $r$, resp. $de-r$.
\medskip

Similarly, the special fiber $M^l_I\otimes_{\O_K}k'$ of the unramified local model
$M^l_I$ can be considered as a closed subscheme of the affine flag variety
$\flag_I\otimes_kk'$ via a natural closed immersion
$$
i^l: M^l_I\otimes_{\O_K}k'\ \to\ \flag_I\otimes_kk'\ .
$$
In fact, by [G1], $M^l_I\otimes_{\O_K}k'$
can be identified with the scheme-theoretic
union of a finite number of Schubert varieties in $\flag_I$
and is reduced (see [G1]). This union
is stable under the action of $P_I$.

 \medskip

Suppose now that $R$ is a $k'$-algebra and that $\{\F^j_t\}_{j,t}$ gives a
$\Spec R$-valued point of $\M_I\otimes_{\O_K}k'$. For $j=1,\ldots, e$, let
$$
 \L_{t}^j\subset \tilde\Lambda_{i_t}\otimes_{k[[\Pi]]}R[[\Pi]] 
$$
 be the inverse
image of $\F^j_t\subset \Lambda_{i_t}\otimes_{\O_{F_0}}R\simeq
\ti\Lambda_{i_t}\otimes_{k[[\Pi]]}R[[\Pi]]/(\Pi^e)$ under 
$$
\ti\Lambda_{i_t}\otimes_{k[[\Pi]]}R[[\Pi]]\to
\ti\Lambda_{i_t}\otimes_{k[[\Pi]]}R[[\Pi]]/(\Pi^e)\ . 
$$ 
We obtain a
$R[[\Pi]]$-lattice chain $\L^j_I$ $$ \L^j_0\subset \L^j_1\subset \cdots
\subset \L^j_{m-1}\subset \Pi^{-1}\L_0^j $$ which provides us with a
$\Spec R$-valued point of the affine flag variety $\flag_I$. In this way
we obtain morphisms of Ind-schemes
\begin{eq}
F^j: \M_I\otimes_{\O_K}k'\ \to\ \flag_I\otimes_kk'
\end{eq}
and
\begin{eq}
F=(F^j)_j: \M_I\ \otimes_{\O_K}k'\ \to\ \prod_{j=1}^e \flag_I\otimes_kk'\ .
\end{eq}
The morphism $F$ is a closed immersion. Actually, the $R$-rational points
of $\M_I\otimes_{\O_K}k'$ correspond to collections of lattice chains for
$j=1,\ldots, e$, $$ \L^j_0\subset \L^j_1\subset \cdots \subset
\L^j_{m-1}\subset \Pi^{-1}\L_0^j $$ which fit into a  diagram 
$$
\matrix{\ti\Lambda_{i_0, R}&\subset&\ti\Lambda_{i_1, R}&\subset &\cdots
&\subset& \ti\Lambda_{i_{m-1}, R}&\subset&\Pi^{-1}\ti\Lambda_{i_0, R}\cr
\cup&&\cup&&&&\cup&&\cup\cr \Pi^{1-e}\L^1_0&\subset
&\Pi^{1-e}\L^1_1&\subset& \cdots&\subset &\Pi^{1-e}\L^1_{m-1}&\subset
&\Pi^{-e}\L^1_0\cr \cup&&\cup&&&&\cup&&\cup\cr
\vdots&&\vdots&&&&\vdots&&\vdots\cr \cup&&\cup&&&&\cup&&\cup\cr
\L^e_0&\subset &\L^e_1&\subset& \cdots&\subset &\L^e_{m-1}&\subset
&\Pi^{-1}\L^e_0\cr \cup&&\cup&&&&\cup&&\cup\cr
\vdots&&\vdots&&&&\vdots&&\vdots\cr \cup&&\cup&&&&\cup&&\cup\cr
\L^1_0&\subset &\L^1_1&\subset& \cdots&\subset &\L^1_{m-1}&\subset
&\Pi^{-1}\L^1_0\cr \cup&&\cup&&&&\cup&&\cup\cr \Pi^e\ti\Lambda_{i_0,
R}&\subset&\Pi^e\ti\Lambda_{i_1, R}&\subset &\cdots &\subset&
\Pi^e\ti\Lambda_{i_{m-1}, R}&\subset&\Pi^{e-1}\ti\Lambda_{i_0, R},\cr } 
$$
and are such that $\L^{j}_t/\L^{j-1}_t$ for $j=2,\ldots,e$ (resp. $\Pi^{j-e}\L_t^j/\Pi^{j+1-e}\L_t^{j+1}$ for $j=1,\ldots,e-1$)
are $R$-locally free of rank
$r_{j}$ (resp. $d-r_{j+1}$), while $\L^1_t/\Pi^e\ti\Lambda_{i_t, R}$ and
$\ti\Lambda_{i_t, R}/\Pi^{1-e}\L_t^1$ are $R$-locally free of rank $r_1$,
resp. $d-r_1$.
\medskip

In what follows, for simplicity,
we will use a bar to denote the special fiber of a scheme
(or of a morphism of schemes) over $\Spec\O_K$
or over $\Spec\O_E$.

We will see that the special fiber $\ov{\M_I}$ can be naturally identified
with the geometric convolution of the reduced subschemes $\ov{M}^l_I$,
$l=1,\ldots, e$, of the affine flag variety $\flag_I\otimes_kk'$.  More
precisely, we will see below that the special fiber of the diagram
(\ref{diagram}) relates to a convolution diagram for the $P_I$-equivariant
subschemes $\ov{M}^l_I$, $l=1,\ldots, e$, defined as by Lusztig, Ginzburg
etc. ([Lu]):

\begin{eq} \label{convdiagram}
\matrix{&&& U &&&\cr
&&&&&&\cr
&&p_1\swarrow&&\searrow p_2&&&&\cr
&&&&&&\cr
&&\ov {M}^1_I\times\cdots\times \ov{M}^e_I\ &&
\ov{M}^1_I\ti\times\cdots\ti\times \ov{M}^e_I &\ \buildrel p_3\over\to\
\ov{M^{\rm naive}_I\otimes_{\O_E}{\O_K}} \subset \flag_I\otimes_kk'\ .
& \cr}
\end{eq}
\smallskip

Let us explain how the diagram (\ref{convdiagram}) is obtained (e.g [Lu]).
For simplicity of notation, we set $G=\GL_d(k((\Pi)))$ and let $\pi: G\to
\flag_I=G/P_I$ be the natural quotient morphism (of Ind-schemes). We also
set $Z_l=\ov M^l_I\subset \flag_I\otimes_kk'$ and denote by $\ti Z_l$ the
inverse image of $Z_l$ under $\pi\otimes_kk'$. Often we will omit from the
notation the base change from $k$ to $k'$; this should not cause any
confusion. Now set $$ U=\ti Z_1\times\cdots\times\ti Z_{e-1}\times
Z_e\subset G\times\cdots \times G\times G/P_I\ $$ and let $$ p_1: U\ \to\
Z_1\times\cdots\times Z_{e-1}\times Z_e $$ be the coordinate-wise
projection $(g_1,\cdots , g_{e-1},g_eP_I)\mapsto (g_1P_I,\ldots,
g_{e-1}P_I,g_eP_I)$. The morphism $p_1$ is a $(P_I)^{e-1}$-torsor for the
action given by
\begin{eq}\label{action3}
(v_1,\ldots, v_{e-1})\cdot (g_1,\ldots, g_{e-1},g_e P_I)=(g_1v_1^{-1}, \ldots, g_{e-1}v_{e-1}^{-1}, g_eP_I)\ .
\end{eq}

The convolution
$Z_1\ti\times\cdots\ti\times Z_e$ is defined as the
quotient of $U$ by the free action of $(P_I)^{e-1}$ given by
\begin{eq}\label{action4}
\ \ \ \ \ \ \ (v_1,.., v_{e-1})* (g_1,.., g_{e-1},g_e P_I)=(g_1v_1^{-1}, v_1g_2v_2^{-1},..,
v_{e-2}g_{e-1}v_{e-1}^{-1}, v_{e-1}g_e P_I)\ .
\end{eq}
We denote by
$$
p_2: U\to Z_1\ti\times\cdots\ti\times Z_e
$$
the quotient morphism.

Finally, the morphism $\ti Z_1\times\cdots\times \ti Z_{e-1}\times Z_{e}\to G/P_I$
given by $(g_1,\ldots, g_{e-1},g_eP_I)\mapsto g_1g_2\cdots g_eP_I$ factors through
the quotient to give
$$
p_3\ :\ Z_1\ti\times\cdots\ti\times Z_e\to G/P_I\ .
$$

\medskip
Let us now explain how the above convolution diagram (\ref{convdiagram})
relates to the diagram (\ref{diagram}): There is an isomorphism
$Z_1\ti\times\cdots\ti\times Z_e\simeq \ov{\M}_I$ given by $$ (g_1,\ldots,
g_e)\mapsto (\Pi^{e-1}g_1\cdot\ti\Lambda_I,
\Pi^{e-2}g_1g_2\cdot\ti\Lambda_I,\ldots, (g_1g_2\cdots g_e)\cdot
\ti\Lambda_I)\ . $$ In fact, an $R$-valued point $(g_1,\ldots, g_{e-1},
g_eP_I)$ of $\ti Z_1\times\cdots\times \ti Z_{e-1}\times Z_e$ determines a
pair consisting of a point $$ (\L^1_I, \L^2_I,\ldots,
\L^e_I)=(\Pi^{e-1}g_1\cdot\ti\Lambda_{I, R},
\Pi^{e-2}g_1g_2\cdot\ti\Lambda_{I,R} , \ldots, (g_1g_2\cdots g_e)\cdot
\ti\Lambda_{I,R}) $$ of $\ov{\M}_I$ and a collection, for $j=2,\ldots, e$,
of isomorphisms of chains $$ \sigma^j_I\ :\
\Pi^{-e+j-1}\L^{j-1}_I/\L^{j-1}_I\ \simeq\
\ti\Lambda_{I,R}/\Pi^{e-j+1}\ti\Lambda_{I,R}\simeq \Lambda^{j}_{I, R}\ .
$$ The isomorphisms $\sigma^j_I$ are given via the inverses of the maps
given by the action of $g_1\cdots g_{j-1}$ $$ \ti\Lambda_{I,R}\to
g_1\cdots g_{j-1}\cdot \ti\Lambda_{I,R}=\Pi^{-e+j-1}\L^{j-1}_I\ . $$ The
pair $((\L^j_I)_j ,\ (\sigma^j_I))_{j=2}^e$ corresponds to a point in the
special fiber $\ov{\widetilde \M}_I$.  Hence, we obtain a morphism

\begin{eq}\label{um}
u: U\to \ov{\widetilde \M}_I
\end{eq}
and (after the identification $Z_1\ti\times\cdots\ti\times Z_e=
\ov{M}^1_I\ti\times\cdots\ti\times \ov{M}^e_I\simeq \ov{\M}_I$) a diagram
\begin{eq} \label{combdiagram}
\matrix{&&&  U\ \ &&&\cr
&&& \downarrow u &&&\cr
&&&&&&\cr
&&& \ov{\widetilde \M}_I &&&\cr
&&&&&&\cr
&&\ov{p}_I\swarrow&&\searrow \ov{q}_I&&&&\cr
&&&&&&\cr
&&\ov {M}^1_I\times\cdots\times \ov{M}^e_I &&
\ov{M}^1_I\ti\times\cdots\ti\times \ov{M}^e_I & \buildrel p_3\over\to
\ov{M^{\rm naive}_{I}\otimes_{\O_E}{\O_K} }\subset \flag_I\otimes_kk'\ .
& \cr}
\end{eq}

It is easy to see that we have $p_1=u\cdot \ov{p}_I$ and $p_2=u\cdot
\ov{q}_I$.

There is a natural surjective group scheme homomorphism
\begin{eq}\label{changegps}
(P_I)^{e-1}\to  \prod_{l=2}^e{\rm Aut}_k(\ti\Lambda_I/\Pi^{e-l+1}\ti\Lambda_I)
=\prod_{l=2}^{e}\ov{\Gg}^{(l)}_I\ ;
\end{eq}
Denote its kernel by ${\cal K}$. Then the morphism $u$ realizes
$\ov{\widetilde\M_I}$ as the quotient of $U$ by the action of ${\cal
K}\subset (P_I)^{e-1}$ given by (\ref{action3}). Then, the torsor
$\ov{p}_I$ is identified with the $\prod_{l=2}^{e}\ov{\Gg}^{(l)}_I$-torsor
obtained from $p_1$ by taking the quotient by ${\cal K}$. Similarly, and
at the same time, the morphism $u$ realizes $\ov{\widetilde\M_I}$ as the
quotient of $U$ by the action of ${\cal K}\subset (P_I)^{e-1}$ given by
(\ref{action4}). Then, the torsor $\ov{q}_I$ is identified with the
$\prod_{l=2}^{e}\ov{\Gg}^{(l)}_I$-torsor obtained from $p_2$ by taking the
quotient.

\bigskip

\medskip

\section{The canonical local model for $G=\Res_{F/F_0}\GL_d$} \label{canonical}
\setcounter{equation}{0}

We continue with the assumptions and the notation of the previous
sections.
\begin{Definition}\label{deflocN}
{\rm The canonical model $M^{\rm can}_I:=M^{\rm can}(\O_F, \Lambda_I, {\bf r})$
for the group $G=\Res_{F/F_0}GL_d$, the coweight $\mu$ given by ${\bf r}$, and
the lattice chain $\Lambda_I$,
is the scheme theoretic image of the morphism
$$
\pi'_I: \M_I\to M^{\rm naive}_I\otimes_{\O_E}\O_K\to M^{\rm naive}_I\
$$
which is obtained by composing the morphism $\pi_I$ with the
base change morphism.}
\end{Definition}
\medskip
 By definition, the canonical local model $M^{\rm can}_I$ is a
closed subscheme of the naive local model $M^{\rm naive}_I$. Using
Proposition \ref{generic} we see that $M^{\rm can}_I$ and $M^{\rm
naive}_I$ have the same generic fiber. The scheme $M^{\rm can}_I$ is flat
over $\Spec\O_E$ since, by Theorem \ref{flat}, $\M_I$ is  flat over
$\Spec\O_K$. Therefore, $M^{\rm can}_I$ is the (flat) scheme theoretic
closure of the generic fiber $M_I\otimes_{\O_E}E$ in $M^{\rm naive}_I$.

 \begin{Remark} \label{extra}
{\rm a) By Theorem \ref{flat}, $\M_I$ is  flat over
$\Spec\O_K$ and hence reduced (its generic fiber being reduced).
Therefore, since $\pi'_I$ is proper, $M^{\rm can}_I$ can also be described 
as the reduced induced closed subscheme structure on the closed subset ${\rm Im}(\pi'_I)$
of the scheme $M^{\rm naive}_I$.

b) Our definition does not provide a description of $M^{\rm can}_I$ as 
a moduli scheme. On the other hand, we can observe 
that $\M_I$ and $M^{\rm naive}_I$ are moduli schemes,
the morphism $\pi'_I: \M_I\to M^{\rm naive}_I$ has a moduli description, and  
$M^{\rm can}_I$ has the following property with respect to $\pi'_I$:
It is the maximal reduced closed subscheme $Z$ of $M^{\rm naive}_I$ with the property
that, for every algebraically closed field $\Omega$, each $\Omega$-valued point of 
$Z$ lifts via $\pi'_I$ to an $\Omega$-valued
point of $\M_I$.
}\end{Remark}

In [PR] \S 8, we have defined the {\sl local model} $M^{\rm loc}_I$
for $\O_F$, $\Lambda_I$ and ${\bf r}$ as follows. For every $t\in\{0,\ldots, m-1\}$, we consider
the standard (naive) local model
$M^{\rm naive}(\Lambda_{i_t}):=M^{\rm naive}(\O_F, \Lambda_{i_t}, {\bf r})$
associated to the lattice $\Lambda_{i_t}$ (and $(F, V, \mu)$). There
is a morphism
\begin{eq}
\pi_{i_t}: M^{\rm naive}_I\to M^{\rm naive}(\Lambda_{i_t})\ ,
\end{eq}
obtained by $\{\F_i\}_{i=0}^{m-1}\mapsto \F_t$.
In [PR], we set
\begin{eq}\label{defloc}
M^{\rm loc}_I:=\bigcap_{i_t\in I}\pi^{-1}_{i_t}(M^{\rm loc}(\Lambda_{i_t}))
\end{eq}
(scheme theoretic intersection in $M^{\rm naive}_I$) where $M^{\rm
loc}(\Lambda_{i_t})\subset M^{\rm naive}(\Lambda_{i_t})$ are the (flat)
local models of EL-type which were studied in [PR]. By the above remarks,
we have 
$$
 M^{\rm loc}(\Lambda_{i_t})=M^{\rm can}_{\{i_t\}}. 
 $$
  The recent
results of G\"ortz imply now the following theorem.

\begin{thm}\label{thmintersection}
\ \

(a)  $M_I^{\rm loc}$ is flat over ${\cal O}_E$ and we have $M^{\rm
can}_I=M^{\rm loc}_I$.

(b)  The special fiber $M^{\rm can}_I\otimes_{\O_E}k$ is reduced;
its irreducible components are normal and with rational singularities.
\end{thm}

\begin{Remark} \label{conj}
{\rm (a) The flatness property in (a) above was conjectured in [PR], \S
8.

\smallskip

(b) Denote by $\mu_i$ the miniscule coweight $(1^{r_i}, 0^{d-r_i})$ of
$\GL_d$. By [G1], the special fiber $\overline M^l_I$ can be identified
with the union of Schubert cells $\bigcup_{w\in {\rm Adm}_I(\mu_l)} \O_{w}$
of the partial affine flag variety $\flag_I\otimes_kk'$. Here ${\rm
Adm}_I(\mu_l)$ denotes the $\mu_l$-admissible set inside $\tilde
W_I\backslash \tilde W/\tilde W_I$. Here $\tilde W$ denotes the extended
affine Weyl group of $\GL_d(k((\Pi)))$ and $\tilde W_I$ the subgroup of
$\tilde W$ which corresponds to the parahoric subgroup $P_I$; see [KR].
By Theorem \ref{thmintersection} (b) and the discussion in \S
\ref{relaflag}, the special fiber $\overline {M^{\rm can}_I}$ can be identified 
(up to nilpotent elements) with the
image of the convolution morphism 
$$ 
\left(\bigcup_{w\in {\rm
Adm}_I(\mu_1)} \O_{w}\right)\ti\times\cdots \ti\times\left(\bigcup_{w\in
{\rm Adm}_I(\mu_e)} \O_{w}\right)\ \to \flag_I\otimes_kk'\ . 
$$
 This image
is equal to the union $\bigcup_{w\in {\rm Adm}_I(\mu)} \O_{w}$ with
$\mu=\mu_1+\cdots +\mu_e$.

 (c) In this part of the remark we use a bar to 
 denote the special fiber of a scheme or a morphism of schemes over $\Spec\O_K$.
It follows from the definition of $M^{\rm can}_I$ that the scheme theoretic image $\pi_I(\M_I)$ is a
closed subscheme of $M^{\rm can}_I\otimes_{\O_E}\O_K$. We can easily see that 
the $\O_K$-schemes $M^{\rm can}_I\otimes_{\O_E}\O_K$ and $\pi_I(\M_I)$ have the same
generic fiber. Since $M^{\rm can}_I\otimes_{\O_E}\O_K$ is flat it follows that $\pi_I(\M_I)=M^{\rm can}_I\otimes_{\O_E}\O_K$. Similarly,
consider $\overline {\pi_I}(\overline{\M_I})\subset \overline{\pi_I(\M_I)}=\overline{M_I^{\rm can}\otimes_{\O_E}\O_K}$; these two schemes agree up to nilpotents.
By Theorem \ref{thmintersection} (b) the special fiber $\overline{M_I^{\rm can}\otimes_{\O_E}\O_K}$
is reduced (recall that the residue field $k$ is assumed perfect) and so $\overline {\pi_I}(\overline{\M_I})= \overline{M_I^{\rm can}}\otimes_k k'$.

}
\end{Remark}

\begin{Proof} Note that each morphism $\pi_{i_k}$ induces an isomorphism
between the generic fibers:
$$
\pi_{i_k}\otimes_{\O_E}E: M^{\rm naive}_I\otimes_{\O_E}E
\buildrel\sim\over\to M^{\rm naive}(\Lambda_{i_k})\otimes_{\O_E}E\ .
$$
Therefore,
$
M^{\rm loc}_I\otimes_{\O_E}E=M^{\rm naive}_I\otimes_{\O_E}E=
M^{\rm can}_I\otimes_{\O_E}E\ .
$
Since $M^{\rm can}_I$ is the scheme theoretic closure
of its generic fiber in $M^{\rm naive}_I$, we obtain
\begin{eq} \label{canincl}
M^{\rm can}_I\subset M^{\rm loc}_I
\subset M^{\rm naive}_I\ 
\end{eq}
where the inclusions are inclusions of closed subschemes. In what follows,
for simplicity, we will use a bar to denote the special fiber of a scheme
or a morphism of schemes over $\Spec\O_E$. Definition (\ref{defloc}) implies
$$ 
\overline {M^{\rm loc}_I}=\bigcap_{i_t\in
I}\bar\pi^{-1}_{i_t}(\overline{M^{\rm loc}(\Lambda_{i_t})})\ . 
$$ 
As we have
seen above, $\overline{M^{\rm naive}_I}$, resp. $\overline{M^{\rm
naive}(\Lambda_{i_t})}$, can be identified with a closed subscheme of the
affine flag, resp. affine Grassmannian, variety for $\GL_d$ over $k$. The
morphisms $\bar\pi_{i_t}$ can then be identified with the restrictions of
natural (smooth) projection morphisms from the affine flag variety to the
affine Grassmannian. By [PR], the special fibers $\overline{M^{\rm
loc}(\Lambda_{i_t})}$ are reduced and they are identified with Schubert
varieties in the affine Grassmannian; therefore the inverse images under the smooth morphisms $\pi_{i_t}$ are
also (reduced) Schubert varieties in the affine flag variety. By [G1] (see
also [F1]) all Schubert varieties in the affine flag variety are normal,
simultaneously Frobenius split, and with rational singularities; therefore
arbitrary intersections of Schubert varieties in the affine flag variety
are also reduced unions of Schubert varieties. We conclude that $\overline
{M^{\rm loc}_I}$ is reduced and that its irreducible components are normal
with rational singularities. Therefore, to show that $M^{\rm loc}_I$ is
flat and hence that $M^{\rm can}_I=M^{\rm loc}_I$, it will be enough to
show that the generic points of the irreducible components of $\overline
{M^{\rm loc}_I}$ lift to characteristic zero. This statement has recently
been shown by G\"ortz ([G3] Proposition 5.1) by using results of Haines and Ng\^o
([HN2]). Hence part (a) follows. Part (b) now follows from (a) and the
above description of the special fiber $\overline{{M}^{\rm loc}_I}$.
\endproof
\end{Proof}

\begin{Remark} \label{newproof}
{\rm As was observed by T. Haines, the use of the lifting theorem of 
G\"ortz [G3] can be avoided as follows. The proof of the second part of Remark \ref{conj}
shows that $\overline {({M}^{\rm can}_I)}_{\rm red}$ 
(the reduced induced closed subscheme structure on $\overline {{M}^{\rm can}_I}$)
is the union of Schubert varieties 
corresponding to $w$ in ${\rm Adm}_I(\mu)$. On the other hand, 
it follows from the definition of $M^{\rm loc}_I$,
that $\overline {{M}^{\rm loc}_I}$ (which is already known to be reduced by the first part
of the proof of Theorem \ref{thmintersection})
is the union of Schubert varieties corresponding to $w$ in ${\rm Perm}_I(\mu)$.
Here ${\rm Perm}_I(\mu)$ denotes the $\mu$-permissible set which however
has been shown to be identical with ${\rm Adm}_I(\mu)$
(Haines and Ng\^o [HN2] for $I=\{0,\ldots, d-1\}$,
G\"ortz [G3] in the remaining cases). 
The closed immersions
$$
\overline {({M}^{\rm can}_I)}_{\rm red}\subset  \overline {{M}^{\rm can}_I}\subset
\overline {{M}^{\rm loc}_I}
$$
are thus all isomorphisms, and so $\overline {{M}^{\rm can}_I}=
\overline {{M}^{\rm loc}_I}$ and $\overline {{M}^{\rm can}_I}$ is reduced.
The flatness of $M^{\rm can}_I$ now implies that $M^{\rm can}_I=M^{\rm loc}_I$
and the rest follows.

}
\end{Remark}

 \part{ }

\section{The ``naive" local models for $G=\Res_{F/F_0}{\rm GSp}_{2g}$} \label{naivesymlp}
\setcounter{equation}{0}

We continue with the notation of \S \ref{notations}. Let $(V, \{\ ,\ \})$
be the standard symplectic vector space over $F$ of dimension $2g$ with
basis $e_1,\ldots, e_g,f_1,\ldots, f_g$, i.e
\begin{eq}
\{e_i,e_j\}=\{f_i,f_j\}=0,\quad \{e_i,f_j\}=\delta_{ij}\ .
\end{eq}
Let $<v, w>=\Tr_{F/F_0}\{v,w\}$. Then, since $F/F_0$ is separable,
$<\ ,\ >$ is a non-degenerate alternating form on $V$ with values in $F_0$
which, for all $a\in F$, satisfies
\begin{eq}
<av, w>=<v, aw>\ .
\end{eq}
If $\Lambda$ is an $\O_F$-lattice in $V$, we set
$$
\Lambda^*:=\{v\in V\ |\ \{v, \lambda\}\in\O_F, \ \hbox{\rm for all $\lambda\in \Lambda$}\}\ ,
$$
$$
\hat\Lambda:=\{v\in V\ | <v, \lambda>\in\O_{F_0}, \ \hbox{\rm for all $\lambda\in \Lambda$}\}\ ,
$$
for the dual (``complementary") $\O_F$-lattices with respect
to the forms $\{\ ,\ \}$ and $<\ ,\ >$ respectively.

Now let $\delta$ be an $\O_F$-generator of the inverse
different ${\cal D}^{-1}_{F/F_0}$ (if $F/F_0$ is tamely ramified,
we can take $\delta=\pi^{1-e}$). Set
$$
\Lambda_0={\rm Span}_{\O_F}\{e_1,\ldots, e_g,\delta f_1,\ldots, \delta f_g\}\subset V .
$$
Then the $\O_F$-lattice $\Lambda_0$ is self-dual with respect
to the form $<\ ,\ >$, i.e
\begin{eq}
\hat\Lambda_0=\Lambda_0\ .
\end{eq}
Indeed, $<e_i, a\delta f_j>=0$ if $i\neq j$, and $<e_i, a\delta
f_i>=\Tr_{F/F_0}(a\delta)$; this is in $\O_{F_0}$ exactly when $a$ is in
$\O_F$.

For $0\leq r\leq g$, let
$$
\Lambda_r=
{\rm Span}_{\O_F}\{\pi^{-1} e_1,\ldots, \pi^{-1} e_r, e_{r+1}, \ldots, e_g, \delta f_1,\ldots ,\delta f_g\}\ .
$$
We obtain a chain of inclusions of $\O_F$-lattices
\begin{eq}
\pi\Lambda_0\subset \hat\Lambda_r\subset \hat\Lambda_0=\Lambda_0\subset \Lambda_r\subset \pi^{-1}\Lambda_0\ .
\end{eq}
In fact, we have
$$
\hat\Lambda_r={\rm Span}_{\O_F}\{e_1,\ldots ,e_g, \pi\delta f_1,\ldots ,\pi\delta f_r,
\delta f_{r+1},\ldots , \delta f_{g}\}\ .
$$
We can extend $\Lambda_0\subset \Lambda_1\subset \cdots \subset \Lambda_g$
to a complete $\O_F$-lattice chain $\{\Lambda_i\}_{i\in\Z}$ in $V$ by
setting
$$
\Lambda_i=\pi^{-t}\Lambda_j, \ \hbox{\rm for $i=2gt+j$, $0\leq j\leq g$},
$$
$$
\Lambda_{i}=\pi^{-t}\hat\Lambda_{-j}, \ \hbox{\rm for $i=2gt+j$, $-g\leq j< 0$},
$$
The essential part of this $\O_F$-lattice chain is
$$
\hat\Lambda_{g}\subset \cdots\subset\hat\Lambda_1\subset \hat\Lambda_0=\Lambda_0\subset\Lambda_1\subset\cdots \subset \Lambda_g=\pi^{-1}\hat\Lambda_g\ .
$$
The lattice chain $\{\Lambda_i\}_{i\in\Z}$ is ``self-dual" (for every $i$ there
is a $j$ such that $\hat\Lambda_i=\Lambda_j$, in fact we have $\hat\Lambda_i=\Lambda_{-i}$) and ``complete"
(for every $i$, ${\rm dim}_k(\Lambda_{i+1}/\Lambda_i)=1$).
We will sometimes write

\begin{eq}\label{forms}
<\ ,\ >_{\pm i}\ :\ \Lambda_{\pm i}\times \Lambda_{\mp i}\ \to\ \O_{F_0}
\end{eq}
for the corresponding perfect form. These sets of forms
are alternating in the sense that
$$
<v ,w >_{\pm i}=-<w, v>_{\mp i}\ .
$$

Now fix a subset $I=\{i_0<\cdots <i_{m-1}\}\subset \{0,1,\ldots, g\}$
and consider the self-dual periodic $\O_F$-lattice chain $\Lambda_I$ given by taking
all lattices of the form $\pi^n\Lambda_{i_k}$, $\pi^n\hat\Lambda_{i_t}$ for
$n\in \Z$, $t=0,\ldots, m-1$. An essential part of the lattice chain $\Lambda_I$
is
$$
\hat\Lambda_{i_{m-1}}\subset \cdots \subset \hat\Lambda_{i_0}\subset \Lambda_{i_0}\subset \cdots
\subset \Lambda_{i_{m-1}}\subset \pi^{-1}\hat\Lambda_{i_{m-1}}\ .
$$

The standard (``naive") local model $N^{\rm naive}_I$ associated
by Rapoport-Zink [RZ], Definition 3.27 to the reductive group $G=\Res_{F/F_0}{\rm GSp}(V, <\ ,\ >)$,
the cocharacter $\mu$ given by $\{(1^g, 0^g)\}_\phi$
and the parahoric subgroup which is the stabilizer of the
$\O_F$-lattice chain $\Lambda_I$, is by definition the $\O_{F_0}$-scheme representing the following functor
on $({\rm Schemes}/\O_{F_0})$:

For every $\O_{F_0}$-scheme $S$, $N^{\rm naive}_I(S)$ is the set of
collections $\{\F_{i_t},\F_{-i_t}\}_{t=0,\ldots, m-1}$  of
$\O_F\otimes_{\O_{F_0}}\O_S$-submodules of $\Lambda_{i_t, S}$, resp.\
$\hat\Lambda_{i_t, S}$ which fit into a commutative diagram $$
\matrix{\hat\Lambda_{i_{m-1}, S}&\to&\cdots&\to& \hat\Lambda_{i_0,
S}&\to&\Lambda_{i_0, S}&\to& \cdots& \to& \Lambda_{i_{m-1}, S}&\buildrel
\pi\over \to & \hat\Lambda_{i_{m-1}, S}\cr \cup&& &&\cup&&
\cup&&&&\cup&&\cup\cr \F_{-i_{m-1}}&\to &\cdots &\to &\F_{-i_0}&\to&
\F_{i_0}&\to&\cdots &\to&\F_{i_{m-1}}&\to &\F_{-i_{m-1}}\cr} $$ and are
such that:
\smallskip

a) $\F_{i_t}$, resp. $\F_{-i_t}$, is   Zariski locally on $S$ a $\O_S$-direct summand
of $\Lambda_{i_t, S}$, resp. $\hat \Lambda_{i_{t}, S}$, of rank $eg$,
\smallskip

b) the compositions
$
\F_{-i_t}\subset \hat\Lambda_{i_t, S}\to \hat\F_{i_t}$\ , $\F_{i_t}\subset
\Lambda_{i_t, S}=\hat{\hat\Lambda}_{i_t,S}\to \hat\F_{-i_t}$, where
$\hat\F_{\pm i_t}=Hom_{\O_S}(\F_{\pm i_t}, \O_S)$ and the second maps are
the duals of the inclusions $\F_{i_t}\subset \Lambda_{i_{t}, S}$, resp.\
$\F_{-i_t}\subset \hat\Lambda_{i_t, S}$, are the zero maps.
\smallskip

c) For every $a\in\O_F$, and $t=0,\ldots, m-1$, we have
$$
\det(a\ |\ \F_{\pm i_t})=\prod_{\phi}\phi(a)^{g}
$$
where again this identity is meant as an identity of polynomial
functions on $\O_F$.
\medskip

\begin{Remark} \label{sect}
{\rm For $\F\subset \Lambda_{i_t,S}$, we set
$\F^\perp:=\ker(\hat\Lambda_{i_t, S}\to \hat\F)\subset \hat\Lambda_{i_t,
S}$. For $\Gg\subset \hat\Lambda_{i_t,S}$ we set
$\Gg^\perp:=\ker(\Lambda_{i_t,S}=\hat{\hat\Lambda}_{i_t,S}\to
\hat\Gg)\subset \Lambda_{i_t, S}$. If $\F$, resp.\ $\Gg$ are locally
$\O_S$-direct summands of $\Lambda_{i_t, S}$, resp.\ $\hat\Lambda_{i_t,
S}$, then $\F^\perp$, resp.\ $\Gg^\perp$ are locally $\O_S$-direct
summands of $\hat\Lambda_{i_t, S}$, resp.\ $\Lambda_{i_t, S}$. Condition
(b) implies that $$ \F_{-i_t}\subset {\rm ker}(\hat\Lambda_{i_t, S}\to
\hat\F_{i_t})=(\F_{i_t})^\perp\ ,\quad \F_{i_t}\subset {\rm
ker}(\Lambda_{i_t, S}\to \hat\F_{-i_t})=(\F_{-i_t})^\perp\ . $$ Since by
(a), $\F_{\pm i_t}$, $(\F_{\pm i_t})^\perp$ all have rank $eg$, we obtain
$\F_{-i_t}=(\F_{i_t})^\perp$, $\F_{i_t}=(\F_{-i_t})^\perp$.
\smallskip

Hence, $N^{\rm naive}_I(S)$ is in bijection with the set of collections
$\{\F_t\}_t$ of $\O_F\otimes_{\O_{F_0}}\O_S$-submodules $\F_t\subset
\Lambda_{i_t,S}$, which are, Zariski locally on $S$, $\O_S$-direct
summands of $\Lambda_{i_t, S}$ of rank $eg$ and which satisfy:
\smallskip

i) For every $a\in\O_F$, $t=0,\cdots, m-1$,
$$
\det(a\ |\ \F_{t})=\prod_{\phi}\phi(a)^{g}
$$
(as always this identity is meant as an identity of polynomial
functions on $\O_F$),
\smallskip

ii) The inclusions $\F_t\subset\Lambda_{i_t, S}$, $\F^\perp_t\subset\hat\Lambda_{i_t, S}$
fit into a commutative diagram
$$
\matrix{\hat\Lambda_{i_{m-1}, S}&\to&\cdots&\to& \hat\Lambda_{i_0, S}&\to&\Lambda_{i_0, S}&\to& \cdots&
\to& \Lambda_{i_{m-1}, S}&\buildrel \pi\over \to & \hat\Lambda_{i_{m-1}, S}\cr
\cup&& &&\cup&& \cup&&&&\cup&&\cup\cr
\F_{m-1}^\perp&\to &\cdots &\to &\F_{0}^\perp&\to& \F_{0}&\to&\cdots &\to&\F_{m-1}&\to &\F_{m-1}^\perp\ .\cr}
$$}
\end{Remark}

\section{The splitting model for $G=\Res_{F/F_0}{\rm GSp}_{2g}$} \label{splitSp}
\setcounter{equation}{0}

We continue with the notation of the previous section. Consider the
functor $\N_I=\N(\O_F, \Lambda_I, g)$ on (Schemes/$\Spec\O_K$) which to a
$\O_K$-scheme $S$ associates the set of collections $\{\F^j_{i_t},
\F^j_{-i_t}\}_{t,j=1,\ldots, e}$ of
$\O_F\otimes_{\O_{F_0}}\O_S$-submodules of $\Lambda_{i_t, S}$, resp.\
$\hat\Lambda_{i_t, S}$ which fit into a commutative diagram $$
\matrix{\hat\Lambda_{i_{m-1}, S}&\to&\cdots&\to& \hat\Lambda_{i_0,
S}&\to&\Lambda_{i_0, S}&\to& \cdots& \to& \Lambda_{i_{m-1}, S}&\buildrel
\pi\over \to & \hat\Lambda_{i_{m-1}, S}\cr \cup&& &&\cup&&
\cup&&&&\cup&&\cup\cr \F^e_{-i_{m-1}}&\to &\cdots &\to &\F^e_{-i_0}&\to&
\F^e_{i_0}&\to&\cdots &\to&\F^e_{i_{m-1}}&\to &\F^e_{-i_{m-1}}\cr \cup&&
&&\cup&& \cup&&&&\cup&&\cup\cr \F^{e-1}_{-i_{m-1}}&\to &\cdots &\to
&\F^{e-1}_{-i_0}&\to& \F^{e-1}_{i_0}&\to&\cdots
&\to&\F^{e-1}_{i_{m-1}}&\to &\F^{e-1}_{-i_{m-1}}\cr \cup&& &&\cup&&
\cup&&&&\cup&&\cup\cr \vdots&&&&\vdots&&\vdots&&&&\vdots&&\vdots\cr \cup&&
&&\cup&& \cup&&&&\cup&&\cup\cr \F^1_{-i_{m-1}}&\to &\cdots &\to
&\F^1_{-i_0}&\to& \F^1_{i_0}&\to&\cdots &\to&\F^1_{i_{m-1}}&\to
&\F^1_{-i_{m-1}}\cr} $$ and are such that:

a) $\F^j_{i_t}$, resp. $\F^j_{-i_t}$, is  Zariski locally on $S$ a $\O_S$-direct summand
of $\Lambda_{i_t, S}$, resp. $\hat \Lambda_{i_{t}, S}$, of rank $jg$
and satisfies, for all $a\in \O_F$,
$$
(a\otimes 1-1\otimes\phi_j(a))\F^j_{\pm i_t}\subset \F^{j-1}_{\pm i_t}\ .
$$
\smallskip

b) the compositions $$ \F^j_{-i_t}\subset \hat\Lambda_{i_t, S}\to
\hat\F^j_{i_t}\ , \ \ \ \F^j_{i_t}\subset \Lambda_{i_t,
S}=\hat{\hat\Lambda}_{i_t,S}\to \hat\F^j_{-i_t}\ , $$ are the zero maps.
\medskip

By Remark \ref{sect}, we see that the above conditions
imply that for every $t$, $\F^e_{-i_t}=(\F^e_{i_t})^\perp$,
$\F^e_{i_t}= (\F^e_{-i_t})^\perp$. Hence,
we obtain chains
$$
(0)\subset \F^1_{-i_t}\subset \cdots\subset \F^e_{-i_t}= (\F^e_{i_t})^\perp\subset
\cdots\subset (\F^1_{i_t})^\perp\subset \hat\Lambda_{i_t,S}\ ,
$$
$$
(0)\subset \F^1_{i_t}\subset \cdots\subset \F^e_{i_t}= (\F^e_{-i_t})^\perp\subset
\cdots\subset (\F^1_{-i_t})^\perp\subset \Lambda_{i_t,S}\ .
$$
\smallskip

c) In addition to (a) and (b), we require that, for all $j=1,\ldots, e-1$,
and every $a\in \O_F$,
\smallskip

$$
\prod_{j+1\leq q\leq e}(a\otimes 1-1\otimes \phi_q(a))\ (\F^j_{-i_t})^\perp\subset
\F^{j}_{i_t}\ ,
$$
$$
\prod_{j+1\leq q\leq e}(a\otimes 1-1\otimes \phi_q(a))\ (\F^j_{i_t})^\perp\subset
\F^{j}_{-i_t}\ .
$$
 \smallskip

Obviously the functor $\N_I$ is represented by a projective scheme over
$\Spec\O_K$ which we will also denote by $\N_I$. As in the case of
$G=\Res_{F/F_0}\GL_d$, there is a projective morphism

\begin{eq}
\pi_I: \N_I\to N^{\rm naive}_I\otimes_{\O_{F_0}}\O_K
\end{eq}
given by $\{\F^j_{\pm i_t}\}_{j,t}\mapsto \{\F^e_{\pm i_t}\}_k$.
A construction similar to the one in the proof of
Proposition \ref{generic} shows that, also in this case,
$\pi_I$ induces an isomorphism
\begin{eq}
\pi_I\otimes_{\O_{K}}K: \N_I\otimes_{\O_K}K\buildrel \sim\over\longrightarrow
N^{\rm naive}_I\otimes_{\O_{F_0}}K
\end{eq}
on the generic fibers.

\bigskip

Now suppose that $\{\F^j_{\pm i_t}\}_{j,t}$ is
an $S$-valued point of $\N_I$.
As in the case of $G={\rm Res}_{F/F_0}\GL_d$,
for $l=1, \ldots ,e$, let us set
$$
\Psi^l_{\pm i_t, S}=\ker(Q^{l}(\pi)\ |\ \Lambda_{\pm i_t, S}/\F^{l-1}_{\pm i_t})\ ;
$$
this is an $\O^{(l)}_K\otimes_{\O_K}\O_S$-module. We also set
$$
\Upsilon^l_{\pm i_t, S}:={\rm ker}(\pi-a_l\ |\ \Lambda_{i_t, S}/{\F^{l-1}_{\pm i_t}})=
{\rm ker}(\pi-a_l\ |\ \Psi^l_{\pm i_t, S})\ .
$$

The proof of Proposition \ref{bundles} implies that
the $\O^{(l)}_K\otimes_{\O_K}\O_S$-module $\Psi^l_{\pm i_t, S}$
is, locally on $S$, free of rank $2g$ and that its formation
commutes with base change in $S$. Similarly, $\Upsilon^l_{\pm i_t, S}$
is a locally free $\O_S$-module of rank $2g$ whose formation
commutes with base change in $S$.

\begin{lemma} \label{perp}
Suppose that $\{\F^j_{\pm i_t}\}_{j,t}$ is an $S$-valued point of $\N_I$. Then for
$l=1,\ldots, e$, $k=0,\ldots, m-1$, we have
$$
Q^l(\pi)^{-1}(\F^{l-1}_{\pm i_t})=(\F^{l-1}_{\mp i_t})^\perp
$$
where the left hand side is the inverse image of the submodule
$\F^{l-1}_{\pm i_t}\subset \Lambda_{\pm i_t, S}$
under $\Lambda_{\pm i_t,S}\to \Lambda_{\pm i_t, S}$ given
by multiplication by $Q^l(\pi)$.
\end{lemma}

\begin{Proof}
The proof of Proposition \ref{bundles} (b) shows that
$Q^l(\pi)^{-1}(\F^l_{\pm i_t})\subset \Lambda_{\pm i_t, S}$ is locally an
$\O_S$-direct summand of rank $g(2e-l+1)$. Observe that the condition (c)
in the definition of the splitting model $\N_I$ translates to $$
(\F^{l-1}_{\mp i_t})^\perp\subset Q^l(\pi)^{-1}(\F^{l-1}_{\pm i_t})\ . $$
Now $(\F^{l-1}_{\mp i_t})^\perp$ and $Q^l(\pi)^{-1}(\F^{l-1}_{\pm i_t})$
have the same $\O_S$-rank and they are both locally $\O_S$-direct summands
of $\Lambda_{\pm i_t, S}$. Hence, they are equal.\endproof
\end{Proof}
\medskip

Suppose that $\{\F^j_{\pm i_t}\}_{j,t}$ is an $S$-valued point of $\N_I$.
Lemma \ref{perp} implies that
\begin{eq}\label{psiperp}
\Psi^l_{\pm i_t, S}=(\F^{l-1}_{\pm i_t, S})^\perp/\F^{l-1}_{\pm i_t, S}\ .
\end{eq}

Therefore, there are perfect $\O_S$-bilinear  forms
\begin{eq} \label{forms}
<\ ,\ >^l_{\pm i_t}\ :\ \Psi^l_{\pm i_t, S}\times \Psi^l_{\mp i_t, S}\ \to\ \O_S
\end{eq}
induced by the forms $<\ ,\ >_{\pm i_t}$. These satisfy $$
<v, w>^l_{\pm i_t}=-<w, v>^l_{\mp i_t}, \ \hbox{\ \rm and\ }<av, w>^l_{\pm i_t}=<v, aw>^l_{\pm i_t}\ ,
$$ for all $a\in \O^{(l)}_K$ (i.e the pairings respect the action of
$\O^{(l)}_K$).
\smallskip

For $l=1,\ldots,e$, consider the
chain of free $\O_K^{(l)}$-modules $\Lambda^l_I$
obtained from the free $\O_F\otimes_{\O_{F_0}}\O_K$-module
chain $\Lambda_{I,\O_K}:=\Lambda_I\otimes_{\O_{F_0}}\O_K$ by extending scalars via
$$
\phi^l: \O_F\otimes_{\O_{F_0}}\O_K\ \to\ \O_K[T]/(Q^l(T))=\O_K^{(l)}\ .
$$
We will define perfect $\O_{K}$-bilinear alternating pairings
$$
<\ ,\ >^l_{\pm i_t}\ :\ \Lambda^l_{\pm i_t}\times \Lambda^l_{\mp i_t}\ \to\ \O_{K}
$$ which respect the action of $\O^{(l)}_K$ as follows: Using
(\ref{kerim1}), we see that there are canonical isomorphisms
\begin{eq}\label{caniso2}
\Lambda^l_{\pm i_t}\ \simeq\ {\rm Im}(\ Q_l(T)
\ |\ \Lambda_{\pm i_t, \O_K} )=\ker(Q^l(T)\ |\ \Lambda_{\pm i_t,\O_K} )\ .
\end{eq}
Suppose that $v\in \Lambda^l_{\pm i_t}$, $w\in \Lambda^l_{\mp i_t}$. Via
(\ref{caniso2}) we can identify $v$ with an element of $\ker(Q^l(T)\ |\
\Lambda_{\pm i_t,\O_K} )\subset \Lambda_{\pm i_t,\O_K}$ and choose $\ti
w\in \Lambda_{\mp i_t,\O_K} $ such that $$ Q_l(T)\cdot \tilde w=w\ . $$ We
set

\begin{eq}\label{deriv}
<v,w>^l_{\pm i_t}=<v, \ti w>_{\pm i_t}\ .
\end{eq}
\smallskip

It is easy to see that this is independent of the choice of $\ti w$. It
provides us with perfect $\O_{K}$-bilinear forms which respect the action
of $\O^{(l)}_K$ and satisfy $$
<w, v>^l_{\mp i_t}=-<v, w>^l_{\pm i_t}\
$$
(i.e they are alternating).

Let us set $V^l$ for the $K[T]/(Q^l(T))$-module obtained from the
$F\otimes_{F_0}K$-module $V\otimes_{F_0}K$ by extending scalars via
$\phi^l\otimes_{\O_K}K\ :\ F\otimes_{F_0}K\to K[T]/(Q^l(T))$. Then, for
all $t=0,\ldots, m-1$, $\Lambda^l_{\pm i_k}\subset V^l$ and the pairings
$<\ ,\ >^l_{\pm i_t}$ are all restrictions of a single perfect
$K$-bilinear alternating pairing $$
<\ ,\ >^l\ :\ V^l\times V^l\ \to \ K
$$ which respects the action of $\O^{(l)}_K\otimes_{\O_K}K=K[T]/(Q^l(T))$.
It is easy to see that, under this form, $\Lambda^l_{-i_t}$ is dual
(``complementary") to $\Lambda^l_{i_t}$. In this sense, the chain
$\Lambda^l_I$ is a periodic self-dual chain of free $\O^{(l)}_K$-modules
in $V^l$.

\medskip

Consider the chain of $\O^{(l)}_K\otimes_{\O_K}\O_S$-modules $\Psi^l_{I, S}$:
\begin{eq}
\ \ \ \ \cdots\to\Psi^l_{-i_{m-1}, S}\to\cdots \to \Psi^l_{-i_{0}, S}\to
\Psi^l_{i_0, S}\to\cdots \to \Psi^l_{i_{m-1}, S}\buildrel \pi\over
\to \Psi^l_{-i_{m-1}, S}\to\cdots\
\end{eq}
over $S$ with the morphisms induced by the commutative diagram in the
definition of $\N_I$, and with the bilinear forms (\ref{forms}).

\begin{prop} \label{bundlessympl1} \  a) The pairings (\ref{forms})
provide the chain $\Psi^l_{I, S}$ with the structure of a polarized chain
of $\O^{(l)}_K\otimes_{\O_K}\O_S$-modules $\Psi^l_{I, S}$ of type
$(\Lambda^l_I)$ in the sense of [RZ] Def. 3.14, p. 75 ($\O^{(l)}_K$ is
not a maximal order in $\O^{(l)}_K\otimes_{\O_K}K$, however the definition
still makes sense).

b) Zariski locally on $S$, the polarized chain of
$\O^{(l)}_K\otimes_{\O_K}\O_{S}$-modules
$\Psi^l_{I, S}$ is (symplectically) isomorphic to the
polarized chain of
$\O^{(l)}_K\otimes_{\O_K}\O_{S}$-modules $\Lambda^l_{I, S}:=\Lambda^l_I\otimes_{\O_K}\O_{S}$
which is obtained from the $\O^{(l)}_K$-chain $\Lambda^l_I$.
\end{prop}

\begin{Proof}
To show (a) we have to show that the chain $\Psi^l_{I,S}$ satisfies the
conditions of [RZ], Def. 3.14 p. 75 (see also Def. 3.6 and Cor. 3.7).
Assuming (a), part (b) of the proposition follows from a simple extension
of [RZ], Prop. A 21 to the case at hand.

Now the only condition in loc. cit. that does not follow immediately from
the definitions is the requirement (corresponding to condition (2) of Def.
3.6) that Zariski locally on $S$ the quotient of two successive modules in
the chain $\Psi^l_{I,S}$ is $\O^{(l)}_K\otimes_{\O_K}\O_S$-isomorphic to
the quotient of the two corresponding successive modules of the chain
$\Lambda^l_{I,S}$. This can be shown exactly as the corresponding
statement in the proof of Proposition \ref{bundles}.\endproof
\end{Proof}

\medskip

\medskip

For $l=1,\ldots, e$, there is a natural isomorphism $$ V\otimes_{F,\phi_l}
K\simeq {\rm Im}(Q^{l+1}(T)\ |\ V^{l})=\ker(T-a_l\ |\ V^{l})\ . $$ A
construction analogous to (\ref{deriv}) allows us to define a perfect
$K$-bilinear alternating form $$
<\ ,\ >_l\ :\ V\otimes_{F,\phi_l} K\times V\otimes_{F,\phi_l} K\ \to\ K\ .
$$
Now set
\begin{eq} \label{xisympl}
\ \ \ \ \ \ \Xi^l_{\pm i_k}=\Lambda_{\pm i_k}\otimes_{\O_F,\phi_l}\O_K \ \simeq\ {\rm Im}(Q^{l+1}(T)\ |\ \Lambda^l_{\pm i_k}
 )
=\ker(T-a_l\ |\ \Lambda^l_{\pm i_k})\ .
\end{eq}
Once again, we can see that
\begin{eq}
<\ ,\ >_{l, \pm i_t}\ :\ \Xi^l_{\pm i_t}\times \Xi^l_{\mp i_t}\to \O_K
\end{eq}
defined by restricting $<\ ,\ >_l$ to the lattices $\Xi_{\pm i_t}$, $\Xi_{\mp i_t}$
give a system of perfect $\O_K$-bilinear alternating
forms. By construction, we have
\begin{eq}\label{constr}
<v, w>_{l, \pm i_t}=<v, \ti w>^l_{\pm i_t},
\end{eq}
where we regard $v$ as an element of $\ker(T-a_l\ |\ \Lambda^l_{\pm i_t})$
and where $\ti w\in \Lambda^l_{\mp i_t}$ satisfies $Q^{l+1}(T)\cdot \ti
w=w$.

Hence, for each $l=1,\ldots, e$, we obtain a self-dual $\O_K$-lattice chain
$\Xi^l_I$ in the $K$-vector space $V\otimes_{F,\phi_l}K$ by using the lattices
$\Xi^l_{\pm i_t}=\Lambda_{\pm i_t}\otimes_{\O_F, \phi_l}\O_K$.
The essential part of this chain is:
$$
\Xi^l_{-i_{m-1}}\subset \cdots \subset \Xi^l_{-i_0}\subset \Xi^l_{i_0}\subset \cdots
\subset \Xi^l_{i_{m-1}}\subset a_l^{-1}\Xi^l_{i_{m-1}}\ .
$$
\medskip

Now denote by $N^l_I$ the ``unramified" local model $N^l_I:=N(\O_K,
\Xi^l_I, <\ ,\ >)$ over $\Spec\O_K$ defined in [RZ] for $G={\rm
GSp}(V\otimes_{F,\phi_l}K, <,>_l)$ (a group over $K$), the cocharacter
$\mu$ given by $(1^g, 0^g)$ and the self-dual lattice chain $\Xi^l_I$. By
definition, $N^l_I$ is the projective scheme over $\Spec\O_K$ which
classifies collections $\{\F_{\pm i_t}\}_t$ of $\O_S$-submodules $\F_{\pm
i_t}\subset \Xi^l_{\pm i_t, S}:=\Xi^l_{\pm i_t}\otimes_{\O_K}\O_S$ which
fit into a commutative diagram $$ \matrix{\Xi^l_{-i_{m-1}, S}&\to & \cdots
&\to& \Xi^l_{-i_0, S}&\to & \Xi^l_{i_0, S}& \to &\cdots& \to &
\Xi^l_{i_{m-1}}&\buildrel a_l\over \to &\Xi^l_{i_{m-1}}\cr
\cup&&&&\cup&&\cup&&&&\cup&&\cup\cr
\F_{-i_{m-1}}&\to&\cdots&\to&\F_{-i_0}&\to&\F_{i_0}&\to&\cdots&\to&\F_{i_{m-1}}&
\to&\F_{-i_{m-1}}\ \cr} $$ where $\F_{\pm i_t}$ are Zariski locally
$\O_S$-direct summands of $\Xi^l_{i_t,S}$ of rank $g$ and which satisfy $$
\F_{-i_t}=\F_{i_t}^\perp,\quad \F_{i_t}=\F_{-i_t}^\perp\ . $$
\smallskip

Now let us denote by ${\cal H}^{(l)}_I$, resp. ${\cal H}^l_I$, the group scheme
over $\Spec\O_K$ whose $S$-points are the $\O^{(l)}_K\otimes_{\O_K}\O_S$-module,
resp. $\O_S$-module, automorphisms
of the polarized chain  $\Lambda^{l}_I\otimes_{\O_K}\O_S$, resp.   $\Xi^l_I\otimes_{\O_K}\O_S$,
which respect the forms $<\ ,\ >^l_{\pm i_t}$, resp. $<\ ,\ >_{l, \pm i_t}$,
up to a similitude which is the same for all indices $t$. These groups
are extensions of the multiplicative group by the group scheme of symplectic
$\O^{(l)}_K\otimes_{\O_K}\O_S$-module,
resp. $\O_S$-module, automorphisms
of the polarized chains
$\Lambda^{l}_I\otimes_{\O_K}\O_S$, resp.  $\Xi^l_I\otimes_{\O_K}\O_S$.
An argument as in the proof of [RZ] Prop. A.21 shows that the latter
group schemes are smooth  over $\Spec\O_K$.
Therefore, ${\cal H}^{(l)}_I$ and ${\cal H}^l_I$ are also smooth group schemes over $\Spec\O_K$.

\medskip

Now for an $S$-valued point of $\N_I$ given by $\{\F^j_{\pm i_t}\}_{j,t}$
and $l=1,\ldots, e$, $k=0,\ldots, m-1$, we set $$ \Upsilon^l_{\pm i_t,
S}=\ker(\pi-a_l\ |\ \Psi^l_{\pm i_t, S})=\ker(\pi-a_l\ |\ (\F^{l-1}_{\mp
i_t, S})^\perp/\F^{l-1}_{\pm i_t, S})\ . $$ Notice that there is a
canonical $\O_S$-homomorphism $$ {\rm Im}(Q^{l+1}(\pi)\ |\ \Psi^{l}_{\pm
i_t, S})\to \ker(\pi-a_l\ |\ \Psi^l_{\pm i_t, S})=\Upsilon^l_{\pm i_t, S}\
. $$ It follows from Proposition \ref{bundlessympl1} and (\ref{xisympl})
that this is an isomorphism.

A construction as in (\ref{constr}) now allows us to use this isomorphism
and the forms $<\ ,\ >^l_{\pm i_t}:\Psi^l_{\pm i_t, S}\times \Psi^l_{\mp
i_t, S}\ \to\ \O_S\ $ to derive $\O_S$-bilinear alternating forms $$
<\ ,\ >_{l, \pm i_t}:\Upsilon^l_{\pm i_t, S}\times \Upsilon^l_{\mp i_t, S}\ \to\ \O_S\ .
$$

\begin{prop} The $\O_S$-modules $\Upsilon^l_{\pm i_t, S}$ are locally free
of rank $2g$ and the $\O_S$-bilinear alternating
forms
$$
<\ ,\ >_{l,\pm i_t}\ :\ \Upsilon^l_{\pm i_t, S}\times \Upsilon^l_{\mp i_t, S}\ \to\ \O_S\ .
$$
are perfect.
Furthermore, the resulting chain of $\O_S$-modules
over $S$
\begin{eq}
\ \ \ \cdots\to\Upsilon^l_{-i_{m-1}, S}\to\cdots \to \Upsilon^l_{-i_{0}, S}\to \Upsilon^l_{i_0, S}\to\cdots \to \Upsilon^l_{i_{m-1}, S}\buildrel a_l\over
\to \Upsilon^l_{-i_{m-1}, S}\to\cdots
\end{eq}
is a polarized chain of type $(\Xi^l_I)$ and is, Zariski locally on $S$, (symplectically) isomorphic to
the polarized chain of $\O_S$-modules $\Xi^l_I\otimes_{\O_K}\O_S$
obtained from the self-dual lattice chain $\Xi^l_I$.
\end{prop}

\begin{Proof} This follows from Proposition
\ref{bundlessympl1} and the above discussion.\endproof
\end{Proof}

 \bigskip

Let $\widetilde {\N_I}$ denote the scheme over $\Spec\O_K$ whose
$S$-points correspond to pairs $$ \widetilde {\N_I}(S):=\left(\{\F^j_{\pm
i_t}\}_{j, t}\ ,\ \{\sigma^j_I\}\right)_{j=2}^e $$ where $\{\F^j_{\pm
i_t}\}_{j,t}$ is an $S$-valued point of $\N_I$ and for $l=2,\ldots, e$, $$
\sigma^l_I\ :\ \Psi^l_{I,S}\ \buildrel \sim\over\to\ \Lambda^l_{I, S} $$
is a symplectic (up to similitude) isomorphism of polarized
$\O^{(l)}_K\otimes_{\O_K}\O_S$-chains. The natural projection morphism
$q_I: \widetilde {\N_I}\to \N_I$ is a torsor for $\prod_{l=2}^e{\cal
H}^{(l)}_I$.

\medskip

Notice that an isomorphism $\sigma^l_I$ as above, induces
a symplectic (up to a similitude) isomorphism of chains of $\O_S$-modules
$$
\tau^l_I: \Upsilon^l_{I, S}=\ker(\pi-a_l\ |\ \Psi^l_{I, S}) \buildrel \sim\over\to\
\ker(\pi-a_l\ |\ \Lambda^l_{I, S}) \simeq   \Xi^l_{I, S}\ .
$$
Similarly, for $l=1$, $\Psi^1_{\pm i_t, S}=\Lambda_{\pm i_t, S}$, and we obtain
a canonical symplectic isomorphism
$$
v_I: \Upsilon^1_{I, S}  \buildrel \sim\over\to\
 \Xi^1_{I, S}\ .
$$

Now if $\{\F^j_{\pm i_t}\}_{j,t}$ is an $S$-valued point of $\N_I$, then
since $$ \F^{l-1}_{\pm i_t}\subset \F^l_{\pm i_t}\subset (\F^{l-1}_{\mp
i_t})^\perp\ ,\quad (\pi-a_l)\F^l_{\pm i_t}\subset \F^{l-1}_{\pm i_t}\ ,
$$ we can consider $\F^l_{\pm i_t}/\F^{l-1}_{\pm i_t}$ as an
$\O_S$-submodule of $\Upsilon^l_{I,S}=\ker(\pi-a_l\ |\ (\F^{l-1}_{\mp
i_t})^\perp/\F^{l-1}_{\pm i_t}\ )$. In fact, $\F^l_{\pm i_t}/\F^{l-1}_{\pm
i_t}$ is locally a direct $\O_S$-summand and $$ \left(\F^l_{\pm
i_t}/\F^{l-1}_{\pm i_t}\right)^\perp=\F^l_{\mp i_t}/\F^{l-1}_{\mp i_t} $$
under the ``derived" forms $<\ ,\ >_{l,\pm i_t}:\Upsilon^l_{\pm i_t,
S}\times\Upsilon^l_{\mp i_t, S}\to \O_S$. Therefore,
\begin{eq}
v_{\pm i_t}(\F^1_{\pm i_t})\ \subset\ \Xi^1_{\pm i_t}\otimes_{\O_K}\O_S\
,\ \ \hbox{resp.}\ \sigma^l_{\pm i_t}(\F^l_{\pm i_t}/\F^{l-1}_{\pm i_t})\
\subset\ \Xi^l_{\pm i_t}\otimes_{\O_K}\O_S
\end{eq}
provide us with $S$-valued points of the ``unramified" local models
$N^1_I$, resp.\ $N^l_I$ for $l=2,\ldots, e$. As in the case of
$\Res_{F/F_0}\GL_d$ we obtain a morphism of schemes $$ p_I: \widetilde
{\N_I}\ \to\ \prod_{j=1}^eN^l_I\  . $$ This is again a torsor for the
smooth group scheme $\prod_{l=2}^e{\cal H}^{(l)}_I$. Again, as in the case
of $\Res_{F/F_0}\GL_d$, we have obtained a diagram of morphisms of schemes
over $\Spec\O_K$:

\begin{eq} \label{diagramSympl}
\matrix{&&&\widetilde\N_I\ \ &&&&&&&\cr
&&&&&&&\cr
&&\ \ p_I\swarrow&&\searrow q_I\ \ &&&&\cr
&&&&&&&&&\cr
&&\prod_{l=1}^eN^l_I\ \ \ \ &&\ \ \ \ \N_I&\buildrel\pi_I\over\to&N^{\rm naive}_I\otimes_{\O_{F_0}}\O_K\cr}
\end{eq}

\no in which both of the slanted arrows are torsors
for the smooth group scheme $\prod_{l=2}^e{\cal H}^{(l)}_I$.
Once again, since by [G2]
the schemes $N^l_I$ are flat over $\Spec\O_K$, the existence
of such a diagram implies:

\begin{thm}\label{flatsympl}
The scheme $\N_I$ is flat over $\Spec\O_K$.\endproof
\end{thm}

\section{Affine flag varieties for the symplectic group} \label{affineSympl}
\setcounter{equation}{0}

In this section, we will use the notations and terminology of \S \ref{aFlag}.
Let us consider
$$
\ti\Lambda_0=k[[\Pi]]^{2g}=k[[\Pi]]\ti e_1\oplus\cdots\oplus k[[\Pi]]\ti e_g\oplus
k[[\Pi]]\ti f_1\oplus\cdots\oplus k[[\Pi]]\ti f_g
$$
with the $k[[\Pi]]$-bilinear alternating form $<\ ,\ >: \ti\Lambda_0\times \ti\Lambda_0\to k[[\Pi]]$
given by
$$
<\ti e_i, \ti e_j>=0\ ,\quad <\ti f_i, \ti f_j>=0\ , \quad <\ti e_i, \ti f_j>=\delta_{ij}\ .
$$ For $0\leq r\leq g$, we  introduce the $k[[\Pi]]$-lattices
$\ti\Lambda_r$ in $\ti\Lambda_0\otimes_{k[[\Pi]]}k((\Pi))$ by $$
\ti\Lambda_r={\rm Span}_{k[[\Pi]]}\{\Pi^{-1}\ti e_1,\cdots ,\Pi^{-1}\ti
e_r,\ti e_{r+1},\cdots, \ti e_g,\ti f_1,\cdots ,\ti f_g\}\ . $$ Set
$\ti\Lambda_{-r}=\widehat{\ti\Lambda_r}:=\{ v\in k((\Pi))^{2g}\ |\ <v,
w>\in k[[\Pi]], \ \hbox{\rm for all $w\in \ti\Lambda_r$}\}$. It is easy to
see that  $$ \ti\Lambda_{-r}={\rm Span}_{k[[\Pi]]}\{\ti e_1, \cdots, \ti
e_g,\Pi \ti f_1,\cdots ,\Pi \ti f_r, \ti f_{r+1},\cdots,\ti f_g\}\ . $$
Now consider a subset $I=\{i_0<\cdots <i_{m-1}\}\subset \{0,\ldots, g\}$.
>From this we obtain the lattice chain $\ti\Lambda_I$ in $k((\Pi))^{2g}$
\begin{eq} \label{latticeSympl}
\ti\Lambda_{-i_{m-1}}\subset\cdots \subset \ti\Lambda_{-i_0}\subset
\ti\Lambda_{i_0}\subset\cdots\subset \ti\Lambda_{i_{m-1}}\subset \Pi^{-1}\Lambda_{-i_{m-1}}\ .
\end{eq}
By adding all the multiples $\Pi^m\Lambda_{\pm i_t}$, $m\in\Z$, to the above lattice chain,
we obtain the corresponding  periodic lattice chain. In what follows, we will
sometimes
use the same symbol to denote both a lattice chain and its corresponding
periodic lattice chain. This should not cause any confusion.
By definition, a lattice chain (which is not necessarily periodic) is
 self-dual if the dual of every lattice in the chain
appears in the corresponding periodic lattice chain. It is
clear that $\ti\Lambda_I$ is a self-dual lattice chain.
\medskip

The
partial affine flag variety $\sflag_I$ associated to the symplectic
similitude group ${\rm GSp}_{2g}$
and the subset $I$ is the ind-scheme over $k$ which represents the functor which
to a $k$-algebra $R$ associates the set of  self dual
$R[[\Pi]]$-lattice chains
\begin{eq} \label{affinelattice}
\L_{-i_{m-1}}\subset\cdots \subset \L_{-i_0}\subset
\L_{i_0}\subset\cdots\subset \L_{i_{m-1}}\subset \Pi^{-1}\L_{-i_{m-1}}
\end{eq}
in $R((\Pi))^{2g}=k((\Pi))^{2g}\otimes_{k}R$, such that each successive
quotient of the above chain is a locally free $R$-module of rank equal to
the $k$-dimension of the corresponding quotient in (\ref{latticeSympl}).
The Ind-group scheme ${\rm GSp}_{2g}(k((\Pi)))$ acts on $\sflag_I$ and we
can identify (${\rm GSp}_{2g}(k((\Pi)))$-equivariantly) $\sflag_I$ with
the fpqc quotient $$ \sflag_I={\rm GSp}_{2g}(k((\Pi)))/P_I $$ where $P_I$
is the parahoric subgroup scheme of ${\rm GSp}_{2g}(k((\Pi)))$ whose
$k$-valued points stabilize the lattice chain $\ti\Lambda_I$ of
(\ref{latticeSympl}).
\medskip

Fix an integer $r$. We may also consider the
partial affine flag variety $\sflag^r_I$ associated to the symplectic
 group ${\rm Sp}_{2g}$
and the subset $I$. This is the ind-scheme over $k$ which represents the functor which
to a $k$-algebra $R$ associates the set of  self dual
$R[[\Pi]]$-lattice chains

\begin{eq} \label{affinelattice}
\L_{-i_{m-1}}\subset\cdots \subset \L_{-i_0}\subset
\L_{i_0}\subset\cdots\subset \L_{i_{m-1}}\subset \Pi^{-1}\L_{-i_{m-1}}
\end{eq}
in $R((\Pi))^{2g}=k((\Pi))^{2g}\otimes_{k}R$, such that
\smallskip

i) each successive quotient of the above chain is a locally free
$R$-module of rank equal to the $k$-dimension of the corresponding
quotient in (\ref{latticeSympl}),
\smallskip

ii) we have $\widehat{\L}_{i_0}=\Pi^r\L_{-i_0}$.
\medskip

 The Ind-group scheme $\Sp_{2g}(k((\Pi)))$ acts
on $\sflag_I^r$. Sending the lattice chain
$\L_{I}$ to $\Pi^m\L_{I}$ gives an $\Sp_{2g}(k((\Pi)))$-equivariant
isomorphism
$$
\sflag_I^r\ \buildrel \sim\over\to\
\sflag^{r-2m}_I  \ .
$$
The Ind-schemes $\sflag^r_I$ are all closed Ind-subschemes
of $\sflag_I$.
In fact, $\sflag^r_I$ for different $r$ are all isomorphic
as Ind-schemes (but not necessarily $\Sp_{2g}(k((\Pi)))$-equivariantly).

\section{Local models and symplectic affine flag varieties} \label{relationsympl}
\setcounter{equation}{0}

Let us identify $\O_F\otimes_{\O_{F_0}}k$ and $k[[\Pi]]/(\Pi^e)$ via
the isomorphism $\O_F\otimes_{\O_{F_0}}k\simeq k[[\Pi]]/(\Pi^e)$
given by $\pi\otimes 1\mapsto \Pi$. Consider
the $k[[\Pi]]/(\Pi^e)$-isomorphism $\Lambda_0\otimes_{\O_{F_0}}k\simeq
\ti\Lambda_0\otimes_{k[[\Pi]]}k[[\Pi]]/(\Pi^e)$
given by
 $e_i\mapsto \ti e_i$, $\delta f_j\mapsto \ti f_j$. This isomorphism
is compatible with the symplectic forms on both sides. In fact, there
are obvious similar isomorphisms
\begin{eq} \label{iso}
\Lambda_i\otimes_{\O_{F_0}}k\simeq \ti\Lambda_i\otimes_{k[[\Pi]]}k[[\Pi]]/(\Pi^e)
\end{eq}
which induce a (symplectic) isomorphism between the
polarized $k[[\Pi]]/(\Pi^e)$-chains $\Lambda_I\otimes_{\O_{F_0}}k$
and $\ti\Lambda_I\otimes_{k[[\Pi]]}k[[\Pi]]/(\Pi^e)$.
\smallskip

Suppose that $\{\F_{\pm i_t}\}_t$ corresponds to an $\Spec R$-valued point
of the special fiber $N^{\rm naive}_I\otimes_{\O_{F_0}}k$ of the naive
local model. Set $\ti\Lambda_{\pm i_t, R}=\ti\Lambda_{\pm
i_t}\otimes_{k[[\Pi]]}R[[\Pi]]$. Let
\begin{eq}
\L_{\pm i_t}\subset \ti\Lambda_{\pm i_t, R}
\end{eq}
be the inverse image of $\F_{\pm i_t}\subset \Lambda_{\pm i_t}\otimes_{\O_{F_0}}R
\simeq \ti\Lambda_{\pm i_t}\otimes_{k[[\Pi]]}{R[[\Pi]]/(\Pi^e)}$
under
$$
\ti\Lambda_{\pm i_t, R}\to \ti\Lambda_{\pm i_t}\otimes_{k[[\Pi]]}R[[\Pi]]/(\Pi^e)\ .
$$

We obtain an $R[[\Pi]]$-lattice chain $\L_I$
$$
 \L_{-i_{m-1}}\subset\cdots \subset \L_{-i_0}\subset
\L_{i_0}\subset\cdots\subset \L_{i_{m-1}}\subset \Pi^{-1}\L_{-i_{m-1}}\ $$
which satisfies property (i) of the definition of the (partial) symplectic
affine flag variety. We claim that $\widehat{\L}_{i_t}=\Pi^{-e}\L_{-i_t}$.
This will establish that the chain above is self-dual and satisfies
property (ii) with $r=-e$. Now we have $$ \L_{-i_t}\subset
\ti\Lambda_{-i_t, R}=\widehat{\ti\Lambda}_{i_t, R}\subset
\widehat{\L}_{i_t}\ . $$ Here the quotients $\ti\Lambda_{-i_t,
R}/\L_{-i_t}$ and therefore $\widehat{\L}_{i_t}/\ti\Lambda_{-i_t, R}$ are
$R$-locally free of rank $eg$. Hence, $\widehat{\L}_{i_t}/\L_{-i_t}$ is
$R$-locally free of rank $2eg$; this is the same as the $R$-rank of
$\Pi^{-e}\L_{-i_t}/\L_{-i_t}$. By our definitions, $<\L_{-i_t}, \L_{i_t}>\
\subset \Pi^e R[[\Pi]]$ and so $\Pi^{-e}\L_{-i_t}\subset
\widehat{\L}_{i_t}$. Since the formation of $\widehat{\L}_{-i_t}$ and
$\Pi^{-e}\L_{-i_t}$ from the $R[[\Pi]]$-lattice $\L_{-i_t}$ commutes with
base change we obtain that $\Pi^{-e}\L_{-i_t}=\widehat{\L}_{i_t}$.
Therefore, the $R[[\Pi]]$-lattice chain $\L_I$ gives an $R$-valued point
of $\sflag_I^{-e}\subset \sflag_I$.
\smallskip

We have therefore obtained a morphism
\begin{eq}
i: N^{\rm naive}_I\otimes_{\O_{F_0}}k\to \sflag_I^{-e}\subset \sflag_I
\end{eq}
which is a closed immersion of Ind-schemes.

Similarly, the special fiber $N^l_I\otimes_{\O_K}k'$ of the ``unramified"
local model $N^l_I$ can be considered as a closed subscheme of the
symplectic affine flag variety $\sflag_I\otimes_kk'$. In fact, by [G2],
$N^l_I\otimes_{\O_K}k'$ is reduced and can be identified with the
scheme-theoretic union of a finite number of Schubert varieties in
$\sflag_I\otimes_kk'$.

Recall that ${\cal H}_I$ is the group scheme over $\Spec\O_{F_0}$ whose
$S$-valued points give the symplectic automorphisms up to similitude of
the polarized chain $\Lambda_I\otimes_{\O_{F_0}}\O_S$. The above
symplectic isomorphism between the polarized
$\O_{F}\otimes_{\O_{F_0}}k=k[[\Pi]]/(\Pi^e)$-chains
$\Lambda_I\otimes_{\O_{F_0}}k$ and
$\ti\Lambda_I\otimes_{k[[\Pi]]}k[[\Pi]]/(\Pi^e)$ allows us to identify the
special fiber $\overline{\cal H}_I$ with the group scheme giving the
symplectic similitude automorphisms of
$\ti\Lambda_I\otimes_{k[[\Pi]]}k[[\Pi]]/(\Pi^e)$. This is a factor group
of the parahoric group scheme $P_I$ giving the symplectic similitude
isomorphisms of $\ti\Lambda_I$. The closed immersion $i$ is equivariant
for the action of $P_I$ in the sense that the action of ${\cal H}_I$
stabilizes the image of $i$, that the action on this image factors through
$P_I\to\overline{\cal H}_I$ and  that $i$ is $\overline{\cal
H}_I$-equivariant.
\medskip

Suppose now that $\{\F_{\pm i_t}^j\}_{j,t}$ corresponds to an $\Spec
R$-valued point of the special fiber $\N_I\otimes_{\O_{K}}k'$ of the
splitting model. For $j=1,\ldots, e$ let
\begin{eq}
\L^j_{\pm i_t}\subset \ti\Lambda_{\pm i_t, R}
\end{eq}
be the inverse image of $\F^j_{\pm i_t}\subset \Lambda_{\pm
i_t}\otimes_{\O_{F_0}}R \simeq \ti\Lambda_{\pm
i_t}\otimes_{k[[\Pi]]}{R[[\Pi]]/(\Pi^e)}$ under $$ \ti\Lambda_{\pm i_t,
R}\to \ti\Lambda_{\pm i_t}\otimes_{k[[\Pi]]}R[[\Pi]]/(\Pi^e)\ . $$ As
above, we obtain an $R[[\Pi]]$-lattice chain $\L^j_I$ $$
 \L^j_{-i_{m-1}}\subset\cdots \subset \L^j_{-i_0}\subset
\L^j_{i_0}\subset\cdots\subset \L^j_{i_{m-1}}\subset
\Pi^{-1}\L^j_{-i_{m-1}} \ .$$ Using a similar argument as above, one can
see that it satisfies properties (i) and (ii) of the definition with
$r=-2e+j$ and therefore gives an $R$-valued point of the symplectic affine
flag variety $\sflag^{-2e+j}_I$. We obtain morphisms:
\begin{eq}
F^j: \N_I\otimes_{\O_K}k'\ \to\ \sflag^{-2e+j}_I\otimes_k k'\subset \sflag_I\otimes_k k'\ .
\end{eq}
and
\begin{eq}
F=(F^j)_j\ :\ \N_I\otimes_{\O_K}k'\ \to\
\prod_{j=1}^e\sflag^{-2e+j}_I\otimes_k k'\subset \prod_{j=1}^e\sflag_I\otimes_k k'\ .
\end{eq}
The morphism $F$ is a closed immersion.  Exactly as in the case of ${\rm
Res}_{F/F_0}{\rm GL}_d$ we can see that the special fiber $\overline
{\N_I}:=\N_I\otimes_{\O_K}k'$ can be naturally identified with the
geometric convolution of the reduced subschemes ${N}^j_I\otimes_{\O_K}k'$
of the symplectic affine flag variety $\sflag_I\otimes_kk'$. Similarly,
the special fiber of the diagram (\ref{diagramSympl}) relates to the
convolution diagram for these subschemes (the analogue of
(\ref{convdiagram})) in exactly the same fashion as it was explained in \S
\ref{relaflag} for $G={\rm Res}_{F/F_0}{\rm GL}_d$.

\section{The canonical local model for $G=\Res_{F/F_0}{\rm GSp}_{2g}$} \label{canonicalsymplectic}
\setcounter{equation}{0}

\begin{Definition} \label{cansympl}
{\rm The canonical local model $N^{\rm can}_I:=N^{\rm can}(\O_F,
\Lambda_I, <\ ,\ >)$ for the group $G=\Res_{F/F_0}{\rm GSp}_{2g}$ and the
self-dual lattice chain $\Lambda_I$ is the scheme theoretic image of the
morphism $$ \pi'_I\ :\ \N_I\to N^{\rm naive}_I\otimes_{\O_{F_0}}\O_K\to
N^{\rm naive}_I $$ which is obtained by composing the morphism $\pi_I$
with the base change morphism.}
\end{Definition}

Since $\pi_I'$ is proper, the canonical local model $N^{\rm can}_I$ is a closed
subscheme of the naive local model $N^{\rm naive}_I$. Since $\pi_I\otimes_{\O_{K}}K:
\N_I\otimes_{\O_{F_0}}K\to N^{\rm naive}_I\otimes_{\O_{F_0}}K$ is
an isomorphism,
 $N^{\rm can}_I$ and $N^{\rm naive}_I$ have the same generic fiber.
The scheme $N^{\rm can}_I$ is flat over $\Spec\O_{F_0}$ since
$\N_I$ is flat over $\Spec\O_K$.
Therefore, $N^{\rm can}_I$ is the (flat) scheme
theoretic closure of the generic fiber $N^{\rm naive}_I\otimes_{\O_{F_0}}F_0$
in $N^{\rm naive}_I$.

Suppose now that $I=\{0\}$ or that $I=\{g\}$. In this case, the self-dual
lattice chain $\Lambda_I$ consists of $\{\pi^t\Lambda_0\}_{t\in \Z}$,
resp. $\{\pi^t\Lambda_g\}_{t\in\Z}$ (we have $\hat\Lambda_0=\Lambda_0$,
$\hat\Lambda_g=\pi\Lambda_g$) and the subgroup of $G(F_0)={\rm
GSp}_{2g}(F)$ which stabilizes $\Lambda_I$ is a special maximal parahoric
subgroup. Then it follows that the unramified local models $N^l_I$ are
smooth Lagrangian Grassmannians over $\Spec\O_K$. Hence, we deduce that
$\prod_{l=1}^eN^l_I$ is irreducible and smooth over $\Spec\O_K$. Since
$\prod_{l=2}^e{\cal H}^{(l)}_I$ is a smooth group scheme with
geometrically connected fibers, we conclude, using the diagram
(\ref{diagramSympl}), that, in this case, $\N_I$ is also irreducible and
smooth over $\Spec\O_K$; therefore the special fiber
$\N_I\otimes_{\O_K}k'$ is irreducible. As a result, the special fiber
$N^{\rm can}_I\otimes_{\O_{F_0}}k$ of the canonical local model $N^{\rm
can}_I$ is irreducible. More generally, suppose that $I=\{i_0\}$ consists
of one index only. This is the case in which the subgroup of $G(F_0)={\rm
GSp}_{2g}(F)$ which stabilizes $\Lambda_I$ is a maximal parahoric
subgroup. Then by [G2], the geometric special fibers of the unramified
local models $N^l_I$ are irreducible. As above, we conclude that the
special fiber $\overline{N}^{\rm can}_I=N^{\rm can}_I\otimes_{\O_{F_0}}k$
is once again irreducible. In fact, we can then show more:
\medskip

\begin{thm} \label{symplemax}
Suppose that $I=\{i_0\}$ consists of one index only.
 Then:
\smallskip

(i) $N^{\rm can}_I$  is normal and Cohen-Macaulay.
\smallskip

(ii) The special fiber $\overline{N}^{\rm can}_I$ is integral and normal
with rational singularities. It can be identified with the Schubert
variety $\overline{\cal O}_{e\mu_1}$ in $\sflag_I$, where $\mu_1$ is the
coweight $(1^g, 0^g)$ of ${\rm GSp}_{2g}$.

\end{thm}

\begin{Proof}
This follows closely the arguments in [PR], proofs of Propositions 5.2--5.3
(see also loc. cit. Remark 5.5). For simplicity of notation, we will
drop the subscript $I$ and write $N$ instead of $N^{\rm naive}$. We will
also use a bar to denote the special fiber of a scheme
over $\O_K$ or over $\O_{F_0}$, depending on the context.
Consider the
proper morphism
$$
\pi\ :\ \N\to N\otimes_{\O_{F_0}}\O_K\ .
$$
Let $N'=\Spec(\pi_*(\O_\N))$ and consider the scheme-theoretic image
$\pi(\N)\subset N\otimes_{\O_{F_0}}\O_K$. Since $\N$ is flat over $\Spec\O_K$
the same is true for $\pi(\N)$.  Let $\varpi$ be a uniformizer of $\O_K$. The cohomology exact sequence
obtained by applying $\pi_*$ to
$$
0\to \O_\N\buildrel \varpi\over\to \O_\N\to \O_{\overline\N}\to 0
$$
gives an injective homomorphism
$$
\O_{N'}/\varpi\O_{N'}\to \overline \pi_*(\O_{\overline \N})\ .
$$
This fits in a commutative diagram:
\begin{eq}\label{com}
\matrix{\O_{\pi(\N)}/\varpi\O_{\pi(\N)}&\to &\O_{\overline \pi(\overline \N)}\cr
\downarrow&&\downarrow\cr
\O_{N'}/\varpi\O_{N'}&\to&\overline \pi_*(\O_{\overline \N})\cr}\ .
\end{eq}
By the definition of the scheme theoretic image the upper horizontal
homomorphism is surjective. Since by the discussion before the statement
of the theorem, $\overline \N$ is reduced and irreducible, the same is
true for the scheme-theoretic image $\overline \pi(\overline \N)\subset
\overline N$. Let $\mu_1$ be the miniscule coweight $(1^g, 0^g)$ of ${\rm
GSp}_{2g}$. The special fibers $\overline N^l_I$ of the corresponding
unramified models can be identified with the Schubert variety $\overline
{\O}_{\mu_1}$ in the affine partial flag variety $\sflag_I\otimes_kk'$
(see  [G2]). By \S \ref{relationsympl}, \S \ref{relaflag} the morphism
$\overline \pi: \overline \N\to \overline \pi(\overline \N)\subset
\overline N$ can be identified with the convolution morphism $$
m_{(\mu_1,\ldots, \mu_1)}: \overline \O_{\mu_1}\ti\times\cdots\ti\times
\overline\O_{\mu_1}\ \to\ \overline \O_{e\mu_1}\subset \sflag_I\otimes_kk'
$$ This morphism is birational on its image. The scheme $\overline
\pi(\overline \N)$ can be identified with the Schubert variety
$\overline\O_{e\mu_1}$ in $\sflag_I\otimes_kk'$; it is therefore normal
with rational singularities ([Fa], [G2]). Since $\overline \pi$ is proper,
the natural morphism
$
\Spec(\overline \pi_*(\overline \N))\to \overline \pi(\overline \N)
$
is finite, and now since by the above $\overline \pi(\overline \N)$ is normal,
$
\Spec(\overline \pi_*(\overline \N))\to \overline \pi(\overline \N)
$
is actually an isomorphism. We conclude that in the diagram (\ref{com})
above, the right vertical homomorphism is an isomorphism. An argument as
in [PR] proof of Proposition 5.2 now implies that the homomorphisms
$\O_{\pi(\N)}/\varpi\O_{\pi(\N)}\to \O_{\overline \pi(\overline \N)}$ and
$\O_{N'}/\varpi\O_{N'}\to \overline \pi_*(\O_{\overline \N})$ which appear
in (\ref{com}) are also isomorphisms. Therefore, the special fibers of
$N'$ and $\pi(\N)$ coincide and they are both equal to $\overline
\pi(\overline \N)$ which by the above is integral, normal and with
rational singularities. In fact, we can see as in loc. cit. that
$N'=\pi(\N)$ and that $\pi(\N)$ is normal and Cohen-Macaulay. To deduce
the claims of the theorem for $N^{\rm can}$ we can now proceed along the
lines of [PR], proof of Proposition 5.3: Recall that the canonical local
model is the scheme-theoretic image of the morphism $$ \pi': \N\to
N\otimes_{\O_{F_0}}\O_K\to N\ , $$ i.e $N^{\rm can}=\pi'(\N)$. An argument
as in  loc. cit. now shows that $$ N^{\rm
can}\otimes_{\O_{F_0}}\O_K=\pi(\N)\ ,\ \  \pi(\N)/{\rm Gal}(K/F_0)=N^{\rm
can}\ , $$ and the desired statements for $N^{\rm can}$  follow (see loc.
cit. for more details).
\endproof

\medskip

\begin{Remark} {\rm It follows that $\overline N^{\rm can}_I=\overline{\cal O}_{e\mu_1}$
 is the union of all the Schubert strata (cells)
in $\sflag_I$ which correspond to double cosets in the extended affine
Weyl group which, in the Bruhat order, are $\leq $ to the coset given by
the coweight $\mu=e\mu_1$. The set of these cosets is exactly the
$\mu$-admissible set as defined in ([KR]).}
\end{Remark}

We now consider general index sets $I$. For $I=\{i_0,\ldots,
i_{m-1}\}\subset \{0,\ldots, g\}$, and $i_t\in I$, we can consider the
morphism $$ \pi_{i_t} : N^{\rm naive}_I\ \to\ N^{\rm naive}_{\{i_t\}}\ $$
obtained by $\{\F_{\pm i_n}\}_{n=0}^{m-1}\mapsto \F_{\pm i_t}$. As in the
case of $G={\rm Res}_{F/F_0}\GL_d$ (see \S \ref{canonical} and [PR], \S
8), we can consider the scheme theoretic intersection in $N_I^{\rm
naive}$,
\begin{eq}\label{int}
N_I^{\rm loc}:= \bigcap_{i_t\in I}\pi^{-1}_{i_t}(N^{\rm can}_{\{i_t\}})\ .
\end{eq}

\begin{thm} \label{symple}
(a)  $N_I^{\rm can}= N_I^{\rm loc}$.

\smallskip
\noindent
 (b) The special fiber $N_I^{\rm can}\otimes_{{\cal O}_{F_0}}k$
is reduced and its irreducible components are normal with rational
singularities. It can be identified with the union in $\sflag_I$ of the
Schubert cells ${\cal O}_w$ with $w$ in the $e\mu_1$-admissible set in
$\tilde W_I\backslash \tilde W /\tilde W_I$.
\end{thm}

Here $\tilde W$ denotes the extended affine Weyl group of ${\rm
GSp}_{2g}(k[[\Pi]])$ and $\tilde W_I$ the subgroup corresponding to the
parahoric subgroup $P_I$.

\medskip

\begin{Proof}
We consider the chain of closed embeddings of ${\cal O}_{F_0}$-schemes
with identical generic fibers, $$N_I^{\rm can}\subset N_I^{\rm loc}\subset
N_I^{\rm naive}\ \ .$$ By [G3], Prop.\ 6.1 all generic points of the
special fiber of $N_I^{\rm naive}$ can be lifted to the generic fiber. In
other words, the above inclusions induce bijections on the underlying
topological spaces. In fact, these bijections follow directly from [G3] Theorem 7.2
which states that the 
$\mu$-admissible and $\mu$-permissible sets coincide, cf. Remark \ref{newproof}.
On the other hand, by Theorem \ref{symplemax} the
special fiber of $N_{\{ i_k\}}^{\rm can}$ is reduced and hence may be
identified with a Schubert variety in a symplectic Grassmannian. Now the
same argument as in the proof of Theorem \ref{thmintersection} implies
that the special fiber of (\ref{int}) is reduced with all its irreducible
components normal and with rational singularities. It follows that
$N_I^{\rm loc}$ is flat over ${\rm Spec}\ {\cal O}_{F_0}$ and hence
$N_I^{\rm can}= N_I^{\rm loc}$. The last statement of (b) follows as
Remark  \ref{conj}, (b) from Section \ref{relationsympl}.
\endproof
\end{Proof}

\begin{Remark}
{\rm  It seems plausible to expect that $N_I^{\rm can}=N_I^{\rm naive}$,
i.e.\ that $N_I^{\rm naive}$ is flat over ${\rm Spec}\ {\cal O}_{F_0}$,
comp.\ [G3].
\par

Let $I=\{ 0\}$. The conjecture above may be reduced to a question on a
certain space of matrices. Let

\begin{eqnarray*}
P=\{ & A=\left(\matrix{a&b\cr 0&~^ta\cr}\right)\in M_{2ge}\ ;\ a,b\in
M_{ge}\ ,\ ~^tb=-b\ ,\\ & {\rm char}_a(T)=
(\prod\limits_{i=1}^e(T-a_i))^g\ ,Q(A)=0\ \ \}\ .
\end{eqnarray*}

The question is whether $P$ is flat over $\Spec\ {\cal O}_{F_0}$.

The relation to the previous conjecture is given by the following diagram
analogous to [PR], (1.3), $$ N^{\rm naive}_{\{ 0\}}
\buildrel\pi\over\longleftarrow \tilde N^{\rm naive}_{\{ 0\}}
\buildrel\phi\over\longrightarrow P\ \ .$$ Here $$\tilde N^{\rm naive}_{\{
0\}}(S)= \{ ({\cal F}\subset\Lambda_{0,S},\alpha)\}\ \ ,$$ where ${\cal
F}$ defines a point of $N^{\rm naive}_{\{ 0\}}(S)$ and where $\alpha$ is a
symplectic automorphism of $\Lambda_{0,S}$ which carries ${\cal F}$ into
the Lagrangian subspace ${\cal F}_0$ of $\Lambda_0$ generated over ${\cal
O}_F$ by $e_1,\ldots, e_g$. Then $\pi$ is a torsor under the Siegel
parabolic in ${\rm Sp}_{2ge}\simeq {\rm Sp}(\Lambda_0,<\ ,\ > )$ and
$\phi$ is a smooth morphism, given by $$\phi(({\cal F},\alpha))
=\alpha^{-1}\cdot\pi\cdot\alpha\ \ ,$$ which we express as a matrix in
terms of the ${\cal O}_S$-basis
 $e_1,\ldots, e_g, \pi
e_1,\ldots, \pi e_g,\ldots,\break \pi^{e-1}e_1,\ldots, \pi^{e-1}e_g,
\delta f_1,\ldots, \delta f_g, \pi \delta f_1,\ldots, \pi\delta
f_g,\ldots, \pi^{e-1}\delta f_1,\ldots, \pi^{e-1}\delta f_g$ of
$\Lambda_{0,S}$. }
\end{Remark}

\part{}

\section{Nearby cycles} \label{nearby}
\setcounter{equation}{0}

In this section, we will assume that the residue field $k$ of $\O_{E}$ is
finite. Our aim is to describe the sheaves of nearby cycles for the local
models $M^{\rm can}_I\otimes_{\O_E}\O_K$ and $N^{\rm
can}_I\otimes_{\O_E}\O_K$ as convolutions of the sheaves of nearby cycles
associated to the ``unramified" local models $M^j_I$ and $N^j_I$
respectively (see below for a precise statement). For simplicity, we will
restrict our discussion mostly to the case of $G={\rm Res}_{F/F_0}\GL_d$,
i.e to the models $M^{\rm can}_I$; the case of $G={\rm Res}_{F/F_0}{\rm
GSp}_{2g}$ is similar.

Fix a prime number $\ell$ which is invertible in $\O_E$ and a square root
of the cardinality $|k|$ in $\overline\Q_{\ell}$. Let $\O$ be a discrete
valuation ring
which is a finite flat extension of $\O_E$ with fraction field $L$
contained in $F_0^{\rm sep}$. If $X$ is a scheme of finite type over
$\Spec\O$ with constant relative dimension $d$ denote by $$
R\Psi^X_L=R\Psi^X_{\O}\overline\Q_{\ell}[d]\hbox{$\left({d\over 2}\right)$} $$ the
(adjusted) complex of nearby cycles of $X$ over $\Spec\O$. This is an
element in the derived category of complexes of $\ov\Q_{\ell}$-sheaves on
the geometric special fiber $X\otimes_{\O}\bar k$  with bounded
constructible cohomology sheaves and continuous ${\rm Gal}(F^{\rm
sep}_0/L)$-action which lifts (i.e is compatible with) the action of ${\rm
Gal}(F^{\rm sep}_0/L)$ on $X\otimes_{\O}\bar k$ through the Galois group
of the residue field of $\O$. If $X$ has smooth generic fiber 
then by [I], Theorem 4.2 and Cor. 4.5,
$R\Psi^X_L$ is a Verdier self dual perverse $\overline\Q_{\ell}$-sheaf  on
$X\otimes_{\O}\bar k$. In fact, under this assumption, G\"ortz-Haines  ([GH] Appendix Theorem 10.1) 
show using de Jong's alteration theorem,
Weil II and the calculations in [RZ2] that 
$R\Psi^X_L$ is also mixed.

For simplicity, if $X$ is a scheme over $\O_E$ with smooth generic fiber, we will write $R\Psi^X_K$
instead of $R\Psi^{X\otimes_{\O_E}\O_K}_K$ for the (adjusted) complex of
nearby cycles of $X\otimes_{\O_E}\O_K$ over $\O_K$. Again, this is a
perverse $\overline\Q_{\ell}$-sheaf on $X\otimes_{\O_E}\bar
k$ with an action of ${\rm Gal}(F^{\rm sep}_0/K)$; it is isomorphic to the
complex of $\overline\Q_{\ell}$-sheaves $R\Psi_E^X$ with the ${\rm
Gal}(F^{\rm sep}_0/E)$-action restricted to the subgroup ${\rm Gal}(F^{\rm
sep}_0/K)$.
\smallskip

By \S \ref{relaflag} and Remark \ref{conj} part (b), the special fiber
$\overline M^{\rm can}_I$ can be naturally identified with a reduced
finite union of Schubert varieties in the partial affine flag variety
$\flag_I$. On the other hand, for each $j=1,\ldots, e$, the special fiber
$\overline M^j_I$ of the unramified local model $M^j_I$ over $\Spec\O_K$
can also be identified with a finite union of Schubert varieties in
$\flag_I\otimes_kk'$. In this way, we can regard $$ \quad R\Psi^{M^j_I}_K,
\ \ R\Psi^{M^{\rm can}_I}_K $$ as perverse $\overline\Q_{\ell}$-sheaves
 on $\flag_I\otimes_k\bar k$ with compatible ${\rm Gal}(F^{\rm
sep}_0/K)$-actions which are $P_I$-equivariant. By Remark \ref{conj} (b),
these perverse sheaves are supported on the union of Schubert cells
corresponding to the $\mu_j$-admissible, resp. $\mu$-admissible cosets,
where $\mu=\mu_1+\cdots +\mu_e$.

For each $j=1,\ldots, e$, we now let $\Phi_j$ be a $P_I$-equivariant
perverse $\overline\Q_{\ell}$-sheaf on $\overline M^j_I$ with compatible
${\rm Gal}(F^{\rm sep}_0/K)$-action. The convolution construction of
Ginzburg, Lusztig, etc. (see for example [Lu]) allows us to construct  an
element $$ \Phi_1\star\cdots\star \Phi_e $$ in the derived category of
complexes of $\ov\Q_{\ell}$-sheaves on $\flag_I\otimes_k\bar k$ supported
on $\ov M^{\rm can}_I$ with bounded constructible cohomology sheaves and
compatible ${\rm Gal}(F^{\rm sep}_0/K)$-action. (In what follows, for
simplicity of notation, we will use a bar to denote the geometric special
fiber over $\bar k$ and omit the base change from the notation). The
construction proceeds as follows ([Lu] 1.2 and 1.3). Consider the diagram
obtained by the convolution diagram (\ref{convdiagram}) by base changing
from $k'$ to $\bar k$:

\begin{eq}
\matrix{&&& U_{\bar k} &&&\cr
&&&&&&\cr
&&p_1\swarrow&&\searrow p_2&&&&\cr
&&&&&&\cr
&&\ov {M}^1_I\times\cdots\times \ov{M}^e_I\ &&
\ov{M}^1_I\ti\times\cdots\ti\times \ov{M}^e_I &\ \buildrel p_3\over\to\
\ov{M^{\rm naive}_I\otimes_{\O_E}{\O_K}} \subset \flag_I\otimes_k\bar k\ .
& \cr}
\end{eq}

Consider the pull back of the exterior tensor product
$p_1^*(\Phi_1\boxtimes\cdots \boxtimes\Phi_e)$; since $p_1$ is a smooth
morphism, this is a perverse $\ov \Q_{\ell}$-sheaf up to a shift by the
relative dimension of $p_1$. By its definition,
$p_1^*(\Phi_1\boxtimes\cdots \boxtimes\Phi_e)$ is equivariant for the
action (\ref{action3}); however, since the complexes of sheaves $\Phi_j$
are $P_I$-equivariant, it is also equivariant for the action
(\ref{action4}). Recall that $p_2$ is a $P_I$-torsor for the action
(\ref{action4}) (which is actually locally trivial in the Zariski
topology). Therefore, by descent (see also [BBD] Theorem 4.2.5), there
is a perverse $\ov\Q_{\ell}$-sheaf with compatible ${\rm Gal}(F^{\rm
sep}_0/K)$-action $$ \Phi_1\ti\boxtimes\cdots\ti\boxtimes\Phi_e $$ on
$\ov{M}^1_I\ti\times\cdots\ti\times \ov{M}^e_I$, which is unique up to
unique isomorphism,  such that $$
p^*_2(\Phi_1\ti\boxtimes\cdots\ti\boxtimes\Phi_e)=p_1^*(\Phi_1\boxtimes\cdots
\boxtimes\Phi_e)\ . $$ We now set $$ \Phi_1\star\cdots\star
\Phi_e:=R{p_3}_*(\Phi_1\ti\boxtimes\cdots\ti\boxtimes\Phi_e)\ . $$
\medskip

\begin{thm} \label{convo}
\ \smallskip

(a) The sheaf $R\Psi^{M^1_I}_K   \star\cdots \star
R\Psi^{M^e_I}_K$ on $\flag_I\otimes_k\bar k$ is mixed perverse and
Verdier self dual.
\smallskip

(b) There is an isomorphism of perverse $\overline\Q_{\ell}$-sheaves with
${\rm Gal}(F^{\rm sep}_0/K)$-action 
$$
 R\Psi^{M^{\rm can}_I}_K \simeq
R\Psi^{M^1_I}_K   \star\cdots \star R\Psi^{M^e_I}_K 
$$ 
on
$\flag_I\otimes_k\bar k$.
\end{thm}

\bigskip

\begin{Proof} Recall the diagram
(\ref{diagram})
\begin{eq}
\matrix{&&&\widetilde\M_I\ \ &&&&&&&\cr
&&&&&&&\cr
&&\ \ p_I\swarrow&&\searrow q_I\ \ &&&&\cr
&&&&&&&&&\cr
&&\prod_{l=1}^eM^l_I\ \ \ \ &&\ \ \ \ \M_I&\buildrel\pi_I\over\to&M^{\rm naive}_I\otimes_{\O_E}\O_K\ .\cr}
\end{eq}
By the K\"unneth formula [I], Theorem 4.7, we have an isomorphism of
perverse $\ov\Q_{\ell}$-sheaves with compatible ${\rm Gal}(F^{\rm
sep}_0/K)$-action on the geometric special fiber of $M^1_I\times\cdots
\times M^e_I$,
\begin{eq}\label{prod}
R\Psi^{M^1_I\times\cdots\times M^e_I}_K\ \simeq\
R\Psi^{M^1_I}_K   \boxtimes\cdots \boxtimes
R\Psi^{M^e_I}_K \ .
\end{eq}
This induces an isomorphism between the pull-backs
\begin{eq}\label{prod2}
\ov p_I^*(R\Psi^{M^1_I\times\cdots\times M^e_I}_K)\ \simeq\
\ov p_I^*(R\Psi^{M^1_I}_K   \boxtimes\cdots \boxtimes
R\Psi^{M^e_I}_K)\ .
\end{eq}
>From the definitions, and using the comparisons of the special fiber of
the diagram (\ref{diagram}) with the convolution diagram
(\ref{convdiagram}) explained at the end of \S \ref{relaflag}, we obtain
an isomorphism
\begin{eq}\label{prod3}
\ov p_I^*(R\Psi^{M^1_I}_K   \boxtimes\cdots \boxtimes
R\Psi^{M^e_I}_K)\ \simeq\ \ov q_I^*(R\Psi^{M^1_I}_K   \ti \boxtimes\cdots \ti \boxtimes
R\Psi^{M^e_I}_K)\ .
\end{eq}
Since both $p_I$ and $q_I$ are smooth, $p^*_I$ and $q_I^*$ commute with
the nearby cycle functor. Therefore, we obtain an isomorphism $$ \ov
p_I^*(R\Psi^{M^1_I}_K   \boxtimes\cdots \boxtimes R\Psi^{M^e_I}_K)\
\simeq\ R\Psi^{{\ti\M}}_K\ \simeq\ \ov q_I^*(R\Psi^{\M}_K) $$ which by
[BBD] Theorem 4.2.5   and (\ref{prod3}) gives an isomorphism of perverse
$\ov\Q_{\ell}$-sheaves with compatible ${\rm Gal}(F^{\rm
sep}_0/K)$-action
\begin{eq}\label{prod4}
  R\Psi^{\M}_K\ \simeq\ R\Psi^{M^1_I}_K   \ti \boxtimes\cdots \ti \boxtimes
R\Psi^{M^e_I}_K\ .
\end{eq}
We now notice that since $\pi_I:\M\to M^{\rm can}_I\otimes_{\O_E}\O_K\subset
 M^{\rm naive}_I\otimes_{\O_E}\O_K$ is
  proper and since $\pi_I$ induces an isomorphism on the generic fibers, there
is a canonical isomorphism
\begin{eq}
R\ov {\pi_I}_* (R\Psi^{\M}_K)\ \simeq\  R\Psi^{M^{\rm can}_I}_K\ .
\end{eq}
Hence, by (\ref{prod4}) there is an isomorphism
\begin{eq}
R\Psi^{M^{\rm can}_I}_K\ \simeq\
R\ov {\pi_I}_*  (R\Psi^{M^1_I}_K   \ti \boxtimes\cdots \ti \boxtimes
R\Psi^{M^e_I}_K)=R\Psi^{M^1_I}_K   \star\cdots \star
R\Psi^{M^e_I}_K
\end{eq}
with the last equality given by the identification of $p_3$
with $\ov \pi_I$.
This establishes both parts (a) and (b) of the Theorem.
\endproof

We note that the factors $R\Psi^{M^l_I}_K$ are known perverse sheaves,
{\bf at least if $0\in I$,}
thanks to the result of Haines and Ng\^o regarding the unramified case
[HN1].

\begin{Remark} {\rm 

(a) When $I=\{0\}$, the nearby cycles $R\Psi^{M^{\rm can}_I}_E$ are pure of weight $0$, since the splitting model is smooth in this case, comp. [PR]. In the general case, it is an interesting problem to determine the weights occurring in  $R\Psi^{M^{\rm can}_I}_E$ and their multiplicities, comp [GH]. 

(b) The same arguments applied to the local models $N^{\rm can}_I$
for the group $G={\rm Res}_{F/F_0}{\rm GSp}_{2g}$ show that
\begin{eq}
R\Psi^{N^{\rm can}_I}_K\simeq R\Psi^{N^1_I}_K   \star\cdots \star
R\Psi^{N^e_I}_K
\end{eq}
as perverse $\ov\Q_{\ell}$-sheaves with ${\rm Gal}(F^{\rm
sep}_0/K)$-action.

 (c) Theorem \ref{convo} determines the 
nearby cycles of $M^{\rm can}_I$ over $K$. To obtain the 
nearby cycles $R\Psi^{M^{\rm can}_I}_E$ over the reflex field $E$ 
one needs to specify in addition the corresponding -via the
isomorphism of Theorem \ref{convo} (b)- ${\rm Gal}(K/E)$-action on the convolution product 
$$
R\Psi^{M^1_I}_K   \star\cdots \star R\Psi^{M^e_I}_K\ .
$$
Let us identify $\sigma\in {\rm Gal}(K/E)$ 
with a permutation of the set $\{1,\ldots, e\}$ via 
the action of $\sigma$ on the set of embeddings $K\to F_0^{\rm sep}$.
We then expect that the action of $\sigma$ on the convolution
product should be given using the ``commutativity isomorphisms"
of [HN1] Proposition 22. (This in turn 
is a version of the isomorphism of [Ga] Theorem 1 (b).)
In the case that $I=\{0\}$ this issue is 
discussed in some more detail in [PR] Remark 7.4.

}
\end{Remark}

\end{Proof}
\end{Proof}
\part{}

\section{Splitting and local models in the general PEL case} \label{general}
\setcounter{equation}{0}

In this section, we explain the construction of
splitting models in the general (ramified) PEL case.
As we shall see this also suggests a general construction of local models.
We take $F_0=\Q_p$ in the notation used elsewhere in this paper.
Specifically, we will use the following notation (following closely [RZ], see 1.38):
\smallskip

$\bullet$ $F$ a finite direct product of finite field extensions of $\Q_p$,

$\bullet$ $B$ a finite central algebra over $F$,

$\bullet$ $V$ a finite dimensional (left) $B$-module,

$\bullet$ $(\ ,\ )$ a nondegenerate alternating $\Q_p$-bilinear form on $V$,

$\bullet$ $b\mapsto b^*$ an involution on $B$ which satisfies
$(bv, w)=(v, b^*w)$, \ $v$, $w\in V$,

$\bullet$ $\O_B$ a maximal order of $B$ invariant under $*$.

\smallskip

If $W$ is a right $B$-module, we define a left $B$-module on $W$
by restriction of scalars $*: B\to B^{\rm opp}$. With this convention
the dual vector space $V^*={\rm Hom}_{\Q_p}(V,\Q_p)$ is a left $B$-module
and the form $(\ ,\ )$ induces an isomorphism of $B$-modules
$$
\psi: V\to V^*\ .
$$
In the same way, for an $\O_B$-lattice $\Lambda$ in $V$, the ${\Z_p}$-module
$\Lambda^*={\rm Hom}_{\Z_p}(\Lambda, \Z_p)$
becomes a left $\O_B$-module. The image of $\Lambda^*$ under the map
$$
\Lambda^*\to V^*\buildrel\psi^{-1}\over \simeq V
$$
is the ``dual" lattice $\hat\Lambda$ of $\Lambda\subset V$
with respect to $(\ ,\ )$. The form $(\ ,\ )$ induces a perfect bilinear pairing
$$
(\ ,\ ): \Lambda\times\hat\Lambda\to \Z_p\ .
$$

\smallskip

Let $F_1$ be the $\Q_p$-algebra which consists of the
$*$-invariant elements of $F$. For simplicity
we will assume that $F_1$ is a field; the local models in
the general case are products of local models for
cases in which $F_1$ is a field. We will denote by $\tau$ the automorphism of $F$
obtained by restricting the involution $*$. There are three cases:
\smallskip

(I) $F=F_1\times F_1$ and $\tau(a_1, a_2)=(a_2, a_1)$,

(II) $F=F_1$,

(III) $F$ is a quadratic field extension of $F_1$.
\smallskip

The existence of the $*$-linear form $(\ ,\ )$ implies that, even in case I, $V$ is a free $F$-module;
we will denote its rank by $d$.

Let $G$ be the algebraic group over $\Q_p$, whose points with values
in a $\Q_p$-algebra $R$ are given by:
$$
G(R)=\{g\in \GL_B(V\otimes_{\Q_p}R)\ |\ (gv, gw)=c(g)(v, w), c(g)\in R\}\ .
$$

Let us fix in addition
\smallskip

$\bullet$ a cocharacter $\mu: {\Gm}_{N}\to G_N$
defined over the finite extension $N$ of $\Q_p$, given
up to conjugation.
\smallskip

We assume that the corresponding
eigenspace  decomposition of $V\otimes_{\Q_p}N$ is given by
$$
V\otimes_{\Q_p}N=V_0\oplus V_1
$$
(i.e the only weights are $0$ and $1$) and that
the composition $c\circ \mu: {\Gm}_N\to{\Gm}_N$ is the identity. This implies that
both $V_0$ and $V_1$ are totally isotropic for the
form on $V\otimes_{\Q_p}N$ obtained by $(\ ,\ )$ by extending scalars
(by [RZ] Definition 3.18 and 3.19 (b) these conditions
correspond to the situation describing moduli
of $p$-divisible groups). Notice that this implies
that the pairing $(\ ,\ )$ induces an isomorphism
\begin{eq}\label{dual}
V_0\simeq V_1^*={\rm Hom}_{N}(V_1, N)
\end{eq}
where $V_1^*$ becomes a left $B$-module as above, by first regarding
it naturally as a right $B$-module and then composing
with the involution $*: B\to B^{\rm opp}$.
As usual let $E$ be the field of definition
of the conjugacy class of $\mu$. We shall also fix
\smallskip

$\bullet$ $\L$ a selfdual periodic multichain of  $\O_B$-lattices in $V$
([RZ] Definition 3.13).
\smallskip

Recall that ``selfdual" means that if $\Lambda$ is in $\L$
then the dual lattice $\hat\Lambda$ is also in $\L$.
As in loc. cit. we can consider $\L$ as a category with
morphisms given by inclusions of lattices.

Now let $\Phi$ be the
set of $\Q_p$-algebra homomorphisms of $F$ in $\bar\Q_p$.
For $a\in F$
let
$$
\det(T\cdot I-a\ |\ V_1)=\prod_{\phi\in\Phi}(T-\phi(a))^{r_\phi}
$$
so that the cocharacter
$$
\mu_{\bar\Q_p}: {\Gm}_{\bar\Q_p}\to G_{\bar\Q_p}\subset {\GL}_B(V\otimes_{\Q_p}\bar\Q_p)\subset
\GL_F(V\otimes_{\Q_p}\bar\Q_p)=\prod_{\phi\in\Phi}\GL(V\otimes_{F,\phi}\bar\Q_p)
$$
 is given, up to conjugation, by
$\{(1^{r_\phi}, 0^{d-r_{\phi}})\}_{\phi\in\Phi}$ with $d$ the $F$-rank of $V$.
We can think
of the automorphism $\tau$ of $F$ as giving a permutation of $\Phi$ by $\phi\mapsto \phi\cdot \tau$.
For every $\phi\in \Phi$ we have
\begin{eq}
r_{\phi}+r_{\phi\cdot\tau}=d\ .
\end{eq}
 Indeed, by (\ref{dual}),
the sum $r_{\phi}+r_{\phi\cdot\tau}$ is the multiplicity of the
eigenvalue $\phi(a)$ for the action of $a\in F$ on $V\otimes_{\Q_p}\bar\Q_p$.
This is equal to $d$ since $V$ is $F$-free of rank $d$.
\smallskip

Set $m=[F_1:\Q_p]$ and let $n$ be the $\Q_p$-dimension of $F$.
We choose an ordering of the $\Q_p$-algebra homomorphisms $\phi_i: F\to \bar\Q_p$,
$1\leq i\leq n$, which in the case that $F\neq F_1$ has the
property that any two embeddings $\phi$, $\phi'$ with the same restriction
to $F_1$ are successive. Denote by $K$ the Galois closure of $F$ in
$\bar\Q_p$. Then $E\subset K$.

Suppose now that $S$ is an $\O_K$-scheme. In what follows
undecorated tensor products are meant to be over $\Z_p$.
If $b$ is a unit of $B$ which normalizes $\O_B$
and $\Lambda\in \L$ then by the definitions $b\Lambda\in \L$.
For such a $b$, conjugation by $b^{-1}$ defines
an isomorphism $\O_B\to \O_B$, $x\mapsto b^{-1}xb$.
If $M$ is an $\O_B\otimes\O_S$-module we denote by
$M^b$ the $\O_B\otimes\O_S$-module obtained by restriction of scalars
with respect to this isomorphism. Left multiplication by $b$
induces a $\O_B\otimes\O_S$-linear homomorphism
$
b: M^b\to M
$.

\smallskip

Let us now define a functor $\M$ on the category of $\O_N$-schemes.

\begin{Definition} \label{gensplit}
{\rm A point of $\M$ with values in an $\O_N$-scheme $S$ is given by
the following data.

1. For each $i=1,\ldots ,n+1$, a functor from the category
of the multichain $\L$ to the category of $\O_B\otimes \O_S$-modules
on $S$
$$
\Lambda\mapsto F^i_\Lambda\ , \ \Lambda\in\L .
$$

2. For $i=1,\ldots ,n+1$, a morphism of functors
$$
j^i_\Lambda: F^i_\Lambda \to \Lambda\otimes\O_S .
$$

\no We are requiring that the following conditions are satisfied:
\smallskip

a) For each $\Lambda\in\L$, $i=1,\ldots, n+1$, the homomorphism $j^i_\Lambda$ is injective
(and so it identifies $F^i_\Lambda$ with a $\O_B\otimes\O_S$-submodule of $\Lambda\otimes\O_S$).
Both $F^i_\Lambda$ and the quotient $(\Lambda\otimes\O_S)/F^i_\Lambda$ are
finite locally free $\O_S$-modules.
\smallskip

b) If $b$ is a unit of $B$ which normalizes $\O_B$ there are ``periodicity"
$\O_B\otimes\O_S$-linear isomorphisms
$$
\theta_{b, \Lambda}: (F^i_\Lambda)^b\buildrel\sim\over\longrightarrow F^i_{b\Lambda}
$$
which make the diagrams
$$
\matrix{(F^i_\Lambda)^b&\buildrel j^i_\Lambda\over\to &(\Lambda\otimes\O_S)^b\cr
\theta_{b,\Lambda}\downarrow\ \ \ &&\downarrow b\cr
F^i_{b\Lambda}&\buildrel j^i_{b\Lambda}\over\to &b\Lambda\otimes\O_S\cr}
$$
commutative.
\smallskip

c) For the action of $\O_B$ on
$F^1_\Lambda$, we have the following identity of polynomial functions
$$
{\rm det}_{\O_S}(a\ \vert \ F^1_\Lambda)={\rm det}_{\O_S}(a\ \vert \ V_1),\quad a\in\O_B\ .
$$

d) We have $F^{n+1}_\Lambda=(0)$.
For $i=1,\ldots, n$, $F^{i+1}_\Lambda\subset F^i_\Lambda$,
the quotient $F^i_\Lambda/F^{i+1}_\Lambda$ is $\O_S$-locally free
of rank $r_i:=r_{\phi_i}$ and is annihilated
by
$$
a\otimes 1-1\otimes\phi_i(a)\in \O_B\otimes\O_S, \quad \hbox{\rm for all $a\in\O_F$}.
$$

e) Note that (a) implies that $F^i_\Lambda$ is a locally direct $\O_S$-summand
of $\Lambda\otimes\O_S$. We will denote by $(F^i_\Lambda)^\perp$ its orthogonal complement
in $\hat \Lambda\otimes\O_S$
under the perfect pairing
$$
(\ ,\ ): (\Lambda\otimes\O_S)\times (\hat\Lambda\otimes\O_S)\to \O_S\ .
$$
For every $\Lambda\in \L$ and $i=1,\ldots, n+1$,
we  require that $F^i_{\hat \Lambda}\subset (F^i_\Lambda)^\perp$.
\smallskip

f) In addition to the above, we require that:
\smallskip

f1) If $F=F_1$, for every $i=1,\ldots, n$ and $\Lambda\in\L$
$$
\prod_{1\leq k\leq i}(a\otimes1-1\otimes\phi_k(a))((F^{i+1}_{\Lambda})^\perp)\subset F^{i+1}_{\hat\Lambda}
$$
for all $a\in\O_F$.
\smallskip

f2) If $F\neq F_1$, for every $h=1,\ldots, m=[F_1:\Q_p]$ and $\Lambda\in\L$
$$
\prod_{1\leq k\leq 2h}(a\otimes1-1\otimes\phi_k(a))((F^{2h+1}_{\Lambda})^\perp)\subset F^{2h+1}_{\hat\Lambda}
$$
for all $a\in\O_F$.}
\end{Definition}

There is a morphism $\pi: \M\to M^{\rm naive}\otimes_{\O_E}\O_K$, where $M^{\rm naive}$ is the functor
of the ``naive" local model of [RZ] (denoted by $M^{\rm loc}$ in loc. cit.) given
by sending the $S$-point of $\M$ given by
$\Lambda\mapsto (F^i_\Lambda\subset \Lambda\otimes\O_S)_{1\leq i\leq n+1}$
to $\Lambda\mapsto t_\Lambda:=(\Lambda\otimes\O_S)/F^1_\Lambda$. Indeed, the functor $\Lambda\mapsto t_\Lambda$
 satisfies the conditions of loc. cit., Definition 3.27.
For example, (c)
and (e) together with the fact that $F^1_{\hat \Lambda}$, $(F^1_\Lambda)^\perp$
are locally direct $\O_S$-summands of $\hat\Lambda\otimes\O_S$ imply that $F^1_{\hat \Lambda}=(F^1_\Lambda)^\perp$ and so $t_\Lambda$ satisfies   condition (iii)
of loc. cit.

It is clear
that $\M$ is representable by a projective scheme over $\Spec \O_K$ and that
the morphism $\pi$ is projective. We can also see that, on the generic fibers,
$\pi$ induces an isomorphism
$$
\pi\otimes_{\O_K}K : \M\otimes_{\O_K}K\buildrel\sim\over\rightarrow M^{\rm naive}\otimes_{\O_E}K.
$$
Let us use the same symbol $\pi$ for the composed morphism $\pi: \M\to M^{\rm naive}\otimes_{\O_E}O_K\to M^{\rm
naive}$.
The scheme theoretic image $\pi(\M)\subset M^{\rm naive}$
is a closed subscheme of $M^{\rm naive}$ which has the same generic fiber as
$M^{\rm naive}$. One can now set
$$
M^{\rm loc}=\pi(\M).
$$
We believe that, if we exclude the case that the group is orthogonal  and
certain unitary cases, then $M^{\rm loc}$ is a good integral model of its generic fiber.
\footnote{Genestier has
pointed out to us that the orthogonal
case is problematic in this respect.} More precisely, assume that we are either in case (I),
or in case (II) with $*$ an orthogonal involution
(then $G$ is a form of a symplectic group), or in case (III)
with $F/F_1$ unramified. Recall here that an involution of the first kind
on a central simple algebra is called orthogonal resp.\ symplectic, if
after a base change that splits the algebra it becomes the adjoint
involution with respect to a symmetric resp.\ alternating form.
Then it seems that the methods of the present paper prove that $M^{\rm loc}$ is flat over ${\rm Spec}\, \O_E$,
with reduced special fiber, and such that all irreducible components of the special
fiber are normal with rational singularities. Furthermore, let $L$ denote the
completion of the maximal unramified extension of $\Q_p$ and let $\tilde
K=\tilde K_{\cal L}$ be the parahoric subgroup of $G(L)$ which fixes the
lattice chain ${\cal L}\otimes \O_L$ in $V\otimes_{\Q_p}L$. Then $\tilde
K$ acts on $M^{\rm loc}(\overline{{\bf F}}_p)$ and the orbits are
in bijective correspondence with the {\it $\mu$-admissible subset}
${\rm Adm}_{\tilde K}(\mu)$ of $\tilde K\setminus G(L)/\tilde K$.
We refer to [R], section 3, for the definition of the $\mu$-admissible
subset in the general case, cf.\ also [KR].
Our work in the previous sections shows that
all these statements hold true in the following situations 
(and we believe that the general case, as limited above,
may be reduced to these cases):
\medskip

a)  Let $F_1$ be a finite field extension
of $\Q_p$ and consider $B=F_1\times F_1$ with the involution
$(a_1, a_2)^*=(a_2, a_1)$. Let $\O_B=\O_{F_1}\oplus \O_{F_1}$
and take $V=B^d=W_1\oplus W_2$, $W_i=
F_1\cdot e^i_1\oplus\cdots\oplus F_1\cdot e^i_d$, for $i=1$ or $2$,
with the alternating form $(\ ,\ )$ defined by
$$
(e^i_k, e^i_l)=0,\quad (e^1_k, e^2_l)=\delta_{kl}\ , i=1,2;\ k,l=1,\ldots,
d\ \ .
$$
This form identifies $W_2$ with the dual of $W_1$. A selfdual
multichain of lattices in $V$ now is given by
a pair $\L=\{\Lambda_k\}_k$, $\hat\L=\{\hat\Lambda_l\}_l$
where $\L$ is a chain of $\O_{F_1}$-lattices in $W_1$ and $\hat\L$
is the dual chain. In this case,
$$
G=\{(g, c\cdot (g^t)^{-1})\ |\ g\in \GL(W_1), c\in \Gm\}\subset \GL(W_1)\times \GL(W_2)\ .
$$
Therefore, $G\simeq {\rm Res}_{F_1/\Q_p}(\GL_d)\times\Gm$. Let us assume
that $F_1$ is totally ramified over $\Q_p$.
The schemes $\M$ (for various choices
of the cocharacter $\mu$ and the multichain $\L$) can be identified with the splitting models
for ${\rm Res}_{F_1/\Q_p}\GL_d$ of \S \ref{splitGl}. To see this we observe
that by using conditions (e) and (f2) and an argument as in
Lemma \ref{perp} we can show that there is a $1$-$1$
correspondence between submodules
$$
F^i_{\Lambda_k\oplus\hat\Lambda_{l}}=
G^{n+1-i}_{\Lambda_k}\oplus G'^{n+1-i}_{\hat\Lambda_l}\subset (\Lambda_k\otimes\O_S)\oplus
(\hat\Lambda_l\otimes\O_S),\ (i=1,\ldots, n+1)
$$
which correspond to $S$-points of $\M$ and
submodules $G^j_{\Lambda_k}\subset (\Lambda_k\otimes\O_S)$ which correspond
to $S$-points in the splitting model of \S \ref{splitGl}.
Theorem \ref{flat} now implies that the schemes $\M$ are flat
over $\Spec\O_K$.

\medskip

 b) Let $B=F=F_1$ a finite field extension of $\Q_p$
and let  $(V, \{\ ,\ \})$
be the standard symplectic vector space over $F$ of dimension $2g$ with
basis $e_1,\ldots, e_g,f_1,\ldots, f_g$, i.e
\begin{eq}
\{e_i,e_j\}=\{f_i,f_j\}=0,\quad \{e_i,f_j\}=\delta_{ij} .
\end{eq}
We set $(v, w)=\Tr_{F/\Q_p}(\{v,w\})$. In this case,
$G={\rm Res}_{F/\Q_p}{\rm GSp}_{2g}$ and, in case $F_1$ is totally
ramified over $\Q_p$,
the scheme $\M$ can be identified with the splitting model
for ${\rm Res}_{F/\Q_p}{\rm GSp}_{2g}$ of \S \ref{splitSp};
here Theorem \ref{flatsympl} implies the truth of the above conjecture.
\medskip

\begin{Remark}\label{unitary} {\rm An example where the methods
of the previous sections do not directly apply
is provided by the case of a group of unitary similitudes
for a ramified quadratic extension of $\Q_p$.
However, even in this case, there are instances in
which we can show that $M^{\rm loc}$
as defined above, is flat over $\Spec\O_E$. We review some results from
[P]. Let $B=F$ a ramified quadratic extension of $\Q_p$, $p$ odd,
with the involution given by the non-trivial Galois automorphism.
Let $V=F^n$ and denote by $e_i$, $1\leq i\leq n$,
the canonical
$\O_F$-generators of the standard lattice $\Lambda_0:=\O_F^n\subset V$.
Let $\pi$ be a uniformizer of $\O_F$ which satisfies $\pi^*=-\pi$.
We define a non-degenerate alternating
$\Q_p$-bilinear form $(\ ,\ ): V\times V\to \Q_p$
which satisfies $(ax, y)=(x, a^*y)$
for $a\in F$ by setting
$$
(e_i, e_j)=0,\quad (e_i, \pi   e_j)=\delta_{ij}\ ,\ i,j=1,\ldots, n\ \ .
$$
The restriction $(\ ,\ ): \O_F^n\times \O_F^n\to \Z_p$
is a perfect $\Z_p$-bilinear form. Therefore, we have $\hat\Lambda_0=\Lambda_0$
and more generally $\widehat{\pi^n\Lambda_0}=\pi^{-n}\Lambda_0$.
Let $\L$ be the selfdual lattice chain $\{\pi^n\Lambda_0\}_{n\in\Z}$.
Using the duality isomorphism ${\rm Hom}_F(V, F)\simeq {\rm Hom}_{\Q_p}(V, \Q_p)$
given by composing with the trace ${\rm Tr}_{F/\Q_p}: F\to \Q_p$
we see that there exists a unique non-degenerate hermitian form
$\phi: V\times V\to F$
such that
$$
(x,y)={\rm Tr}_{F/\Q_p}(\pi^{-1} \phi(x,y)), \quad x,\ y\in V.
$$
Hence, in this case the group $G$ can be identified with the group
of unitary similitudes of the form $\phi$. Now let $r$, $s$ be two non-negative integers
such that $n=r+s$. Fix  a cocharacter
$\mu_F: {\Gm}_F\to G_F$ such that the corresponding subspace $V_1$ of
$V_F=F^n\otimes F$, when considered as an $F$-module via the first factor, is isomorphic to
$F^r\oplus F^s_\tau$ where $F_\tau$ is the module obtained by $F$ by restriction
of scalars via $\tau: F\to F$. The ``naive" local models that correspond to these choices
have been studied in [P]. As was shown there, when $|r- s|>1$,
they are not flat over $\Spec\O_E$.

Given these choices of PEL data, we can see that $K=F$ and that for any $\Spec\O_F$-scheme $S$
the points $\M(S)$ are now
pairs $(F^2, F^1)$ of $\O_F\otimes\O_S$-submodules of $\Lambda_0\otimes\O_S$
which are locally direct summands as $\O_S$-modules and satisfy
\medskip

i) \ $F^1$ is isotropic for
the form $(\ ,\ )$ on $\Lambda_0\otimes\O_S$;
\medskip

ii) \ $F^2\subset F^1$; and $F^1$, $F^2$ have ranks $n$ and $r$ respectively;
\medskip

iii) \
${\rm det}_{\O_S}(T\cdot I-a\otimes 1\ |\ F^1)=(T-a)^r(T-\tau(a))^s\in \O_S[T]$, for every $a\in\O_F$;
\medskip

iv) \ $(a\otimes 1-1\otimes a)(F^2)=(0)$, \ $(a\otimes 1-1\otimes \tau(a))(F^1)\subset F^2$,
for every $a\in \O_F$
\medskip

(the tensor products are in $\O_F\otimes\O_F$ which maps to $\O_F\otimes\O_S$).
\medskip

For simplicity, let us assume that $r\neq s$; then $K=E=F$. The naive local
model $M^{\rm naive}$ classifies isotropic $\O_F\otimes\O_S$-submodules
$F^1$ of $\Lambda_0\otimes\O_S$ which are locally direct summands
of rank $n$ as $\O_S$-modules and satisfy condition (iii) above;
the morphism $\pi: \M\to M^{\rm naive}$ corresponds
to forgetting $F^2$. We can see that the scheme theoretic image
$M^{\rm loc}:=\pi(\M)\subset M^{\rm naive}$ is contained in the closed
subscheme $M'_{r,s}$ of $M^{\rm naive}$ described by
$$
\wedge^{r+1}(a\otimes 1-1\otimes \tau(a)\ |\ F^1)=(0),
\  \wedge^{s+1}(a\otimes 1-1\otimes a\ |\ F^1)=(0) .
$$
By [P] Theorem 4.5 and its proof, $M'_{r,s}$ is flat over $\Spec\O_E$
when $r=n-1$, $s=1$. Note that the scheme $M'_{r,s}$
has the same generic fiber as $M^{\rm naive}$.
Hence, by the above, if $r=n-1$, $s=1$,
$M^{\rm loc}=M'_{r,s}$ is flat over $\Spec\O_E$.

In fact, the calculations described in loc. cit., 4.16 suggest
that $M'_{r,s}$, and therefore also $M^{\rm loc}$, should be flat over $\Spec\O_E$ for all values of
$r$, $s$. The discussion in loc. cit., 4.16
shows that this flatness statement follows if one knows that the subscheme of $n\times n$-matrices
over $\Ff_p$ defined by
$$
\{A\in {\rm Mat}_{n\times n} |\ A^2=0, \ A=A^t, \ \wedge^{s+1}A=0, \ \wedge^{r+1}A=0, \ \det(T\cdot I-A)\equiv T^n\}
$$
is reduced. This can be viewed as the symmetric matrix version
of a result of Strickland [St] (compare to [PR] Cor. 5.10) and
it can be verified (for various primes $p$)
using Macaulay when $r,s\leq 5$.

In the case considered in this remark, the parahoric subgroup fixing the
lattice chain ${\cal L}$ is a special maximal parahoric. For more general
lattice chains one encounters additional problems.

 }\end{Remark}
  \bigskip
\bigskip

\section{Moduli spaces of abelian varieties}
\setcounter{equation}{0}

In this section we briefly indicate the construction of moduli spaces of
abelian varieties corresponding, in a sense made precise by the diagram
(\ref{15.5}) below, to the splitting and local models of the previous
section. The use of the language of algebraic stacks in (\ref{15.5})
replaces the method of {\it linear modifications} of [P]; its
mathematical content is the same. The reader can refer to [LMB] 
for background on the theory of algebraic stacks.

In this section we will use the following notation, taken from [RZ], ch.\
6. Let $B$ be a semi-simple algebra over $\Q$ and let $\ast$ be a positive
involution on $B$. Let $V$ be a finite-dimensional $\Q$-vector space with
a nondegenerate alternating bilinear form $(\ ,\ )$ with values in $\Q$.
We assume that $V$ is equipped with a $B$-module structure such that
$$(bv,w)= (v,b^\ast w),\ \ v,w\in V,\ \ b\in B\ \ .$$ Let $G\subset
GL_B(V)$ be the closed algebraic subgroup over $\Q$ such that
 $$G(\Q)=\{ g\in GL_B(V)\ \vert\
(gv,gw)=c(g)(v,w),\ c(g)\in\Q\}\ \ .$$ Let ${\cal S}=R_{\C/\R} \G_m$ and
let $h:{\cal S}\to G_\R$ be a homomorphism satisfying the usual Riemann
bilinear relations (cf.\ loc.cit.). We have a corresponding Hodge
decomposition $$V\otimes \C = V_0\oplus V_1$$ and a corresponding
cocharacter $\mu$ of $G$ defined over $\C$. We let $E\subset\overline\Q$
be the corresponding Shimura field. We now fix a prime number $p$ and
choose an embedding $\overline\Q\to\overline\Q_p$. The corresponding
$\nu$-adic completion of $E$ will be denoted $E_\nu$. Let $C^p\subset
G({\bold A}_f^p)$ be an open compact subgroup.

We consider an order $O_B$ of $B$ such that $O_B\otimes\Z_p$ is a maximal
order of $B\otimes \Q_p$. We assume that $O_B\otimes \Z_p$ is invariant
under the involution. We also fix a selfdual periodic multichain ${\cal
L}$ of $O_B\otimes \Z_p$-lattices in $V\otimes\Q_p$ with respect to the
alternating form $(\ ,\ )$.

We recall from loc.cit.\ the definition of a moduli problem ${\cal
A}_{C^p}$ over $(Sch/{\rm Spec}\ {\cal O}_{E_\nu})$. It associates to a
${\cal O}_{E_\nu}$-scheme $S$ the following set of data up to isomorphism:
\begin{itemize}
\item[1.] An ${\cal L}$-set of abelian varieties $A=\{ A_\Lambda\}$.
\item[2.] A $\Q$-homogeneous principal polarization $\overline\lambda$ of
the ${\cal L}$-set $A$.
\item[3.] A $C^p$-level structure
$$\overline\eta :H_1(A,{\bf A}_f^p)\simeq V\otimes {\bf A}_f^p\ {\rm mod}\
C^p\ \ ,$$ which respects the bilinear forms on both sides up to a
constant in $({\bf A}_f^p)^\times$.
\end{itemize}

We require an identity of characteristic polynomials, $$ {\rm{det}}
(T\cdot I-b\ \vert\  {\rm{Lie}}\ A_\Lambda)={\rm{det}} (T\cdot I-b\ \vert\
V_0),\ \ b\in O_B,\ \Lambda\in {\cal L}\ \ . $$ For the definitions of the
terms employed here we refer to loc.cit., 6.3--6.8. We only mention that
$A$ is a functor from the category ${\cal L}$ to the category of abelian
schemes over $S$ up to isogeny of order prime to $p$, with $O_B$-action,
and that a polarization $\lambda$ is a ${\cal O}_B$-linear homomorphism
from $A$ to the dual ${\cal L}$-set $\tilde A$ (for which $\tilde
A_\Lambda=(A_{\Lambda^*})^\wedge)$.

The functor ${\cal A}_{C^p}$ is representable by a quasi-projective scheme
over ${\cal O}_{E_\nu}$, provided that $C^p$ is sufficiently small.

We denote by $M_\Lambda$
the Lie algebra of the universal extension of $A_\Lambda$. Then $\{
M_\Lambda\}$ is a polarized multichain of $({\cal O}_B\otimes
\Z_p)\otimes_{\Z_p}{\cal O}_S$-modules on $S$ of type $({\cal L})$ in the
sense of [RZ], Def.\ 3.14. Let $\tilde{\cal A}_{C^p}$ be the functor which
to $S\in (Sch/{\cal O}_{E_\nu})$ associates the isomorphism classes of
objects $(A, \overline\lambda,\overline\eta)$ of ${\cal A}_{C^p}(S)$ and
an isomorphism of polarized multichains between $\{ M_\Lambda\}$ and
${\cal L}\otimes_{\Z_p}{\cal O}_S$. By [P], Thm.\ 2.2 (a slight extension of [RZ] Thm.\ 3.16), the forgetful
morphism

\begin{eq}\label{15.2}
\pi: \tilde{\cal A}_{C^p}\longrightarrow {\cal A}_{C^p}
\end{eq}
is a principal homogeneous space, locally trivial for the \'etale topology, under the smooth
group scheme ${\cal G}\times_{{\rm Spec}\ \Z_p}{\rm Spec}\ {\cal
O}_{E_\nu}$. Here ${\cal G}=\underline{{\rm Aut}}({\cal L})$ is the group
scheme over ${\rm Spec}\ \Z_p$ with $C_p={\cal G}(\Z_p)$ the  
subgroup of $G(\Q_p)$ fixing the lattice chain ${\cal L}$.

The Lie algebra ${\rm Lie}\ A_\Lambda$ is a factor module $t_\Lambda$ of
$M_\Lambda$. Using the identification of $M_\Lambda$ with
$\Lambda\otimes_{\Z_p}{\cal O}_S$ over $\tilde{\cal A}_{C^p}$ we therefore
obtain a point of the naive local model $M^{\rm naive}$ defined in terms
of the $\Z_p$-data $(B\otimes\Q_p, O_B\otimes\Z_p, V\otimes\Q_p, {\cal
L})$ induced from our global data,

\begin{eq}
\tilde\varphi:\tilde{\cal A}_{C^p}\longrightarrow M^{\rm naive}\ \ .
\end{eq}
Since $\tilde \varphi$ is obviously equivariant for the action of ${\cal
G}\otimes_{\Z_p}{\cal O}_{E_\nu}$, $\tilde\varphi$ corresponds to a
relatively representable morphism of algebraic stacks
\begin{eq}
\varphi: {\cal A}_{C^p}\longrightarrow \left[ M^{\rm naive} /{\cal
G}\otimes_{\Z_p}{\cal O}_{E_\nu}\right]\ \ .
\end{eq}
By [P], Thm.\ 2.2, (a slight extension of [RZ], Prop.\ 3.3), the morphism $\varphi$ is smooth of relative
dimension ${\rm dim}\ G$. Let us form the cartesian product of $\varphi$
with the morphisms ${\cal M}\to M^{\rm loc}\hookrightarrow M^{\rm naive}$,
where ${\cal M}$ denotes the splitting model over ${\cal O}_K$, with $K$
the Galois closure of $E_\nu$,
\begin{eqnarray}\label{15.5}
\begin{array}{ccc}
{\cal A}_{C^p}^{\rm spl} & \longrightarrow & [{\cal M}/{\cal G}_{{\cal
O}_K}]
\\
\\
\big\downarrow && \big\downarrow
\\
\\
{\cal A}_{C^p}^{\rm loc} & \longrightarrow & [M^{\rm loc}/{\cal G}_{{\cal
O}_{E_\nu}}]
\\
\\
\big\downarrow && \big\downarrow
\\
\\
{\cal A}_{C^p} & \longrightarrow & [ M^{\rm naive}/{\cal G}_{{\cal
O}_{E_\nu}}].
\end{array}
\end{eqnarray}
The scheme ${\cal A}_{C^p}^{\rm loc}$ is a closed subscheme of ${\cal
A}_{C^p}$ and is the image of ${\cal A}_{C^p}^{\rm spl}$ in ${\cal
A}_{C^p}$. The scheme ${\cal A}_{C^p}^{\rm loc}$ is a linear modification
of ${\cal A}_{C^p}$ in the sense of [P]; likewise, ${\cal A}_{C^p}^{\rm
spl}$ is a linear modification of ${\cal A}_{C^p}\otimes_{{\cal O}_{E_\nu}}{\cal O}_K$.

The ${\cal O}_K$-scheme ${\cal A}_{C^p}^{\rm spl}$ represents the
following moduli problem on $(Sch/{\cal O}_K)$. It associates to $S$ the
set of isomorphism classes of objects $(A, \overline\lambda,
\overline\eta, {\cal F})$. Here $(A=\{ A_\Lambda\},
\overline\lambda,\overline\eta)$ is an object of ${\cal A}_{C^p}(S)$.

Let $F_1$ be the invariants under $*$ in the center $F$ of $B\otimes
\Q_p$. Then $F_1$ is a direct sum of fields,

\begin{eq}\label{15.6}
F_1=F_{1,1}\oplus\ldots\oplus F_{1,r}\ \ .
\end{eq}
For $k=1,\ldots, r$, let $n_k$ be the degree over $\Q_p$ of the direct summand of $F$
corresponding to the direct summand $F_{1,k}$ of $F_1$. Let $M_\Lambda$ be
the Lie algebra of the universal extension of $A_\Lambda$ and let
$F_\Lambda$ be the kernel of the factor map from $M_\Lambda$ to ${\rm
Lie}\ A_\Lambda$. Then the action of ${\cal O}_{F_1}$ on $M_\Lambda$ and
$F_\Lambda$ induces decompositions $$M_\Lambda = \bigoplus_{k=1}^r
M_{\Lambda, k}\ \ ,\ \ F_\Lambda =\bigoplus_{k=1}^r F_{\Lambda, k}\ \ .$$
The final ingredient ${\cal F}$ of an object of ${\cal A}_{C^p}^{\rm
spl}(S)$ is a collection of functors $\Lambda\mapsto F^i_{\Lambda, k}$ for
$k=1,\ldots, r$ and $i=1,\ldots, n_k$, with functor morphisms
$j^i_{\Lambda, k}: F^i_{\Lambda, k}\to M_{\Lambda, k}$, satisfying for
each $k=1,\ldots, r$ the conditions in Definition \ref{gensplit} when
$\Lambda\otimes {\cal O}_S$ is replaced by $M_{\Lambda, k}$ and
$(F^i_\Lambda, j^i_{\Lambda})$ by $(F^i_{\Lambda, k}, j^i_{\Lambda, k})$,
and such that $F^1_{\Lambda, k}=F_{\Lambda, k}$.

On the other hand, it seems that one cannot hope in general to be able
to describe a ``simple and explicit" moduli problem over $\O_E$ 
that is represented by ${\cal A}^{\rm loc}_{C^p}$.
This is of course a question of finding the appropriate 
conditions on $F_\Lambda={\rm ker}(M_\Lambda\to {\rm Lie}(A_\Lambda))$
that would cut out the closed subscheme ${\cal A}^{\rm loc}_{C^p}\subset
{\cal A}_{C^p}$ (see [PR] Theorem 5.7 for an example in which 
such explicit -but quite complicated- conditions are proposed).

\medskip
\smallskip
\bigskip\bigskip\bigskip

\leftline{Georgios Pappas\hfill Michael Rapoport} \leftline{Dept. of
Mathematics\hfill Mathematisches Institut} \leftline{Michigan State
University\hfill der Universit\"at Bonn} \leftline{E. Lansing\hfill
Beringstr. 1} \leftline{MI 48824-1027\hfill 53115 Bonn}
\leftline{USA\hfill Germany}\leftline{email: pappas@math.msu.edu\hfill
email: rapoport@math.uni-bonn.de}

\end{document}